\documentclass[aop,preprint]{imsart}

\RequirePackage{amsthm,amsmath,amsfonts,amssymb}
\RequirePackage[numbers]{natbib}
\usepackage{hyperref,accents}
\usepackage{tikz}
\usetikzlibrary{calc,arrows,positioning,fit,petri}

\usepackage[normalem]{ulem}%for \sout,\xout
\usepackage{cancel}%for \cancel, \bcancel, \xcancel, \cancelto

\startlocaldefs

\definecolor{lightyellow}{rgb}{1,1,0.7}
\definecolor{greeny}{rgb}{0,.6,0}
\definecolor{orange}{rgb}{1.,0.3,0.1}
\newcommand{\detail}[1]{\par\noi{\bf [Proof detail\ }{#1}
	\hfill{\bf ]}\par\noi\hspace{-4pt}}
\renewcommand{\detail}[1]{}

\newcommand{\dis}{\displaystyle}

\newcommand{\scr}{\scriptstyle}

\newcommand{\noi}{\noindent}
\newcommand{\halmos}{\rule{1ex}{1.4ex}}
\newcommand{\QED}{\nopagebreak{\hspace*{\fill}$\halmos$\medskip}}

\newcommand{\quand}{\quad\mbox{and}\quad}

\newtheoremstyle{mythm}% name
{}%      Space above
{}%      Space below
{\itshape}%         Body font
{}%         Indent amount (empty = no indent, \parindent = para indent)
{\bfseries}% Thm head font
{}%        Punctuation after thm head
{.5em}%     Space after thm head: " " = normal interword space;
{#1 #2 \thmnote{(#3)}}%  Thm head spec (can be left empty, meaning `normal')

\theoremstyle{mythm}
\newtheorem{theorem}{Theorem}[section]
\newtheorem{proposition}[theorem]{Proposition}
\newtheorem{lemma}[theorem]{Lemma}
\newtheorem{exercise}[theorem]{Exercise}
\newtheorem{corollary}[theorem]{Corollary}
\newtheorem{conjecture}[theorem]{Conjecture}

\newtheorem{counterex}[theorem]{Counterexample}
\newtheorem{remark}[theorem]{Remark}

\newcommand{\bt}{\begin{theorem}}
	\newcommand{\et}{\end{theorem}}
\newcommand{\bl}{\begin{lemma}}
	\newcommand{\el}{\end{lemma}}
\newcommand{\bp}{\begin{proposition}}
	\newcommand{\ep}{\end{proposition}}
\newcommand{\bcor}{\begin{corollary}}
	\newcommand{\ecor}{\end{corollary}}
\newcommand{\br}{\begin{remark}\rm}
	\newcommand{\er}{\end{remark}}
\newcommand{\bcon}{\begin{conjecture}}
	\newcommand{\econ}{\end{conjecture}}
\newcommand{\bex}{\begin{exercise}}
	\newcommand{\eex}{\end{exercise}}
\newcommand{\bcou}{\begin{counterex}}
	\newcommand{\ecou}{\end{counterex}}

\newenvironment{Proof}[1][]{\noi\textbf{Proof #1}}{\QED}
\newcommand{\bpro}{\begin{Proof}}
	\newcommand{\epro}{\end{Proof}}

\newcommand{\be}{\begin{equation}}
	\newcommand{\ee}{\end{equation}}
\newcommand{\ba}{\begin{array}}
	\newcommand{\ea}{\end{array}}
\newcommand{\bc}{\be\begin{array}{r@{\,}c@{\,}l}}
	\newcommand{\ec}{\end{array}\ee}

\newcommand{\ga}{\gamma}
\newcommand{\Ga}{\Gamma}
\newcommand{\de}{\delta}
\newcommand{\De}{\Delta}
\newcommand{\eps}{\varepsilon}

\newcommand{\sig}{\sigma}

\newcommand{\tet}{\theta}
\newcommand{\Tet}{\Theta}
\newcommand{\si}{\ensuremath{\sigma}}

\newcommand{\Ai}{{\cal A}}
\newcommand{\Bi}{{\cal B}}

\newcommand{\Di}{{\cal D}}

\newcommand{\Fi}{{\cal F}}
\newcommand{\Gi}{{\cal G}}
\newcommand{\Hi}{{\cal H}}

\newcommand{\Ni}{{\cal N}}

\newcommand{\Pc}{{\cal P}}

\newcommand{\Ri}{{\cal R}}

\newcommand{\Wi}{{\cal W}}
\newcommand{\Xc}{{\cal X}}

\newcommand{\R}{{\mathbb R}}
\newcommand{\N}{{\mathbb N}}
\newcommand{\Z}{{\mathbb Z}}

\newcommand{\Q}{{\mathbb Q}}

\renewcommand{\P}{{\mathbb P}}
\newcommand{\E}{{\mathbb E}}

\newcommand{\li}{\langle}
\newcommand{\re}{\rangle}

\newcommand{\down}{\downarrow}
\newcommand{\sub}{\subset}

\newcommand{\asto}[1]{\underset{{#1}\to\infty}{\longrightarrow}}

\newcommand{\Astoo}[1]{\underset{{#1}\to 0}{\Longrightarrow}}

\newcommand{\ti}{\tilde}

\newcommand{\ov}{\overline}

\newcommand{\di}{\mathrm{d}}
\newcommand{\half}{{[0,\infty)}}

\newcommand{\wt}[1]{\widetilde{#1}}
\newcommand{\wN}{\widetilde{N}}

\newcommand{\Rc}{R^2_{\rm c}}

\newcommand{\aro}[2]{\filldraw (#1) circle (0.1cm);
 \draw[very thick,->] (#1)
  let \p1=($(#2,1)$) in +($0.1*(\x1,\y1)$)--+($0.9*(\x1,\y1)$);}

\newcommand{\yu}[1]{#1}

\makeatletter\@addtoreset{equation}{section}
\makeatother

\endlocaldefs

\begin{document}

\begin{frontmatter}
%%%%%%%%%%%%%%%%%%%%%%%%%%%%%%%%%%%%%%%%%%%%%%
%%                                          %%
%% Enter the title of your article here     %%
%%                                          %%
%%%%%%%%%%%%%%%%%%%%%%%%%%%%%%%%%%%%%%%%%%%%%%
\title{Universality of the Brownian net}
%\title{A sample article title with some additional note\thanksref{T1}}
\runtitle{Universality of the Brownian net}
%\thankstext{T1}{A sample of additional note to the title.}

\begin{aug}
%%%%%%%%%%%%%%%%%%%%%%%%%%%%%%%%%%%%%%%%%%%%%%%
%% Only one address is permitted per author. %%
%% Only division, organization and e-mail is %%
%% included in the address.                  %%
%% Additional information can be included in %%
%% the Acknowledgments section if necessary. %%
%% ORCID can be inserted by command:         %%
%% \orcid{0000-0000-0000-0000}               %%
%%%%%%%%%%%%%%%%%%%%%%%%%%%%%%%%%%%%%%%%%%%%%%%
\author[A]{\fnms{Rongfeng}~\snm{Sun}\ead[label=e1]{matsr@nus.edu.sg}},
\author[B]{\fnms{Jan~M.}~\snm{Swart}\ead[label=e2]{swart@utia.cas.cz}}
\and
\author[C]{\fnms{Jinjiong}~\snm{Yu}\ead[label=e3]{jjyu@sfs.ecnu.edu.cn}}
%%%%%%%%%%%%%%%%%%%%%%%%%%%%%%%%%%%%%%%%%%%%%%
%% Addresses                                %%
%%%%%%%%%%%%%%%%%%%%%%%%%%%%%%%%%%%%%%%%%%%%%%
\address[A]{Department of Mathematics,
	National University of Singapore,
	10 Lower Kent Ridge Road, 119076 Singapore.
	\printead[presep={,\ }]{e1}}

\address[B]{The Czech Academy of Sciences,
	Institute of Information Theory and Automation,
	Pod vod\'arenskou v\v e\v z\' i 4,
	18208 Praha 8,
	Czech Republic.
	\printead[presep={,\ }]{e2}}

\address[C]{KLATASDS-MOE, School of Statistics,
	East China Normal University, 3663 North Zhongshan Road,
	Shanghai 200062, China.
	\printead[presep={,\ }]{e3}}
\end{aug}

\begin{abstract}
The Brownian web is a collection of one-dimensional coalescing Brownian motions starting from every point in space and time, while the Brownian net is an extension that also allows branching. We show here that the Brownian net is the universal scaling limit of one-dimensional branching-coalescing random walks with weak binary branching and arbitrary increment distributions that have finite $(3+\eps)$-th moment. This gives the first example in the domain of attraction of the Brownian net where paths can cross without coalescing, which poses fundamental technical challenges not present in the non-crossing case.
\end{abstract}

\begin{keyword}[class=MSC]
\kwd[Primary ]{82C22}
%\kwd{???}
\kwd[; secondary ]{82C41, 60K35}
\end{keyword}

\begin{keyword}
\kwd{Branching-coalescing random walks}
\kwd{Brownian web}
\kwd{Brownian net}
\kwd{invariance principle}
\end{keyword}

\end{frontmatter}
%%%%%%%%%%%%%%%%%%%%%%%%%%%%%%%%%%%%%%%%%%%%%%
%% Please use \tableofcontents for articles %%
%% with 50 pages and more                   %%
%%%%%%%%%%%%%%%%%%%%%%%%%%%%%%%%%%%%%%%%%%%%%%
\tableofcontents

%%%%%%%%%%%%%%%%%%%%%%%%%%%%%%%%%%%%%%%%%%%%%%
%%%% Main text entry area:

\section{Introduction and main result}\label{S:intro}

The Brownian web $\Wi$ is a collection of coalescing Brownian motions on $\R$, starting from every point
in the space-time plane $\R^2$. It originated in Arratia's Ph.D.\ thesis \cite{Arr79}, where he constructed
coalescing Brownian motions starting from every point in $\R$ at time $0$. In a subsequent unpublished manuscript
\cite{Arr81}, his attempt to construct coalescing Brownian motions starting from every point in the space-time plane $\R^2$ was
incomplete due to the presence of special points in $\R^2$ where multiple Brownian paths can emanate. The construction
was completed by T\'oth and Werner in \cite{TW98}, where a variant of Arratia's coalescing Brownian motions played a crucial
role in the study of the so-called true self-repelling motion. They established detailed properties for Arratia's
coalescing Brownian motions, including a complete classification of its special points. Subsequently, Fontes, Isopi, Newman and
Ravishankar \cite{FINR02,FINR04} came across the same object in the study of aging in the zero-temperature dynamics of Ising and Potts
models on $\Z$. They introduced a topology for Arratia's coalescing Brownian motions so that it becomes a random variable taking
values in a Polish space, and they named it the {\em Brownian web} $\Wi$. They also provided a characterisation for the Brownian web,
which they used to prove the weak convergence of systems of coalescing simple random walks to the Brownian web.

In \cite{SS08}, Sun and Swart made an important generalisation of the Brownian web $\Wi$ by allowing paths to branch, which they
call the {\em Brownian net} $\Ni$. To counter the instantaneous coalescence between paths, the branching rate is effectively ``infinite''.
Independently, a different construction of the Brownian net was given by  Newman, Ravishankar, and Schertzer in \cite{NRS10}.

Since their inception, the Brownian web and the Brownian net have found connections to many other topics of research, including population
genetic models \cite{AS11, GSW16, EFS17}, drainage network type models \cite{CDF09, CV14, RSS16, PPR22, RSS23, VZ23, SVZ22},  Poisson webs
\cite{FFW05,FVV15}, the directed spanning forest \cite{CSST21}, oriented percolation \cite{SS11, BGS19, SS19},  planar aggregation models \cite{NT12,NT15}, random maps \cite{KL21}, random walks in i.i.d.\ space-time random environments \cite{SSS14}, traffic models \cite{A17}, determinantal and Pfaffian point processes \cite{TZ11, GPTZ18, GTZ20}. We refer to the survey  \cite{SSS17}  for a more comprehensive list of earlier results. More recently, there have been interesting connections with interface growth models in the Edwards-Wilkinson and KPZ universality classes \cite{Y16, VV26, DDP24}, and a new class constructed from the Brownian web, known as the Brownian castle \cite{CH23a, CH23b}, which admits a further extension that can be constructed from the Brownian net \cite{CHS25+}. There are also intriguing parallels between the special points in the Brownian web and net \cite{SSS09} and those in the {\em directed landscape} \cite{D23}. This makes the Brownian web and Brownian net and their universality classes important objects of study in their own right.

The Brownian web $\Wi$ and the Brownian net $\Ni$ are expected to be the universal scaling limits of one dimensional coalescing
systems, and branching-coalescing systems with weak branching, respectively. The universality of the Brownian web $\Wi$ has been studied
extensively. Models which have been shown to converge to the Brownian web fall under two classes: either paths must coalesce whenever they attempt to cross each other, called the {\em non-crossing} case; or paths can cross each other without coalescing, called the {\em crossing case}. For the non-crossing case, convergence criteria were formulated in \cite{FINR04} and has been verified for many models, including coalescing simple random walks \cite{FINR04}, extremal paths in supercritical oriented percolation \cite{SS11,SS19}, paths in drainage network models \cite{CDF09,RSS16}, Poisson webs \cite{FFW05,FVV15}, the directed spanning forest \cite{CSST21}, etc. The crossing case is much more challenging due to the lack of ordering among paths.The first example in the crossing case that was shown to converge to the Brownian web is the collection of coalescing random walks on $\Z$ with general increment distributions, provided the increments have finite $(3+\eta)$-th moment for some $\eta>0$ \cite{NRS05, BMSV06}. In \cite{NRS05}, the authors formulated a set of convergence criteria tailored for the crossing case, which has also been verified for some drainage network models in \cite{CV14,VZ23}.

Universality results for the Brownian net $\Ni$ are much more limited. For branching-coalescing simple symmetric random walks on $\Z$ with weak branching, i.e., the branching probability is proportional to the space-time diffusive scaling parameter $\eps$, convergence to the Brownian net $\Ni$ was established in \cite{SS08}, which relies heavily on the {\em non-crossing} property that paths must meet before crossing.
In particular, the set of leftmost (resp.\ rightmost) paths starting from each point in the space-time lattice forms a non-crossing system of coalescing random walk paths. The only other model (also in the non-crossing case) that has been shown to converge to the Brownian net is the set of genealogical paths of a spatial $\Lambda$-Fleming-Viot process considered in \cite{EFS17}. There has also been partial progress for a drainage network model with weak branching \cite{SVZ22}, which also belongs to the non-crossing case.

\begin{figure}
	\begin{center}
		
\begin{tikzpicture}[xscale=0.8,yscale=0.7,>=stealth]
	
	% Slightly enlarge the bounding box to accommodate the direction indicator.
	\path[use as bounding box] (3.5,-0.8) rectangle (15.8,8.125);
	
	\begin{scope}
		\clip (3.5,-0.125) rectangle (15.5,8);
		
		%\aro{3,8}{3}\aro{4,8}{-1}\aro{5,8}{-1}\aro{6,8}{0}\aro{7,8}{-1}\aro{8,8}{-3}
		%\aro{9,8}{-3}\aro{10,8}{0}\aro{11,8}{-1}
		\aro{3,7}{2}\aro{4,7}{-1}\aro{5,7}{1}
		\aro{6,7}{1}\aro{7,7}{-1}\aro{8,7}{-3}\aro{9,7}{2}\aro{10,7}{-2}\aro{11,7}{1}
		\aro{3,6}{2}\aro{4,6}{2}\aro{5,6}{1}\aro{6,6}{0}\aro{7,6}{-2}\aro{8,6}{-1}
		\aro{9,6}{2}\aro{10,6}{2}\aro{11,6}{3}\aro{4,5}{1}\aro{5,5}{0}\aro{6,5}{2}
		\aro{7,5}{1}\aro{8,5}{3}\aro{9,5}{-3}\aro{10,5}{-2}\aro{11,5}{-3}\aro{3,4}{1}
		\aro{4,4}{-3}\aro{5,4}{2}\aro{6,4}{1}\aro{7,4}{0}\aro{8,4}{1}\aro{9,4}{1}
		\aro{10,4}{-3}\aro{11,4}{2}\aro{3,3}{1}\aro{4,3}{2}\aro{5,3}{1}\aro{6,3}{1}
		\aro{7,3}{-1}\aro{8,3}{1}\aro{9,3}{3}\aro{10,3}{1}\aro{11,3}{-1}\aro{4,2}{-2}
		\aro{5,2}{-1}\aro{6,2}{-2}\aro{7,2}{-1}\aro{8,2}{2}\aro{9,2}{1}\aro{10,2}{0}
		\aro{11,2}{-2}\aro{4,1}{2}\aro{5,1}{-3}\aro{6,1}{-1}\aro{7,1}{3}\aro{8,1}{-1}
		\aro{9,1}{-2}\aro{10,1}{0}\aro{11,1}{3}\aro{4,0}{3}\aro{5,0}{2}\aro{6,0}{3}
		\aro{7,0}{2}\aro{8,0}{0}\aro{9,0}{-2}\aro{10,0}{0}\aro{11,0}{0}
		%\aro{12,8}{-2}\aro{13,8}{-1}\aro{14,8}{-2}\aro{15,8}{-1}
		\aro{12,7}{-2}\aro{13,7}{-1}\aro{14,7}{-1}
		\aro{15,7}{0}\aro{16,7}{-2}\aro{17,7}{-2}\aro{12,6}{-2}\aro{13,6}{1}\aro{14,6}{1}
		\aro{15,6}{-1}\aro{12,5}{-1}\aro{13,5}{0}\aro{14,5}{-1}\aro{15,5}{-1}\aro{16,5}{-1}
		\aro{17,5}{-2}\aro{12,4}{1}\aro{13,4}{-2}\aro{14,4}{1}\aro{15,4}{-1}\aro{12,3}{1}
		\aro{13,3}{2}\aro{14,3}{0}\aro{15,3}{1}\aro{16,3}{-1}\aro{12,2}{2}\aro{13,2}{1}
		\aro{14,2}{1}\aro{15,2}{0}\aro{12,1}{-2}\aro{13,1}{-2}\aro{14,1}{0}\aro{15,1}{2}
		\aro{16,1}{-1}\aro{12,0}{-2}\aro{13,0}{1}\aro{14,0}{-2}\aro{15,0}{1}\aro{16,0}{-1}
		\aro{11,4}{0}\aro{10,6}{0}\aro{13,6}{-2}\aro{9,3}{2}\aro{9,3}{2}\aro{11,4}{0}
		\aro{11,6}{-2}\aro{4,7}{1}\aro{8,5}{1}\aro{15,7}{-2}\aro{4,4}{1}
		%\aro{6,8}{2}\aro{7,8}{2}\aro{14,8}{0}
		\aro{7,1}{2}\aro{11,7}{1}\aro{13,7}{-2}\aro{16,5}{-1}
		\aro{3,4}{1}
		
	\end{scope}
	
	% Coordinate directions
	\draw[->,very thick] (2.88,-0.8) -- (3.7,-0.8)
	node[right] {\small Space};
	\draw[->,very thick] (2.9,-0.82) -- (2.9,0)
	node[above] {\small Time};
	
\end{tikzpicture}

		\caption[A discrete net.]{{\yu A discrete net with crossing paths.}}
		\label{fig:disnet}
	\end{center}
\end{figure}

An important open problem is to find models with crossing paths {(see Figure \ref{fig:disnet})} that converge to the Brownian net. It is long believed that one such model is branching-coalescing random walks on $\Z$ with weak branching, whose increment distribution satisfies the same assumptions as in \cite{NRS05, BMSV06}. In this paper, we verify this conjecture and thus provide the first example with crossing paths that lies in the domain of attraction of the Brownian net $\Ni$. The fact that such a result is only proved now is due to fundamental difficulties that arise from branching, together with the lack of ordering
among paths. The ideas we develop here should also be useful for proving convergence to the Brownian net for other models in the crossing case.
%More precisely, we consider the system of {\em branching-coalescing random walks} on $\Z$ with weak binary branching and increments that have finite $(3+\eta)$-th moment for some $\eta>0$.

\subsection{Model and Main Result}\label{S:disc}

We now introduce the collection of branching-coalescing random walk paths on $\Z$, which are discrete analogues of the Brownian net, and hence called {\em discrete nets}.

Let $a(\cdot)$ be a probability kernel on $\Z$, which specifies the increment distribution for a random walk on $\Z$.
We assume that
\be\ba{rl}\label{passum}
{\rm (i)} & a(\cdot) \mbox{ defines an irreducible and aperiodic random walk on  } \Z, \\ [5pt]
{\rm(ii)}&\dis\sum_{x\in\Z}xa(x)=0,\\ [5pt]
%{{\rm(iii)}}&{\yu \dis\sig^2:=\sum_{x\in\Z}x^2a(x) <\infty ,}\\
{\rm(iii)}&\dis\sum_{x\in\Z}|x|^{3+\eta}a(x)<\infty\mbox{ for some }\eta>0.
\ec
We will denote ${\varsigma^2:=\sum_{x\in\Z}x^2a(x)}$ throughout the rest of the paper.

Before introducing branching, we first give a graphical construction of coalescing random walks starting from every point of the
space-time lattice $\Z^2=\{(x,t): x, t\in \Z\}$, which we call a {\em discrete web}. Let $\omega=(\omega(x,t))_{(x,t)\in \Z^2}$ be  i.i.d.\ random variables with common distribution $a(\cdot)$. From each $(x,t)\in \Z^2$, we then draw an arrow from $(x,t)$ to $(x+\omega(x,t), t+1)$. Note that for each
$(x, t)\in\Z^2$, there is a unique path $\pi$ starting from $(x,t)$ by following the arrows, and two paths coalesce when they meet at the space-time
lattice point. We will denote by $\sigma_\pi :=t$ the starting time of the path $\pi$, and identify $\pi: [\sigma_\pi, \infty) \to \R$ with a continuous function by linearly interpolating between consecutive integer times. Given $\omega$, let $\pi_{(x,t)}$ denote the random walk path starting from $(x,t)\in \Z^2$. We then call the random set of paths
\be\label{W1}
W:=\big\{\pi_{(x,t)}:(x,t)\in\Z^2\big\}
\ee
a \emph{discrete web with jump distribution $a(\cdot)$}.

To allow for binary branching, we need to allow the possibility of two arrows coming out of each $(x,t)\in \Z^2$. Let
$(\omega^1(x,t), \omega^2(x,t))_{(x,t)\in \Z^2}$ be an i.i.d.\ family of $\Z\times\Z$-valued random variables with common law $a^{(2)}(\cdot, \cdot)$, to be specified in \eqref{a2}. Note that each of the two families $\omega^1=(\omega^1(x,t))_{(x,t)\in \Z^2}$ and $\omega^2=(\omega^2(x,t))_{(x,t)\in \Z^2}$ induces a set of arrows and a corresponding discrete web, which we denote by $W^1$ and $W^2$ respectively.
The paths between $W^1$ and $W^2$ undergo sticky interaction, where stickiness comes from $(x,t)$ with $\omega^1(x,t)= \omega^2(x,t)$. If we consider all paths that can follow both types of arrows, induced by either $\omega^1$ or $\omega^2$, then paths still coalesce when they meet at the same space-time in $\Z^2$, but a path can also branch into two at $(x,t)$ with $\omega^1(x,t)\neq \omega^2(x,t)$. We call the set of all paths obtained by following arrows induced by either $\omega^1$ or $\omega^2$ a {\em discrete net} {(see Figure~\ref{fig:disnet})}.

From now on, we assume that the law $a^{(2)}$ is of the form
\begin{equation}\label{a2}
	a^{(2)}(x_1,x_2):=(1-\eps)1_{\{x_1=x_2\}}a(x_1)+\eps a(x_1)a(x_2),
\end{equation}
which means that with probability $1-\eps$, a single arrow comes out of $(x,t)$ with law $a(\cdot)$; and with probability $\eps$, two independent arrows with common law $a(\cdot)$ are drawn from $(x,t)$, which coalesce at time $t+1$ if they coincide. Let us denote the resulting discrete net by $N_\eps$, which we call a {\em discrete net with (binary) branching probability $\eps$ and jump distribution $a(\cdot)$}. The branching is weak in the sense that we will consider a sequence of such models with $\eps\downarrow 0$.

Our main result is to show that under diffusive scaling of space-time (recall $\varsigma$ from \eqref{passum})
\be\label{Seps}
{S_\eps} : (x, t) \to (\eps x, \varsigma^2\eps^2 t)
\ee
with the same $\eps$ as the branching probability in $N_\eps$, the discrete net $N_\eps$ converges to a standard Brownian net ${\cal N}$.
With a slight overload of notation, for each path $\pi: [\sigma_\pi, \infty)\to\R$, we will use $S_\eps \pi$ to denote the path whose graph is the image of the graph of $\pi$ under $S_\eps$, and for each set of paths $K$, we will use $S_\eps K$ to denote the set of paths obtained by applying $S_\eps$ to each path in $K$. In Section~\ref{S:topo}, we will recall the topology introduced in \cite{FINR04} for the Brownian web $\Wi$. In particular, $N_\eps$
and $\Ni$ are both random subsets of the path space $(\Pi, d)$. Furthermore, $\Ni$ and $\overline{N}_\eps$, the closure of $N_\eps$ in $(\Pi, d)$
(see Remark \ref{R:Neps}), are both compact so that they are random variables taking values in the complete separable metric space $(\Hi, d_\Hi)$ of compact subsets of $(\Pi, d)$. To keep the notation simple, we will often write $N_\eps$ instead of $\overline{N}_\eps$. Our main convergence result can then be formulated as follows.

\bt[Convergence to the Brownian net]\label{T:net}
For $\eps>0$, let $N_\eps$ be the discrete net with (binary) branching probability $\eps$ and jump distribution $a(\cdot)$
satisfying \eqref{passum}. Then
\be\label{nonsimpnet}
S_\eps N_\eps \Astoo{\eps} \Ni\,,
\ee
where $\Ni$ is the standard Brownian net, and $\Rightarrow$ denotes weak convergence of $\Hi$-valued random variables.
\et

\begin{remark}
{In Theorem \ref{T:net}, if we replace the scaling map $S_\eps$ with $\widetilde S_\eps:  (x, t) \to (\eps x/\varsigma, \eps^2 t)$, then $\widetilde S_\eps N_\eps \Rightarrow \Ni_\varsigma$, which is the image of $\Ni$ under the map $(x, t) \to (x/\varsigma, t/\varsigma^2)$ and is called a Brownian net with branching parameter $\varsigma$ (see Theorem \ref{T:Nhop}). Similarly, if the discrete net in Theorem \ref{T:net} has branching probability $b\eps$ for some $b>0$, then $S_\eps N_{b\eps}\Rightarrow \Ni_b$.
}
\end{remark}

\begin{remark}
{Theorem \ref{T:net} implies the first part of \cite[Theorem 1.2]{BMSV06}, which shows that the rescaled discrete web $S_\eps W$ converges in distribution to the  Brownian web $\mathcal W$ under the same assumptions in \eqref{passum}. This follows from the observation that, for every $b>0$ and $\eps>0$, $W$ can be chosen as a subset of $N_{b\eps}$, the discrete net with branching probability $b\eps$, and hence as $\eps\downarrow 0$, every subsequential weak limit of $S_\eps W$ is stochastically dominated by a Brownian net $\Ni_b$ for every $b>0$. Since $\Ni_b$ stochastically decreases to the Brownian web $\Wi$ as $b\downarrow 0$,  while every subsequential weak limit of $S_\eps W$ contains a Brownian web, it follows that all subsequential weak limits of $S_\eps W$ are distributed as the Brownian web $\Wi$.
%Note that a discrete web $W$ is obtained if in (\ref{a2} we focus only on the arrows of the first type $x_1$. Under the finite $(3+\eta)$-th moment assumption on $a(\cdot)$ in \eqref{passum}), it was shown in the first half of \cite[Theorem 1.2]{BMSV06} that the rescaled discrete web $S_{\eps}W$ converges to the Brownian web $\Wi$. This result can also be implied from Theorem~\ref{T:net}.
%Indeed, if we rescale the branching parameter in (\ref{a2}) as $c_\eps\eps$ where $c_\eps=o(\eps)$ as $\eps\to0$, not only the same as the diffusive scale $\eps$, then the limit is the Brownian net with zero left and right drifts, which is the Brownian web $\Wi$ \cite{SS08}.

On the other hand,  \cite[Theorem 1.2]{BMSV06} also shows that if $\sum_{x\in \Z} \frac{|x|^3}{\log^\beta (|x|+2)} a(x)=\infty$ for some $\beta>1$, then the laws of the rescaled discrete web $S_\eps W$ fails to be tight due to the presence of arbitrarily large jumps originating from any finite space-time domain. Since $N_\eps$ contains more paths than $W$, the laws of $(S_{\eps}N_\eps)_{\eps \in (0,1)}$ also fail to be tight for such $a(\cdot)$.
	}

%In Theorem~\ref{T:net}, the finite $(3+\eta)$-th moment assumption on $a(\cdot)$ in \eqref{passum}  is almost optimal. It was shown in \cite[Theorem 1.1~$($ii$)$]{BMSV06} that if $\sum_{x\in \Z} \frac{|x|^3}{\log^\beta (|x|+2)} a(x)=\infty$ for some $\beta>1$, then the laws of the rescaled discrete web $(S_\eps W)_{\eps\in (0,1)}$ fails to be tight due to the presence of arbitrarily large jumps originating from any finite space-time domain. Therefore the laws of $(S_\eps N_\eps)_{\eps \in (0,1)}$ also fail to be tight for such $a(\cdot)$.}
\end{remark}

{
\begin{remark}
We can interpret the discrete net $N_\eps$ as a space-time random environment, such that when we sample a random walk in such a random environment, it follows paths in $N_\eps$ and chooses the two outgoing arrows with equal probability whenever it encounters a binary branching point. The transition kernels of such a random walk in random environment is known as a {\em stochastic flow of kernels} \cite{LR04}. For the discrete net $N_\eps$ with nearest neighbour jumps, it was shown in \cite[Section 1.3]{SSS14} that the associated stochastic flow of kernels converge in the diffusive scaling limit to the so-called Howitt-Warren flow with characteristic measure $\nu({\rm d}x)= c\, \delta_{1/2}({\rm d}x)$, and the joint law of $n$ random walks sampled independently in the same random environment converge to Howitt-Warren sticky Brownian motions. We expect the same to hold in our case. In analogy with the approach developed in \cite{SSS14}, the key step would be to strengthen Theorem \ref{T:net} and show that the set of so-called relevant branching points of $N_\eps$ (see Section \ref{S:RBP}) converge to the set of relevant separation points of~$\Ni$.
\end{remark}
}

\subsection{Proof strategy}\label{S:outline}

We first recall from \cite{SS08} the proof strategy for branching-coalescing simple random walks on $\Z$, which is a non-crossing system where paths must
meet before crossing (see also the survey article \cite{SSS17}). We will then explain the fundamental difficulties posed by crossing systems and how our proof strategy overcomes these difficulties for the system of branching-coalescing random walks we consider.

\subsubsection{{Previous results: the non-crossing Case}}
In \cite{SS08}, the authors considered the discrete net $N_\eps$ on the even space-time lattice $\Z_{\rm even}^2:=\{(x, t)\in \Z^2: x+t \mbox{ is even}\}$, where the jump distribution is given by $a(\pm 1)=1/2$ and each $(x,t)\in \Z^2_{\rm even}$ is a branching point with probability $\eps$.
The key observation behind the construction of the Brownian net $\Ni$ in \cite{SS08} is that, if we consider the discrete webs $W_\eps^{\rm l}, W_\eps^{\rm r}\subset N_\eps$, which consist of the leftmost (reps.\ rightmost) paths in $N_\eps$ starting from every $(x,t)\in \Z^2_{\rm even}$,
then under the diffusive scaling map $S_\eps$ with $\varsigma=1$, $S_\eps W_\eps^{\rm l}$ (resp.\ $S_\eps W_\eps^{\rm r}$) converges in distribution to a Brownian web $\Wi^{\rm l}$ (resp.\ $\Wi^{\rm r}$) consisting of coalescing Brownian motions with drift $-1$ (resp.\ $+1$). Furthermore,
$(S_\eps W_\eps^{\rm l}, S_\eps W_\eps^{\rm r})$ converges jointly to a pair of coupled Brownian webs $(\Wi^{\rm l}, \Wi^{\rm r})$, called the {\em left-right web}, where paths in $\Wi^{\rm l}$ and $\Wi^{\rm r}$ undergo sticky interaction. The fact that paths in $N_\eps$ cannot cross without meeting ensures that there is a well-defined leftmost (resp.\ rightmost) path in $N_\eps$ starting from each $(x,t)\in \Z^2_{\rm even}$.

Note that paths in the discrete net $N_\eps$ can always be constructed by concatenating paths in the pair of discrete webs $(W_\eps^{\rm l}, W_\eps^{\rm r})$.
This motivated the construction of the Brownian net $\Ni$ in \cite{SS08} by concatenating paths in the left-right web $(\Wi^{\rm l}, \Wi^{\rm r})$ and
then taking closure in the path space $(\Pi, d)$ (concatenation between two paths is only allowed at a time of coincidence between the two paths that is
strictly larger than their starting times). This is called the {\em hopping construction} of the {\em standard Brownian net} $\Ni$ from the left-right web $(\Wi^{\rm l}, \Wi^{\rm r})$, where standard refers to the fact that the Brownian paths in $\Wi^{\rm l}$ and $\Wi^{\rm r}$ have drift $-1$ and $+1$ respectively.

In \cite{SS08}, the proof of convergence of the rescaled discrete nets $(S_\eps N_\eps)_{\eps\in (0,1)}$ to the Brownian net $\Ni$ consists of
three steps: a lower bound via hopping, an upper bound via wedges, and tightness via leftmost and rightmost paths.

To show the lower bound that any subsequential weak limit of $(S_\eps N_\eps)_{\eps\in (0,1)}$ contains a copy of $\Ni$, thanks to the hopping construction of the Brownian net, it suffices to show that the
concatenation between any two paths in $(\Wi^{\rm l}, \Wi^{\rm r})$ can be approximated by the concatenation of paths in the rescaled discrete webs
$(S_\eps W^{\rm l}_\eps, S_\eps W^{\rm r}_\eps)$. This is an easy consequence of the convergence $(S_\eps W_\eps^{\rm l}, S_\eps W_\eps^{\rm r})\Rightarrow (\Wi^{\rm l}, \Wi^{\rm r})$.

To show the upper bound that any subsequential weak limit of $(S_\eps N_\eps)_{\eps\in (0,1)}$ contains no additional paths besides $\Ni$, the proof
in \cite{SS08} (also in \cite{EFS17}) relied on the so-called {\em wedge characterization} of the Brownian net given in \cite{SS08} (see also the survey \cite{SSS17}), which uses the dual of the left-right web $(\Wi^{\rm l}, \Wi^{\rm r})$ to identify forbidden regions (called wedges) that cannot be entered by any path in $\Ni$. But the definition of wedges also requires $N_\eps$ to have well-defined leftmost and rightmost paths.

Proving tightness for the family of rescaled discrete nets $(S_\eps N_\eps)_{\eps\in (0,1)}$ amounts to controlling the uniform modulus of continuity
of paths in $S_\eps N_\eps$, see \cite[Prop.~B.1]{FINR04} or Prop.~\ref{P:tig} below. Since all paths in $N_\eps$ starting from each $(x,t)$ are bounded between the unique path in $W^{\rm l}_\eps$, resp.\ $W^{\rm r}_\eps$, starting from $(x,t)$, the uniform modulus of continuity of paths in $N_\eps$ is controlled by the uniform modulus of continuity of paths in $W_\eps^{\rm l} \cup W_\eps^{\rm r}$. The tightness of $(S_\eps N_\eps)_{\eps \in (0,1)}$ then follows directly from the tightness of the rescaled discrete webs $(S_\eps W_\eps^{\rm l})_{\eps\in (0,1)}$ and $(S_\eps W_\eps^{\rm r})_{\eps\in (0,1)}$, which follows from the convergence of $S_\eps W_\eps^{\rm l}\Rightarrow \Wi^{\rm l}$ and $S_\eps W_\eps^{\rm r}\Rightarrow \Wi^{\rm r}$.

Of the three steps above, the upper bound and tightness proof depend crucially on the non-crossing property, which implies that every path in
$N_\eps$ is bounded between the leftmost and rightmost path in $N_\eps$ starting from the same point. These arguments have to be replaced by completely new ideas in the crossing case. Only the idea of giving a lower bound via hopping has a natural extension in the crossing case.

\subsubsection{New result: the crossing case} \label{S:crossingstrat}
We now explain our proof strategy for the branching-coalescing systems we consider, where paths can cross without coalescing.
\medskip

\noindent
{\bf Lower Bound.} Since there is no longer a discrete left-right web $(W^{\rm l}_\eps, W^{\rm r}_\eps)$, we will sample
two discrete webs $(W^1, W^2)$ as defined after \eqref{W1}, so that the discrete net $N_\eps$ can still be constructed by concatenating paths
in $W^1$ and $W^2$. We will show in Theorem \ref{T:stick2} that $(S_\eps W^1, S_\eps W^2)$ converges in distribution to a pair of sticky Brownian webs $(\Wi^1, \Wi^2)$, whose construction and characterisation will be recalled in Section \ref{S:sticky}. Similar to the hopping construction of the Brownian net $\Ni$ from the left right web $(\Wi^{\rm l}, \Wi^{\rm r})$, we can also construct a Brownian net by concatenating paths in $\Wi^1 \cup \Wi^2$ and then taking closure. The lower bound is then obtained by showing that every subsequential weak limit of $(S_\eps N_\eps)_{\eps \in (0,1)}$ contains a copy of $\Ni$ constructed from $(\Wi^1, \Wi^2)$ via the hopping construction. The key step in the proof is the convergence $(S_\eps W^1, S_\eps W^2)\Rightarrow (\Wi^1, \Wi^2)$, which requires finding a suitable harmonic function for a random walk with sticky interaction at the origin (see Prop.~\ref{P:dismart}).
\medskip

\noindent
{\bf Upper Bound.} To show that every subsequential weak limit of $(S_\eps N_\eps)_{\eps \in (0,1)}$ contains no additional paths besides
the Brownian net $\Ni$ constructed in the lower bound, we use a different density argument than in \cite{SS08}. The key observation is that we can dominate the discrete net $N_\eps$ by an auxiliary branching-coalescing system which we call the {\em Bernoulli net} $\wt N_\eps$, see Section \ref{S:bern}. In particular, the Bernoulli net allows branching with arbitrary number of offsprings. What is
special about the Bernoulli net $\wt N_\eps$ is that we can identify explicitly a product invariant law for the branching-coalescing random walks in $\wt N_\eps$,
which converges under diffusive scaling to a Poisson point process with intensity $2$ on $\R$ and is precisely the invariant law for the branching-coalescing Brownian motions in the Brownian net $\Ni$. This convergence can then be used to rule out the existence of additional paths in the limit besides paths in $\Ni$. A key ingredient in our proof is Theorem \ref{T:density}, which shows that for the set of paths in the Bernoulli net $\wt N_\eps$ started from every point in $\Z$ at time $0$, their density at large time $t>0$ still decays at the rate of $1/\sqrt{t}$, which is comparable to coalescing random walks without branching.
\medskip

\noindent
{\bf Tightness.} In contrast to the non-crossing case, where tightness comes for free as a result of tightness of the rescaled left right webs $(S_\eps W^{\rm l}_\eps, S_\eps W^{\rm r}_\eps)_{\eps \in (0,1)}$, in our case, tightness is the most technical part of the proof. The overall strategy follows the multiscale argument developed in \cite{BMSV06}, which proves tightness for the rescaled discrete webs (coalescing random walks) $(S_\eps W)_{\eps \in (0,1)}$ under the same assumption on $a(\cdot)$ as in \eqref{passum}. We go through an exponentially increasing sequence of space-time scales. From one scale to the next, we exploit the coalescence to reduce the density of walks as the scale increases, and at the same time, we control the probability that at least one of the walks makes a large excursion that leads to a poor modulus of continuity. Although this overall strategy is not new, the technical difficulties of incorporating branching are serious and we find it the most difficult part of the proof and took us the longest to complete.

For branching-coalescing random walks with weak branching as we consider, coalescence dominates on small scales. So the density reduction argument can still be applied. However, the difficulty is that even though the branching is weak, a single random walk could still have an arbitrarily large number of
descendants on large time scales. Controlling how far a single random walk and any of its descendants can travel within a short time becomes the key challenge in the crossing case, while in the
non-crossing case, it suffices to consider the large deviation of the leftmost and the rightmost random walk among all the descendants. We overcome this difficulty by controlling the number of {\em relevant branching points} in the discrete net $N_\eps$ (see Theorem \ref{T:kRBP}), which are space-time points where a descendant of
the original random walk (starting at some initial time $S$) branches into two random walks, and each random walk has a descendant such that the two descendants do not meet before some terminal time $T$. The relevant branching points are the branching points that lead to new random walk positions at time $T$. To control the large deviation of all descendants of a single random walk, it suffices to group the descendants according to the choice made at each relevant branching point and then control the large deviation of a random walk sampled from each group.

%\begin{remark}
%Our proof of the upper bound and tightness can be extended to the rescaled Bernoulli net $S_\eps \widetilde N_\eps$. Since the lower bound holds trivially for $S_\eps \widetilde N_\eps$ because it contains more paths than $S_\eps N_\eps$, it follows that Theorem \ref{T:net} also holds for $S_\eps \widetilde N_\eps$.
%\end{remark}

\begin{remark}
	The universality of the Brownian net $\Ni$ should hold beyond the family of branching-coalescing random walks considered in Theorem \ref{T:net}. For example, we could consider branching-coalescing random walks in continuous time that arise naturally as the dual of biased voter models \cite{SSY19, SSY21}. In particular, consider the setting where each random walk jumps at rate $1$ with kernel $a(\cdot)$ satisfying the assumptions in \eqref{passum}, and at rate $\eps>0$, each random walk gives birth to a new particle which immediately makes a jump with kernel $a(\cdot)$, while particles coalesce when they meet. We believe our proof can be extended without much difficulty to this continuous time model. One element becomes even simpler, namely that by a generator calculation, it can be shown that this branching-coalescing random walk admits a product invariant measure. Therefore we no longer need the Bernoulli net introduced in Section \ref{S:bern}.
	
	The Brownian net should also be the scaling limit of branching-coalescing random walks with non-binary branching, provided that the branching probability is of the same order as the spatial scaling parameter $\eps$. If the corresponding discrete net could be stochastically dominated by the auxiliary Bernoulli net constructed in Section~\ref{S:bern}, which allows branching into $k$ particles with probability of the order $\eps^{k-1}$ for $k\geq 2$, then {\yu our proof method still applies}. Otherwise new techniques will be needed, especially to prove the upper bound and tightness.
\end{remark}

%{\yu\begin{remark}
%	If, instead of the binary branching in (\ref{a2}), we allow for multiple branching arrows, our method for proving Theorem~\ref{T:net} remains applicable, provided that the discrete net is stochastically dominated by the auxiliary Bernoulli net constructed in Section~\ref{S:bern}, in which the probability of observing $k$ branching arrows is of order $\eps^{k-1}$.
%	In contrast, if the probability of multiple branching is of order $\eps$, the convergence to the Brownian net is currently open. In particular, the proofs of tightness and of the upper bound based on the Bernoulli net would require substantial modifications.
%\end{remark}}

\subsection{Outline}

The rest of the paper is organized as follows. In Section~\ref{S:background}, we recall some necessary background on the Brownian web, the Brownian net, and sticky Brownian webs. Assuming tightness of the family of rescaled discrete nets $(S_\eps \overline{N}_\eps)_{\eps\in (0,1)}$, we show in Section~\ref{S:lowproof} that every subsequential weak limit of $S_\eps \overline{N}_\eps$ contains a copy of the standard Brownian net $\Ni$. In Section~\ref{S:bern}, we introduce an auxiliary Bernoulli net which admits an explicit invariant measure and dominates the discrete net. In Section~\ref{S:upproof}, we use the domination by the Bernoulli net to show that every subsequential weak
limit of $S_\eps \overline{N}_\eps$ contains no additional paths besides the Brownian net $\Ni$. In Section~\ref{S:tight}, we prove the tightness for the family of rescaled discrete nets $(S_\eps \overline{N}_\eps)_{\eps\in (0,1)}$. Lastly, in Appendix, we sketch
the proof of the equivalence between two alternative constructions of a pair of sticky Brownian webs, and the construction of the Brownian net from sticky Brownian webs.

\section{Background} \label{S:background}

We first recall in Subsection~\ref{S:topo} the space of compact sets of paths in which the Brownian web $\Wi$ and the Brownian net $\Ni$ take their values. Then in Subsections~\ref{S:BW}, \ref{S:sticky} and \ref{S:BN}, we introduce necessary results on the Brownian web, sticky Brownian webs, and the Brownian net, respectively. Proposition \ref{P:equiv} below on the equivalence of two alternative constructions of a pair of sticky Brownian webs and
Theorem \ref{T:Nhop} below on the construction of the Brownian net from sticky Brownian webs have not been proved explicitly in the literature. We will sketch their proof in Appendix.% \ref{S:App}.

\subsection{Space of compact sets of paths} \label{S:topo}

We recall here the Polish space $\Hi$ of compact sets of paths, where the Brownian web and the Brownian net take their values. This space was first introduced in \cite{FINR04}, and we follow the variant defined in \cite[Appendix]{SS08}.

Let $\Rc:=\R^2\cup\{(\pm\infty,t):t\in\R\}\cup\{(\ast,\pm\infty)\}$ be the compactification of $\R^2$ obtained by equipping the set
$\Rc$ with a
topology such that $(x_n,t_n)\to(\pm\infty,t)$ if $x_n\to\pm\infty$ and
$t_n\to t\in\R$, and $(x_n,t_n)\to(\ast,\pm\infty)$ if $t_n\to\pm\infty$
(regardless of the behavior of $x_n$). Such a compactification can be achieved
by taking the completion of $\R^2$ with respect to the metric
\be
\label{rho}
\rho((x_1,t_1),(x_2,t_2))
=|\Tet_1(x_1,t_1)-\Tet_1(x_2,t_2)|\vee|\Tet_2(t_1)-\Tet_2(t_2)|,
\ee
where the map $\Tet=(\Tet_1,\Tet_2)$ is defined by
\be\label{Tet}
\Tet(x,t)=\big(\Tet_1(x,t),\Tet_2(t)\big)
:=\Big(\frac{\tanh(x)}{1+|t|},\tanh(t)\Big).
\ee
We can think of $\Rc$ as the image of $[-\infty,\infty]^2$ under the
map $\Tet$ (see Figure \ref{fig:comp}).

\begin{figure}
\begin{center}
\setlength{\unitlength}{.7cm}
\begin{picture}(10,8)(-5,-4)
\linethickness{.4pt}
\qbezier(0,-3)(0,0)(0,3)
\qbezier(-3,0)(0,0)(3,0)
\linethickness{.4pt}
\qbezier(0,-3)(-6,0)(0,3)
\qbezier(0,-3)(-5.3,0)(0,3)
\qbezier(0,-3)(-4.5,0)(0,3)
\qbezier(0,-3)(-3,0)(0,3)
\qbezier(0,-3)(3,0)(0,3)
\qbezier(0,-3)(4.5,0)(0,3)
\qbezier(0,-3)(5.3,0)(0,3)
\qbezier(0,-3)(6,0)(0,3)
\qbezier(-1.9,-1.8)(0,-1.8)(1.9,-1.8)
\qbezier(-2.2,-1.5)(0,-1.5)(2.2,-1.5)
\qbezier(-2.65,-1)(0,-1)(2.65,-1)
\qbezier(-2.65,1)(0,1)(2.65,1)
\qbezier(-2.2,1.5)(0,1.5)(2.2,1.5)
\qbezier(-1.9,1.8)(0,1.8)(1.9,1.8)

\put(0,-3){\circle*{.25}}
\put(0,3){\circle*{.25}}
\put(0,0){\circle*{.25}}
\put(2.2,1.5){\circle*{.25}}
\put(-2.65,-1){\circle*{.25}}

\put(-.6,-3.5){$(\ast,-\infty)$}
\put(-.6,3.5){$(\ast,+\infty)$}
\put(-1.3,-.5){$(0,0)$}
\put(2.3,1.6){$(+\infty,2)$}
\put(-5.3,-1.1){$(-\infty,-1)$}

\end{picture}
\caption[The compactification $\Rc$ of $\R^2$.]{The compactification
$\Rc$ of $\R^2$ under the map $\Theta$.}\label{fig:comp}
\end{center}
\end{figure}

A \emph{path} in $\Rc$ is any continuous function
$\pi:[\sig_\pi,\infty]\to[-\infty,\infty]\cup\{\ast\}$, with
$\sig_\pi\in[-\infty,\infty]$, such that
$\pi:[\sig_\pi,\infty]\cap\R\to[-\infty,\infty]$ is continuous, and
$\pi(\pm\infty)=\ast$ whenever $\pm\infty\in[\sig_\pi,\infty]$.  Equivalently,
if we identify $\Rc$ with the image of $[-\infty,\infty]^2$ under the map
$\Tet$, then $\pi$ can be identified with its graph, which is a compact subset of $\Tet([-\infty,\infty]^2)$ with the property that
for each $t\in[\Tet_2(\sig_\pi),\Tet_2(\infty)]$, there is exactly one $x$ with $(x,t)$ on the graph of $\pi$. We will often identify $\pi$ with
its graph in the rest of the paper.

Let $\Pi$ denote the space of all paths in $\Rc$, equipped with the metric
\bc\label{dPi}
\dis d(\pi_1,\pi_2)
&:=&\dis|\Tet_2(\sig_{\pi_1})-\Tet_2(\sig_{\pi_2})|\\[5pt]
&&\dis\quad\vee\sup_{t\geq\sig_{\pi_1}\wedge\sig_{\pi_2}}
\big|\Tet_1\big(\pi_1(t\vee\sig_{\pi_1}),t)
-\Tet_1(\pi_2(t\vee\sig_{\pi_2}),t\big)\big|.
\ec
The space $(\Pi,d)$ is complete and separable. Note that paths
converge in $(\Pi,d)$ if and only if their starting times converge and
the functions converge locally uniformly on $\R$.

Let $\Hi$ be the space of compact subsets of $(\Pi, d)$,
equipped with the Hausdorff metric $d_\Hi$. Namely, for $\Xc_1, \Xc_2\in \Hi$,
\be\label{Hmet}
d_{\Hi}(\Xc_1,\Xc_2)=\sup_{x_1\in \Xc_1}\inf_{x_2\in \Xc_2}d(x_1,x_2)
\vee\sup_{x_2\in \Xc_2}\inf_{x_1\in \Xc_1}d(x_1,x_2).
\ee
The metric space $(\Hi, d_\Hi)$ is also complete and separable. For more properties,
see e.g.\ \cite[Appendix B]{SSS14}. Let $(\widehat \Hi, d_{\widehat \Hi})$ be defined similarly to $(\Hi, d_\Hi)$,
except that all paths run backwards in time.

\begin{remark}\label{R:Neps}
	To view the discrete nets $N_\eps$ as random variables taking values in $\Hi$, we need to
	take the closure of $N_\eps$ in $(\Pi, d)$ and show that $\overline{N}_\eps$ is compact.
	This is easily verified if the jump distribution $a(\cdot)$ in \eqref{passum} has a finite
	first moment, see e.g.\ \cite[Lemma 1.1]{NRS05}. It is then easily seen
	that $\overline{N}_\eps\backslash N_\eps$ only contains the trivial paths $\pi$ with
	$\sig_\pi\in\{\pm\infty\}\cup\Z$, and either $\pi(t) \equiv -\infty$ or $\pi(t)\equiv \infty$ for all $t\in [\sigma_\pi, \infty)$.
	To keep the notation simple, we will write $N_\eps$ instead of $\overline{N}_\eps$ in the rest of the paper.
\end{remark}

We will adopt the following notational convention:
\begin{itemize}
	\item For $\Xc\in \Hi$ and $A\subset \Rc$, let $\Xc(A)$ denote the subset of paths in $\Xc$
	with starting points in $A$. Sometimes we will also write $\Xc(A)$ as $\Xc^A$.
	
	\item For $A=\{z\}$, we simply write $\Xc(z)$ instead of $\Xc(\{z\})$.
\end{itemize}

\begin{remark}
The space $(\Hi, d_\Hi)$ was introduced in \cite{FINR04} as the space where the Brownian web $\Wi$ takes its values, because in the collection
of coalescing Brownian motions starting from every point in space and time, there are always random space-time points (called special points) from which multiple Brownian paths originate. Instead of selecting one path from each such special point, which is somewhat arbitrary and does not lead to a nice space for proving convergence, the Brownian web $\Wi$ includes all paths starting from every point. The space $\Hi$ of sets of paths and its accompanying Hausdorff topology is then the natural choice.
\end{remark}

\subsection{The Brownian web} \label{S:BW}

We recall from \cite{FINR04} the characterization
of the Brownian web $\Wi$ as a random variable taking values in the complete separable metric space $(\Hi, d_\Hi)$.

\bt[Characterization of the Brownian web] \label{T:web}
There exists an $\Hi$-valued random variable $\Wi$, called the standard Brownian web, whose
distribution is uniquely determined by the following properties:
\begin{itemize}
	\item[{\rm (i)}] For each deterministic $z\in\R^2$, almost surely there is a
	unique path $\pi_z\in\Wi(z)$.
	
	\item[{\rm (ii)}] For any finite deterministic set of points
	$z_1,\ldots,z_k\in\R^2$, the paths $(\pi_{z_1},\ldots,\pi_{z_k})$ are
	distributed as a collection of coalescing standard Brownian motions.
	
	\item[{\rm (iii)}] For any deterministic countable dense subset $\Di
	\subset\R^2$, almost surely, $\Wi$ is the closure of
	$\{\pi_z:z\in\Di\}$ in $(\Pi, d)$.
\end{itemize}
\et
\noindent
The Brownian web $\Wi$ can be coupled to a dual Brownian web $\widehat \Wi \in \widehat \Hi$, which consists of coalescing Brownian motions running backwards in time and has the same distribution as $\Wi$ if space-time is rotated by $180^\circ$ around the origin. Almost surely, $\Wi$ and $\widehat \Wi$ determine each other by the condition that
paths in $\Wi$ cannot cross paths in $\widehat \Wi$, i.e., if $\pi \in \Wi$ and $\hat\pi \in \widehat \Wi$, then we cannot find $s<t$ with $(\pi(s)-\hat\pi(s))(\pi(t) - \hat \pi(t))<0$. Furthermore, paths in $\Wi$ and $\widehat \Wi$ interact through Skorohod reflection. For further properties and alternative characterizations of $(\Wi, \widehat \Wi)$, see the survey article \cite{SSS17} and the references therein.

The following analogue of Theorem \ref{T:net} for the Brownian web $\Wi$ is due to \cite{NRS05}, which assumed finite $5$-th moment for the random walk
increments in the proof of tightness, and it was then improved in \cite{BMSV06} to the nearly optimal assumption of a finite $(3+\eta)$-th moment.

\bt[Convergence to the Brownian web] \label{T:convweb}
Let $W$ be the discrete web defined in \eqref{W1} with jump distribution $a(\cdot)$
satisfying (\ref{passum}). Then
\be\label{nonsimpweb}
S_\eps W \Astoo{\eps} \Wi,
\ee
where $\Wi$ is the standard Brownian web, and $\Rightarrow$ denotes weak convergence of $\Hi$-valued random variables.
\et

\subsection{Sticky Brownian webs}\label{S:sticky}

The Brownian net $\Ni$ was first constructed in \cite{SS08} (see \cite{NRS10} for a later alternative construction), which generalises the Brownian web $\Wi$ by allowing paths to branch. An intermediate object in the construction of the Brownian net $\Ni$ in \cite{SS08} is a pair of sticky Brownian webs
$(\Wi^{\rm l}, \Wi^{\rm r})$ with drift $-1$ and $+1$ respectively. For our purpose, we will give an alternative construction using a pair
of {\em sticky Brownian webs} $(\Wi^1, \Wi^2)$ with the same drift $0$.

We first motivate the sticky Brownian webs (with drift $0$) through discrete approximation.
Recall the definition of the discrete net $N_\eps$ in Section \ref{S:disc}, which is constructed from an i.i.d.\ family of
random pairs $(\omega^1(x, t), \omega^2(x,t))_{(x,t)\in \Z^2}$ with common law $a^{(2)}$ that satisfies \eqref{a2}. Each
family $(\omega^i(x, t))_{(x,t)\in \Z^2}$, $i=1,2$, determines a discrete web $W^i$, and by Theorem \ref{T:convweb}, $S_\eps W^i$ converges in distribution to a Brownian web $\Wi^i$. In the special case when $\omega^i(x,t)\in \{\pm 1\}$ so that the paths make nearest neighbour jumps, and when we only consider paths starting from the even lattice sites $\Z^2_{\rm even}:=\{(x, t)\in \Z^2: x+t \mbox{ is even}\}$ to avoid parity issues, results from \cite{SSS14} imply that the pair
$(S_\eps W^1, S_\eps W^2)$ converges to a pair of coupled Brownian webs $(\Wi^1, \Wi^2)$, called {\em sticky Brownian webs} because paths between $\Wi^1$ and $\Wi^2$ undergo sticky interaction.

Instead of the discrete webs $(W^1, W^2)$ defined above,  \cite{SS08} considered the discrete webs $W^{\rm l}_\eps$ and $W^{\rm r}_\eps$, induced respectively by
\begin{equation}\label{lromega}
	\omega^{\rm l}(x, t) =  \omega^1(x,t) \wedge \omega^2(x,t), \qquad \omega^{\rm r}(x, t) =  \omega^1(x,t) \vee \omega^2(x,t),
\end{equation}
so that whenever two (nearest neighbour) arrows originate from $(x,t)$, paths in $W^{\rm l}_\eps$ follow the arrow to the left and paths in $W^{\rm r}_\eps$ follow the arrow to the right. It was shown in \cite{SS08} that $(S_\eps W^{\rm l}_\eps, S_\eps W^{\rm r}_\eps)$ converges to a pair of sticky Brownian webs $(\Wi^{\rm l}, \Wi^{\rm r})$ with drift $-1$ and $+1$ respectively, called the {\em left-right Brownian web}. The key observation of \cite{SS08} is that, since the discrete net $N_\eps$ can be constructed by concatenating paths in the two discrete webs $W^{\rm l}_\eps$ and $W^{\rm r}_\eps$ (i.e., switching back and forth between arrows induced by either $\omega^{\rm l}$ and $\omega^{\rm r}$), the limit of $S_\eps N_\eps$ can also be recovered by concatenating paths in the pair of sticky Brownian webs $(\Wi^{\rm l}, \Wi^{\rm r})$. This was called the {\em hopping construction} of the Brownian net $\Ni$. Since the discrete net $N_\eps$ can also be constructed by concatenating paths in the two discrete webs $W^1$ and $W^2$ (i.e., switching back and forth between arrows induced by either $\omega^1$ and $\omega^2$), we can also recover $\Ni$ by concatenating paths in the pair of sticky Brownian webs $(\Wi^1, \Wi^2)$ with drift $0$. This is the content of Theorem \ref{T:Nhop} below.

We will mostly work with the sticky Brownian webs $(\Wi^1, \Wi^2)$ with drift $0$, for which we will need two equivalent characterizations. The first was introduced in \cite{HW09b}, which we call the {\em martingale characterization}. We will use it to prove $(S_\eps W^1, S_\eps W^2) \Rightarrow (\Wi^1, \Wi^2)$ for general jump distribution $a(\cdot)$ satisfying \eqref{passum}. The second was introduced in \cite{SSS14}, which we call the {\em marking construction}. It constructs the pair of sticky Brownian webs $(\Wi^1, \Wi^2)$ and the Brownian net $\Ni$ simultaneously through a Poisson marking procedure that was first developed in \cite{NRS10}. This allows us to recover the Brownian net $\Ni$ from $(\Wi^1, \Wi^2)$ through the {\em hopping construction} sketched above, stated in Theorem \ref{T:Nhop} below. We will show that these two characterizations of the sticky Brownian webs give the same object, so that we can choose the most convenient characterization depending on the context.
\medskip

\noindent
{\bf Martingale Characterization.} Following \cite{HW09a}, $(Y^1,Y^2)=(Y^1_t,Y^2_t)_{t\geq 0}$ is called a pair
of \emph{sticky Brownian motions with parameter
	$\theta\geq0$} if $Y^1$ and $Y^2$ are standard Brownian motions w.r.t.\
some common filtration $(\Fi_t)_{t\geq 0}$, and the following processes are martingales
w.r.t.\ this filtration
\be\ba{rl}\label{HWmart}
{\rm(i)}&\dis Y^1_tY^2_t- \int_0^t1_{\{Y^1_s=Y^2_s\}}\di s,\\[5pt]
{\rm(ii)}&\dis|Y^1_t-Y^2_t|-2\tet\int_0^t1_{\{Y^1_s=Y^2_s\}}\di s.
\ec
By \cite[Prop.~3]{HW09b}, we have the following fact.

\bp[Characterization of sticky Brownian motions]\label{P:stickMP}
For each $\tet\geq 0$ and $(y_1,y_2)\in\R^2$, there exists a pair $(Y^1,Y^2)$
of sticky Brownian motions with parameter $\theta\geq0$ and initial
state $(Y^1_0,Y^2_0)=(y_1,y_2)$, and such a pair is unique in law.
\ep
\noindent
We can readily extend the definition of sticky Brownian motions to the case where $Y^1$
and $Y^2$ start at different times. Note that the smaller the value of $\theta$ is, the stronger
is the stickiness. When $\tet=0$, $(Y^1,Y^2)$ is a pair of coalescing Brownian motions.
The case where $Y^1,Y^2$ have different drifts is treated in \cite[Prop.~14]{HW09b}. In particular,
when $Y^1$ has a smaller drift than $Y^2$ and the pair $(Y^1, Y^2)$ stays ordered after they meet,
$(Y^1, Y^2)$ is equal in distribution to the unique solution of the {\em left-right SDE} in \cite[Sec.~2]{SS08}.

Note that each Brownian web $\Wi^i$ defines in a natural way a filtration $\Fi^i_t$, which is the
\si-field generated by the restriction of all paths
in $\Wi^i$ to times $\leq t$. Let  $\Fi_t:= \Fi^1_t\vee\Fi^2_t:=\sig(\Fi^1\cup\Fi^2)$ be the smallest
\si-field containing  $\Fi^1_t$ and $\Fi^2_t$. Then $\Fi_t$ also forms a filtration. We can now state the
martingale characterization of a pair of sticky Brownian webs $(\Wi^1, \Wi^2)$, which is a variant of
\cite[Thm.~4]{HW09b} and \cite[Thm.~1]{SS19}.

\bt[Martingale characterization of sticky Brownian webs]\label{T:martc}
There exists an $\Hi\times\Hi$-valued random variable $(\Wi^1,\Wi^2)$, called a pair of sticky Brownian webs
with parameter $\theta\geq 0$,  whose distribution is uniquely determined by the following properties:
\begin{itemize}
	\item[{\rm(i)}] Both $\Wi^1$ and $\Wi^2$ are distributed as the standard Brownian web;
	
	\item[{\rm(ii)}] For each deterministic pair $z=(x, t), z'=(x', t')\in \R^2$, the pair of paths
	$(X^1_z, X^2_{z'})$, with $X^1_z\in \Wi^1(z)$ and $X^2_{z'}\in \Wi^2(z')$, is distributed
	as a pair of sticky Brownian motions with parameter $\theta$ relative to the filtration $(\Fi_u)_{u\in \R}$;
	
	\item[{\rm (iii)}] %$($Co-adaptedness$)$
	For each deterministic pair $z=(x,t), z'=(x', t')\in \R^2$, the pair $(X^i_z, X^j_{z'})$,
	with $X^i_z \in \Wi^i(z)$ and $X^j_{z'} \in \Wi^j(z')$ for $i, j=1,2$, is a~Markov~process~relative~to~the~filtration $(\Fi_u)_{u\in \R}$, i.e., $\P\big((X^i_z(s), X^j_{z'}(s))_{s\geq u}\in\cdot\hspace{1pt}\big| \Fi_u\big)=\P\big((X^i_z(s), X^j_{z'}(s))_{s\geq u}\in\cdot\hspace{1pt}\big| X^i_z(u), X^j_{z'}(u) \big)$ for any $u\geq t\vee t'$.
	%if $X^i_z \in \Wi^i(z)$
	%and $X^j_{z'} \in \Wi^j(z')$ for $i, j=1,2$, then the pair $(X^i_z, X^j_{z'})$ is a Markov process relative to the filtration
	%$(\Fi_u)_{u\in \R}$.
\end{itemize}
\et
\noindent
This formulation, in particular, condition (iii), follows that of \cite[Thm.~1]{SS19} for the {\em left-right Brownian web}, cf.~\eqref{lromega}. The
formulation for sticky Brownian webs was discussed in the proof of \cite[Thm.~1]{SS19}.
Notably, condition (iii), which plays the same role as the so-called {\em co-adaptedness} property formulated in \cite[Thm.~4]{HW09b},
%Here, Condition (iii) is first employed in \cite{SS19}, and plays the same role of the so-called {\em co-adaptedness} property formulated in \cite[Thm.~4]{HW09b}.   co-adaptedness is formulated differently than in \cite[Thm.~4]{HW09b} in order to
makes it more clear how \cite[Thm.~76]{H07} can be applied to show that the law of $(\Wi^1, \Wi^2)$ is uniquely determined by conditions
(i)-(ii) plus (iii) (co-adaptedness) (see the proof of \cite[Thm.~1]{SS19} for more details).
\bigskip

\noindent
{\bf Marking Construction.} We now sketch the marking construction of the sticky Brownian webs $(\Wi^1, \Wi^2)$ from
\cite[Section 3]{SSS14}, which is an extension of the marking construction of the Brownian net $\Ni$ in \cite{NRS10}. We refer the reader to these
references for details.

Given a Brownian web $\Wi$ and its dual $\widehat \Wi$, each pair $(\pi, \hat\pi)\in (\Wi, \widehat \Wi)$
interacts through Skorohod reflection. In particular, there is a well-defined {\em intersection local time measure} $\ell$ on the intersection between the graphs of $\pi$ and $\hat\pi$, for $\pi\in \Wi$ and $\hat\pi \in \widehat \Wi$ (see \cite[Prop.~3.4]{SSS14}). Each point of intersection $(x,t)$ between a path $\pi\in \Wi$ with $\sigma_\pi <t$, and $\hat\pi\in \widehat\Wi$ with $\sigma_{\hat\pi}>t$, is called a {\em special point of $\Wi$} of type $(1,2)$ because $\pi$ is the unique
path in $\Wi$ passing through $(x,t)$, and besides $\pi$, there is a second path $\pi' \in \Wi$ that starts from $(x,t)$.
We then define a Poisson point set $S$ on all such $(1,2)$ points with intensity measure $\theta \ell$. Given $\Wi$ and
$S$, we then construct a second set of paths $\Wi'$ such that paths in $\Wi'$ follow paths in $\Wi$, except when a path $p\in \Wi'$ enters a point $(x,t)\in S$, instead of continuing along $\pi\in \Wi$ that passes through $(x,t)$, $p$ switches to the
second path $\pi'\in \Wi$ that starts from $(x,t)$. It was shown in \cite[Theorem~3.5]{SSS14} that $\Wi'$ is well-defined, and the pair $(\Wi, \Wi')$ is called a {\em pair of sticky Brownian webs with coupling parameter $\theta$}. If we include both the path that continues along $\pi$ and the path
that continues along $\pi'$, then we obtain a Brownian net, see \cite[Theorem~3.5 and 6.15]{SSS14}.

We will show in Appendix that the two definitions of sticky Brownian webs above are equivalent.

\bp[Equivalence] \label{P:equiv}
Let $(\Wi^1, \Wi^2)$ be a pair of sticky Brownian webs with parameter $\theta\geq 0$ defined in Theorem \ref{T:martc}, and let $(\Wi, \Wi')$ a pair of sticky Brownian webs with coupling parameter $\theta$ defined via the marking construction above. Then $(\Wi^1, \Wi^2)$ and $(\Wi, \Wi')$ have the same distribution.
\ep

Both the martingale characterization and the marking construction of sticky Brownian webs can be extended to the case with non-zero drifts and parameter $\theta\geq 0$. For simplicity, we have restricted to the case of zero drift. When $(\Wi^1, \Wi^2)$ is a {\em left-right Brownian web}, where $\Wi^1$ is a Brownian web with a smaller drift than $\Wi_2$ and
any pair $(\pi_1, \pi_2)\in (\Wi^1, \Wi^2)$ stay ordered once they meet, it was shown in \cite[Lemma 6.18]{SSS14} that the marking construction of the left-right web is equivalent to a characterization similar to Theorem \ref{T:martc}.

\subsection{The Brownian net} \label{S:BN}

We first briefly recall the construction of the Brownian net $\Ni$ from \cite{SS08}, which generalises the Brownian web $\Wi$ by allowing branching,  and then discuss how this construction should be modified for our setting.

The original construction of $\Ni$, called the {\em hopping construction} in \cite[Thm.~1.3]{SS08}, was motivated by the discrete analogue discussed after \eqref{lromega}, where $\Ni$ was constructed by
concatenating paths in the left-right web $(\Wi^{\rm l}, \Wi^{\rm r})$. The left-right web $(\Wi^{\rm l}, \Wi^{\rm r})$ is characterized by the properties that $\Wi^{\rm l}$, resp.\ $\Wi^{\rm r}$, is a Brownian web with drift $-1$, resp.\ $+1$, and any pair of paths $(l, r)\in (\Wi^{\rm l}, \Wi^{\rm r})$ with
deterministic starting points evolve as a pair of sticky Brownian motions with drifts $(-1, 1)$ and stays ordered after they meet. More generally, the left-right web $(\Wi^{\rm l}, \Wi^{\rm r})$ could have arbitrary drifts $(a, b)$ with $a<b$. We refer the reader to \cite{SS08} for the precise characterization of $(\Wi^{\rm l}, \Wi^{\rm r})$.

We will need a modification of this construction and show that $\Ni$ can also be constructed by concatenation (hopping) between paths in the pair of sticky Brownian webs $(\Wi^1, \Wi^2)$ from Theorem \ref{T:martc}. First, we give the precise definition of hopping between paths in $(\Wi^1, \Wi^2)$, which also applies to $(\Wi^{\rm l}, \Wi^{\rm r})$. Given two paths $\pi, \pi'\in \Pi$,
we call $t\in\R$ an \emph{intersection time} if $\sig_\pi,\sig_{\pi'}<t<\infty$ and $\pi(t)=\pi'(t)$. At such an intersection time
$t$, we can then define a new path $\pi''$ by concatenating the piece of $\pi$ before $t$ with the piece of $\pi'$ after
$t$, i.e., by setting $\pi''(s):=\pi(s)$ for $s\in [\sigma_\pi, t]$ and $\pi''(s)=\pi'(s)$ for $s\in [t,\infty)$. For
any collection of paths $\Ai\sub\Pi$, we let $\Hi_{\rm int}(\Ai)$ denote the
smallest set of paths containing $\Ai$ that is closed under such `hopping'
from one path onto another at intersection times, i.e., $\Hi_{\rm int}(\Ai)$
is the set of all paths $\pi\in\Pi$ of the form
\be\label{hopmeet}
\pi(s) := \pi_k(s) \mbox{ for }  s\in[t_{k-1},t_k], \quad 1\leq k\leq m,
\ee
where $\pi_1,\ldots,\pi_m\in\Ai$, $\sig_{\pi_1}=t_0<\cdots<t_m=\infty$,
and $t_k$ is an intersection time of $\pi_k$ and $\pi_{k+1}$ for each
$k=1,\ldots,m-1$.

In \cite[Thm.~1.3 and Prop.~1.4]{SS08}, the Brownian net $\Ni$ was constructed by setting it to be the closure
\begin{equation}\label{Nconst}
	\Ni :=\ov{\Hi_{\rm int}(\Wi^{\rm l} \cup\Wi^{\rm r})},
\end{equation}
where $(\Wi^{\rm l}, \Wi^{\rm r})$ is a left-right Brownian web with drifts $(-1, 1)$. If $(\Wi^{\rm l}, \Wi^{\rm r})$
is a left-right web with drifts $(-\theta, \theta)$, then this construction defines the {\em Brownian net $\Ni_\theta$ with
	branching parameter $\theta$}, which we will simply denote by $\Ni$ when $\theta=1$. We call $\Ni$ the {\em standard Brownian net}. It is
easily seen that $\Ni_\theta$ is the image of $\Ni$ under the diffusive scaling map $(x,t) \to (x/\theta, t/\theta^2)$.

We have the following alternative hopping construction of the Brownian net,  whose proof is deferred to Appendix.

\bt[Hopping construction of a Brownian net]
Let\label{T:Nhop} $(\Wi^1,\Wi^2)$ be a pair of sticky Brownian webs with parameter $\theta$.
Then
\be\label{Nhop}
\Ni_\theta:=\ov{\Hi_{\rm int}(\Wi^1\cup\Wi^2)}
\ee
is distributed as a Brownian net with branching parameter $\theta$.

If $K$ is a compact subset of $\Rc$ (see definition above \eqref{rho}), then an analogue of \eqref{Nhop} holds, i.e.,
\begin{equation}\label{NKhop}
	\Ni_\theta(K)=\overline{\Hi_{\rm int}(\Wi^1\cup\Wi^2)(K)}.
\end{equation}
\et
As mentioned before, \cite{NRS10} gave an
alternative {\em marking construction} of $\Ni$. There are also alternative
characterizations of the Brownian net, see \cite{SS08} and the
review article \cite{SSS17}. For more details on the coupling between $(\Wi^1, \Wi^2)$ and
$\Ni$ through Poisson marking of $(1,2)$ points of $\Wi^1$, which will be relevant in the proof of Theorem \ref{T:Nhop},
see \cite[Theorem~3.5, 4.4, and 6.15]{SSS14}.

We will need later the fact that the Brownian net is a maximal collection of paths in the sense that it contains every path in its graph (or image set).
%For any set of paths $\Xc\in\Hi$, we define its image set by
%\be\label{imset}
%\cup \Xc:=\{z\in\R^2:\exists \pi\in \Xc{\rm~s.t.~}z\in \pi \}.
%\ee
Let $\Pi_t:=\{\pi\in\Pi:\sig_{\pi}=t \}$ denote the space of all paths starting at time $t$, and for any $\Xc\in\Hi$, let
\be\label{Xt}
\Xc_t:=\Xc\cap\Pi_t
\ee
denote the subset of paths in $\Xc$ starting at time $t$. Then the following result says that every path $\pi\in \Pi_t$ (which is identified with its graph as a subset of $\R^2$) that is contained in the union of the graphs of $p\in \Ni_t$ must also be a path in $\Ni_t$, see \cite[Prop.~1.13]{SS08}.
\bp[Image set property]\label{P:image}
Let $\Ni$ be the standard Brownian net. Then almost surely for all $t\in[-\infty,\infty]$,
\be
\Ni_t=\Big\{\pi\in\Pi_t:\pi\subset \bigcup_{p\in \Ni_t} p \Big\}.
\ee
\ep

\section{The lower bound}\label{S:lowproof}

The main result of this section is Theorem \ref{T:lower}, which shows that every subsequential weak limit of the rescaled discrete net
$S_{\eps}N_\eps$ contains a copy of the Brownian net $\Ni$ constructed as in Theorem~\ref{T:Nhop} by hopping between paths in
a pair of sticky Brownian webs $(\Wi^1, \Wi^2)$. Our strategy is to first show that the pair of discrete sticky webs $(W^1_\eps, W^2_\eps)$ defined as in \eqref{a2} converges to a pair of sticky Brownian webs $(\Wi^1, \Wi^2)$ as characterised in Theorem~\ref{T:martc}.  {\yu We then further
	strengthen this result by showing that the intersection times between paths in the discrete sticky webs also converge to that of the sticky Brownian webs $(\Wi^1, \Wi^2)$.
	Using this fact, we then show that hopping between paths in $\Wi^1 \cup \Wi^2$ at their intersection times can always be approximated by their discrete counterparts.}
	
%	\sout{We then further strengthen this result by showing that the intersection times between paths in the discrete sticky webs also converge to that of the sticky Brownian webs $(\Wi^1, \Wi^2)$, which are then used to show that the hopping between paths in $\Wi^1 \cup \Wi^2$ at their intersection times can always be approximated by their discrete counterparts.}

\subsection{Convergence of discrete sticky webs}\label{S:WWconv}

The main result of this subsection is the convergence of a pair of discrete sticky webs and the intersection times between paths in
the two discrete webs.

\bt[Convergence of discrete sticky webs]\label{T:stick2}
Let $\eps>0$, and let $(W^1_\eps, W^2_\eps)$ be the pair of discrete sticky webs defined as in \eqref{a2}
with branching probability $\eps>0$ and jump distribution $a(\cdot)$ satisfying \eqref{passum}. Then,
\be
(S_{\eps}W^1_\eps,S_{\eps}W^2_\eps)
\underset{\eps\down0}{\Longrightarrow} (\Wi^1,\Wi^2),
\ee
where $(\Wi^1,\Wi^2)$ is a pair of sticky Brownian webs with parameter $1$ as in Theorem \ref{T:martc}.
\et
\bpro
The convergence of the marginals $S_{\eps}W^1_\eps \Rightarrow \Wi^1$ and $S_{\eps}W^2_\eps \Rightarrow \Wi^2$ follows
from Theorem \ref{T:convweb}. To prove their joint convergence, it suffices to show that every subsequential weak limit of $(S_{\eps}W^1_\eps,S_{\eps}W^2_\eps)$ satisfies conditions (i)-(iii) in Theorem \ref{T:martc}. Condition (i) follows
from the convergence of $S_{\eps}W^1_\eps$ and $S_{\eps}W^2_\eps$.

To verify condition (ii), it suffices to show that a pair of paths $(X^1,X^2)\in(W^1_\eps, W^2_\eps)$, whose starting points
converge under diffusive scaling, must converge to a pair of sticky Brownian motions in the sense of Proposition~\ref{P:stickMP}. By
waiting for the earlier path to reach the starting time of the later path, and by translation invariance, it suffices to consider
paths starting at the same time $0$. A stronger form of such a convergence is established in Theorem \ref{T:stBMcon} below.

The verification of condition (iii) follows the same argument as in the proof of Theorem 2 in \cite[Section 4]{SS19}. Recall that the filtration $(\Fi_u)_{u\in \R}$ is generated by paths in $\Wi^1$ and $\Wi^2$, where we may further restrict to $\Wi^1(\Di)$ and $\Wi^2(\Di)$ for a deterministic countable dense set $\Di\subset \R^2$. For any deterministic $z=(x,t), z'=(x', t')\in \R^2$, to prove that $(X^i_z, X^j_{z'})$ with $X^i_z\in \Wi^i(z)$ and $X^j_{z'}\in \Wi^j(z')$ is a Markov process w.r.t.\ $(\Fi_u)_{u\in \R}$, it suffices to show that for any $s\geq \max\{t, t'\}$ and $z_1, \ldots, z_m\in \Di$ with time coordinates less than $s$, conditioned on $(X^i_z(s), X^j_{z'}(s))$ for some $s\geq \max\{t, t'\}$, the law of $(X^i_z(u), X^j_{z'}(u))_{u\geq s}$ in the future does not depend on the past realisation $(X^i_z(u), X^j_{z'}(u), (X^1_{z_k}(u), X^2_{z_k}(u))_{1\leq k\leq m})_{u<s}$.
This property follows readily from approximation by paths in the discrete sticky webs, for which an analogous property holds. This verifies condition (iii).
\epro

We now state a strong form of convergence of a pair of sticky random walks to sticky Brownian motions, which also includes the convergence
of their intersection times. This will allow us to strengthen Theorem \ref{T:stick2} to include the convergence of intersection times
between paths in $W^1_\eps$ and $W^2_\eps$, see Corollary \ref{C:disweb} below.
{\yu This is needed because the hopping operation from~\eqref{hopmeet} is carried out at intersection times between paths, which is not continuous with respect to the topology for path convergence, see \eqref{dPi}. In particular, if $(\pi_n)_{n\geq 1}$ and $(\pi_n')_{n\geq 1}$ are two sequences of paths
in $(\Pi, d)$ such that $\pi_n\to\pi_\infty$, $\pi_n'\to \pi_\infty'$, and $\pi_n$ does not intersect $\pi_n'$, but $\pi'_\infty$ and $\pi_\infty$ intersect at a single time, then $\Hi_{\rm int}(\{ \pi_n, \pi_n'\}) \not\rightarrow \Hi_{\rm int}(\{\pi_\infty, \pi_\infty'\})$.
%if we only have sequences $W^1_\eps, W^2_\eps \in \mathcal H$ such that $(S_{\eps}W^1_\eps,S_{\eps}W^2_\eps) \Rightarrow (\Wi^1,\Wi^2) ~\text{as } \eps\to0,$ then this does not imply that $\Hi_{\rm int}(S_{\eps}W^1_\eps \cup S_{\eps}W^2_\eps) \Rightarrow \Hi_{\rm int}(\Wi^1 \cup \Wi^2) ~\text{as } \eps\to0,$
Such a convergence statement will hold if we can additionally show that the set of intersection times between $\pi_n$ and $\pi_n'$ converges to that between $\pi_\infty$ and $\pi_\infty'$. This can then be used to deduce Theorem~\ref{T:lower} below.
}

%\sout{This is needed because the hopping operation (see \eqref{hopmeet}) is carried out at intersection times between paths, which is not continuous w.r.t.\ the topology for path convergence.}}

For technical convenience, the random walk paths are defined on $[0,\infty)$ by piecewise constant extension
to non-integer times and regarded as random variables taking values in the Skorohod space $D([0,\infty), \R)$.
\bt[Pair convergence including intersection times]\label{T:stBMcon}
Let $X^1$ and $X^2$ be the unique paths in $W^1_\eps$ and $W^2_\eps$ starting at time zero at positions $x^1_\eps$ and
$x^2_\eps$, respectively, such that $(\eps x^1_\eps, \eps x^2_\eps)\to(y^1, y^2) \in\R^2$ as $\eps \downarrow 0$. Set
\be\label{YZeps}
Y_\eps^i(t):= \eps X^i\big(\lfloor t /\varsigma^2\eps^2\rfloor\big)\quad(i=1,2)
\quand
Z_\eps(t):=\varsigma^2\eps^2 \sum_{s=0}^{\lfloor t /\varsigma^2\eps^{2}\rfloor}1_{\{Y_\eps^1(s)=Y_\eps^2(s)\}}.
\ee
Then as $\eps\down0$, the process $(Y_\eps^1,Y_\eps^2,Z_\eps)$ converges in distribution on the Skorohod space
$\Di([0,\infty), \R^3)$ to a limiting process $(Y^1,Y^2,Z)$, where $(Y^1,Y^2)$
is a pair of sticky Brownian motions with parameter $1$ and
starting point $(y^1,y^2)$ characterised in Proposition~\ref{P:stickMP}, and
\be\label{Zt}
Z(t)=\int_0^t1_{\{Y^1(s)=Y^2(s)\}}\di s.
\ee
\et

As a corollary, we can strengthen Theorem \ref{T:stick2}  as follows.
\begin{corollary}\label{C:disweb}
	Let $\Di$ be a countable dense subset of $\R^2$. For each $z=(x,t)\in \Di$ and $i=1,2$, let $X^i_{\eps, z}$ be the unique path in
	$W^i_\eps$ starting from position $\eps \lfloor x\eps^{-1}\rfloor$ at time $\varsigma^2 \eps^2 \lfloor t/\varsigma^2 \eps^2\rfloor$. For $z, z'\in \Di$,
	let $Z_{\eps; z, z'}$ be the rescaled intersection time between $X^1_{\eps, z}$ and $X^2_{\eps, z'}$ defined as in \eqref{YZeps}. Then
	we have
	\be\label{SWZ}
	(S_{\eps}W^1_\eps,S_{\eps}W^2_\eps, (Z_{\eps; z, z'})_{z, z'\in \Di})
	\underset{\eps\down0}{\Longrightarrow} (\Wi^1,\Wi^2, (Z_{z,z'})_{z, z'\in \Di}),
	\ee
	where $Z_{z,z'}(u) = \int_{\max\{t, t'\}}^u 1_{\{Y^1_z(s) =Y^2_{z'}(s) \}} \di s$ is the intersection time between the unique path $Y^1_z\in \Wi^1(z)$
	and $Y^2_{z'}\in \Wi^2(z')$.
\end{corollary}
\noindent
This corollary follows readily from Theorem \ref{T:stick2} and \ref{T:stBMcon} and the fact that for $z, z'\in \Di$, $Z_{z, z'}$ is determined by $Y^1_z\in \Wi^1(z)$ and $Y^2_{z'}\in \Wi^2(z')$.
\medskip

To prove Theorem \ref{T:stBMcon}, we need a preliminary result for $S(t):=X^1(t)-X^2(t)$, which is a random walk
on $\Z$ with stickiness at $0$ and transition kernel
\be\label{tranS}
P(x,y)=\left\{\ba{ll}
\dis(1-\eps)1_{\{y=0\}}+\eps\bar{P}(x,y) \quad&\mbox{if }x=0,\\[5pt]
\bar{P}(x,y)\quad&\mbox{if }x\neq 0,
\ea\right.
\ee
where
\be\label{tranbS}
\bar P(x,x+y):=\sum_{z\in\Z}a(z+y)a(z)\qquad(x,y\in \Z).
\ee
For $t\in \N_0:=\N\cup\{0\}$, let $\bar P^t$ denote the $t$-th power of the transition
matrix $\bar P$, with $\bar P^0(x,y)=1_{\{x=y\}}$, and set $\ov G_t(x,y):=\sum_{s=0}^t\bar P^s(x,y)$. By
\cite[Prop.\ 28.8]{S01}, the limit
\be\label{At}
\bar{A}(x) :=\lim\limits_{t\rightarrow\infty}\big[\bar{G}_t(0,0)-\bar{G}_t(x,0)\big]<\infty
\quad(x\in\Z)
\ee
exists, and $\bar A$ is called the \emph{potential kernel} of the random walk with
transition kernel $\bar P$.

\bp[Compensator of potential kernel]\label{P:dismart}
Let $S$ denote the difference random walk with transition kernel $P$ as in
\eqref{tranS}. Then
\be\label{Amart}
\bar{A}(S(t))-\eps\sum_{s=0}^{t-1}1_{\{S(s)=0\}}
\ee
is a martingale. Moreover
\be\label{pkasym}
\lim\limits_{|x|\rightarrow\infty}\frac{\bar{A}(x)}{|x|}=\frac{1}{2\varsigma^2}.
\ee
\ep
\bpro
Formula (\ref{pkasym}) follows from \cite[Prop.\ 29.2]{S01}. To prove (\ref{Amart}), first denote $\bar{A}_t(x) := \bar{G}_t(0,0)-\bar{G}_t(x,0)$.
Then we have
\bc \label{traneq3}
\dis\bar{P}\bar{A}_t(x)&=& \dis\sum_{\scr{y\in \Z}}\bar{P}(x,y) \bar{A}_t(y) =
\dis\sum_{\scr{y\in \Z}}\bar{P}(x,y)\sum_{s=0}^{t}
\big(\bar{P}^s(0,0)-\bar{P}^s(y,0)\big)\\ [5pt]
&=&\dis\sum_{s=0}^{t}\big(\bar{P}^s(0,0)-\bar{P}^{s+1}(x,0)\big)%\\ [5pt]
=\dis\bar{A}_t(x)+\bar{P}^0(x,0)-\bar{P}^{t+1}(x,0).
\ec
Here $\lim_{t\to\infty}\bar{P}^{t+1}(x,0)=0$ by the null-recurrence of the
random walk with kernel $\bar P$, while
$\lim\limits_{t\rightarrow\infty}\bar{P}\bar{A}_t(x)=\bar{P}\bar{A}(x)$ by
\cite[Prop.\ 13.3]{S01}. Therefore letting $t\to\infty$ in \eqref{traneq3} gives
\be\label{traneq2}
\bar{P}\bar{A}(x)-\bar{A}(x)=1_{\{x=0\}}\qquad(x\in \Z).
\ee
We claim that for the transition kernel $P$ of the difference process $S$,
\bc\label{traneq1}
P\bar{A}(x)-\bar{A}(x)=\eps 1_{\{x=0\}}.
\ec
Indeed, when $x\neq0$, by the definition (\ref{tranS}), we have $P(x,y)=\bar{P}(x,y)$
for all $y\in\Z$. Hence (\ref*{traneq2}) is the same as (\ref*{traneq1}). When $x=0$,
recalling (\ref{tranS}) and using $\bar{A}(0)=0$, we have
\be
P\bar{A}(0)-\bar{A}(0)=(1-\eps)\sum_{y\in\Z}1_{\{y=0\}}\bar{A}(y)+\eps\bar{P}\bar{A}(0)=\eps\bar{P}\bar{A}(0)=\eps 1_{\{x=0\}},
\ee
where in the last equality we used (\ref{traneq2}).

In general, if $S(t)$ is a Markov chain with countable state space and
transition kernel $P$, then for all real functions $f$,
\be\label{stmart}
M(t):=f(S(t))-\sum_{s=0}^{t-1}\Big(Pf(S(s))-f(S(s))\Big)\qquad(t\geq 0)
\ee
is a martingale as long as each term in this expression is in $L^1$.
Applying this to the Markov chain $S$ and the function $\ov A$, and using
(\ref{traneq1}), we see that the process in (\ref{Amart})
is a martingale. Here, to check integrability, we have used that
$\E[|S(t)|]<\infty$ and hence by (\ref{pkasym}), we also have
$\E[\bar{A}(S(t))]<\infty$.
\epro

To prove Theorem \ref{T:stBMcon}, we also need the following facts. Recall that for continuous local martingales $M_1$ and $M_2$, there exists a unique process $\li M_1,M_2\re$ of bounded variation such that
\be
M_1(t)M_2(t)-\li M_1,M_2\re(t)
\ee
is a local martingale (see \cite[Def.\ 1.5.18]{KS91}). Moreover, by \cite[Thm.\ 3.7.1]{KS91} one can define a local time $L(t, x)$ for any real square integrable continuous semimartingale $M$ such that $L(t, x)$ is a measurable function, continuous and nondecreasing in $t$, and for any measurable function $f:\R\to\R$,
\be\label{localtime}
\int_{0}^{t}f(M(s))\di\li M,M\re(s)=\int_{\R}f(x)L(t, x)\di x.
\ee
\medskip

\bpro[of Theorem \ref{T:stBMcon}]
Since $X^1$ and $X^2$ are both random walks with the same jump distribution $a(\cdot)$, by Donsker's invariance principle (see, e.g., \cite[Theorem~8.2]{B99}), $Y_\eps^1$ and $Y_\eps^2$ converge to standard Brownian motions starting from $y^1$ and $y^2$, respectively. This also implies the tightness of $(Y^1_\eps,Y^2_\eps)$. Consequently, to show that $(Y^1_\eps,Y^2_\eps)$ converges to a pair of sticky Brownian motions, it suffices to show that for any subsequential weak limit $(Y^1, Y^2)$,
the processes in (i) and (ii) of (\ref{HWmart}) are martingales. Along the way, we will see that $Z_\eps$, as a function of $(Y^1_\eps,Y^2_\eps)$, also converges weakly to $Z$ defined from $(Y^1, Y^2)$ as in \eqref{Zt}.

Note that the transition kernel of $\big(X^1(t),X^2(t)\big)_{t\in\N}$ is given by
\be
P^{(2)}(x_1,x_2;x_1+z_1,x_2+z_2)=\left\{\ba{ll}
(1-\eps)1_{\{z_1=z_2\}}a(z_1)+\eps a(z_1)a(z_2)&\quad\mbox{if }x_1=x_2,\\[5pt]
a(z_1)a(z_2)&\quad\mbox{if }x_1\neq x_2. \ea
\right.
\ee
Therefore by standard theory (see (\ref{stmart})), the discrete time process
\be\label{covmart1}
X^1(t)X^2(t)-\sum_{k=0}^{t-1}\Ga\big(X^1(k),X^2(k)\big) \qquad (t\geq 0)
\ee
is a martingale, where
\bc
\Ga(x_1,x_2)&=&\dis\sum_{z_1,z_2\in\Z}P^{(2)}(x_1,x_2;x_1+z_1,x_2+z_2)(x_1+z_1)(x_2+z_2)-x_1x_2\\[5pt]
&=&(1-\eps)\varsigma^21_{\{x_1=x_2\}}.
\ec
In terms of $Y^1_\eps$ and $Y^2_\eps$, whose trajectories are piecewise constant, this implies that
\be\label{covmart2}
Y^1_\eps(t)Y^2_\eps(t)-(1-\eps)\int_{0}^{\lfloor t\rfloor_{\eps}}1_{\{Y^1_\eps(s)=Y^2_\eps(s)\}}\di s \qquad (t\geq 0)
\ee
is a martingale, where $\lfloor t\rfloor_{\eps}:=\varsigma^2\eps^2 \lfloor t/\varsigma^2 \eps^2\rfloor$.

On the other hand, since the transition kernel of the difference process $X^1(t)-X^2(t)$ is given by (\ref{tranS}), by Proposition~\ref{P:dismart},
\be\label{diffmart1}
\eps\bar{A}\big(\eps^{-1} [Y^1_\eps(t)-Y^2_\eps(t)]\big)-  \frac{1}{\varsigma^2} \int_{0}^{\lfloor t\rfloor_{\eps}}1_{\{Y^1_\eps(s)=Y^2_\eps(s)\}}\di s
\ee
is also a martingale. It will be more convenient to work with
\be
\widetilde{Z}_\eps(t):=\int_{0}^{t}1_{\{Y^1_\eps(s)=Y^2_\eps(s)\}}\di s,
\ee
for which we note that $Z_\eps(t)=\widetilde{Z}_\eps(\lfloor t\rfloor_\eps)$ for all $t\geq0$.

To prove the convergence of the discrete martingale (\ref{covmart2}) to a continuous one, we adopt an argument in \cite[Prop.\ A.6]{SSS14} (which is
adapted from \cite{HW09a}). First of all, recall that the laws of the processes $(Y^1_\eps,Y^2_\eps)$ on $ D([0,\infty), \R^2)$ are tight. Secondly, since the slope of $\widetilde Z_\eps$ is between zero and one, by the Arzel\`a-Ascoli Theorem, tightness holds also for $(\widetilde Z_\eps)_{\eps>0}$. As a result, by going to a subsequence if necessary, we may assume that the joint processes
$(Y^1_\eps,Y^2_\eps,\widetilde Z_\eps)$ converge weakly in law as $\eps\down0$,
and by Skorohod's representation Theorem (see e.g. \cite[Thm.\ 6.7]{B99}), we can couple $(Y^1_\eps,Y^2_\eps,\widetilde Z_\eps)$ along any sequence $\eps\downarrow 0$ such that the convergence is almost sure. Moreover, since for all $t\geq0$, $(Y^1_\eps(t), Y^2_\eps(t), \widetilde Z_\eps(t))_{\eps>0}$ are uniformly integrable, their convergence is also in $L^1$. Let $(Y^1,Y^2,Z_\ast)$ denote the limiting process. Then, taking the limit
in (\ref{covmart2}) and using the Markov property of $(Y^1_\eps, Y^2_\eps)$, we see that
\be
Y^1(t)Y^2(t)-Z^\ast(t)
\ee
is a martingale, and hence
\be\label{qvar}
\li Y^1,Y^2\re(t)=Z^\ast(t)=\lim_{\eps\down0}
\int_0^t1_{\{Y^1_\eps(s)=Y^2_\eps(s)\}}\di s
\qquad\forall t\geq 0\quad{\rm a.s.}
\ee

On one hand, since for a given $t>0$, the function $(w_1, w_2) \mapsto\int_0^t1_{\{w_1(s)=w_2(s)\}}\di
s$ is upper semicontinuous with respect to the topology on $ D([0,\infty), \R^2)$,
formula (\ref{qvar}) implies that
\be\label{Zleq}
\li Y^1,Y^2\re(t)\leq\int_0^t1_{\{Y^1(s)=Y^2(s)\}}\di s
\qquad(t\geq 0).
\ee
On the other hand, since $\li Y^1,Y^2\re(t)$ is nondecreasing by \eqref{qvar}, almost surely $\di \li Y^1,Y^2\re(t)$ defines a nonnegative measure on $\half$, which implies that
\be\label{Zgeq}
\li Y^1,Y^2\re(t)\geq\int_0^t1_{\{Y^1(s)=Y^2(s)\}}\di\li Y^1,Y^2\re(s)\qquad(t\geq 0).
\ee
To prove the opposite inequality to the one in (\ref{Zleq}), let $L(t, x)$ be the local time of the semimartingale $Y^1-Y^2$.
Applying (\ref{localtime}) to the function $f=1_{\{0\}}$, we find that
\be\label{octi}
\int_0^t1_{\{Y^1(s)=Y^2(s)\}}\di\li Y^1-Y^2,Y^1-Y^2\re(s)= \int_\R 1_{\{0\}}(x) L(t, x)\di x = 0.
\ee
Since $Y^1,Y^2$ are standard Brownian motions, we have
\bc
\dis\li Y^1-Y^2,Y^1-Y^2\re(t)
&=&\dis\li Y^1,Y^1\re(t)+\li Y^2,Y^2\re(t)-2\li Y^1,Y^2\re(t)\\[5pt]
&=&\dis 2t-2\li Y^1,Y^2\re(t).
\ec
Inserting this into (\ref{octi}) yields
\be\label{difplace}
\int_0^t1_{\{Y^1(s)=Y^2(s)\}}\di\li Y^1,Y^2\re(s) = \int_0^t1_{\{Y^1(s)=Y^2(s)\}}\di s.
\ee
Combining the above equality with (\ref{Zgeq}) and (\ref{Zleq}) shows that
\be\label{Zeq}
Z^\ast(t) = \li Y^1,Y^2\re(t)=
\int_0^t1_{\{Y^1(s)=Y^2(s)\}}\di s=Z(t)
\qquad\forall t\geq 0\quad{\rm a.s.}
\ee
Therefore, $Y^1(t)Y^2(t)-\int_0^t1_{\{Y^1(s)=Y^2(s)\}}\di s$ is a martingale. Moreover, it is easy to see that the almost sure convergence of $\widetilde Z_\eps$ to $Z$ implies the almost sure convergence of $Z_\eps$ to $Z$.

To check that the process in (ii) of (\ref{HWmart}) is a martingale, we take the limit $\eps\downarrow 0$ in (\ref{diffmart1}).
Since $\bar{P}$ defined in (\ref{tranbS}) has second moment $2\varsigma^2$, by Proposition~\ref{P:dismart}, for any $\de>0$, there exists $N>0$ such that for all integers $|x|>N$,
\bc
\De(x):=\big|2\varsigma^2\bar{A}(x)-|x|\big|<\de|x|.
\ec
Recall that $S=X^1-X^2$ and denote $C_N:=\sup_{|x|\leq N}\big|2\varsigma^2\bar{A}(x)-|x|\big|$. We see that
\bc
\eps\De\big(S(t\varsigma^{-2}\eps^{-2})\big)
&=&\eps\De\big(S(t \varsigma^{-2}\eps^{-2})\big)1_{\{S(t\varsigma^{-2}\eps^{-2})>N\}} + \eps\De(S(t\varsigma^{-2}\eps^{-2}))1_{\{S(t\varsigma^{-2}\eps^{-2})\leq N\}}\\ [5pt]
&\leq& \eps \de|S(t \varsigma^{-2}\eps^{-2})|+\eps C_N = \de  |Y^1_{\eps}(t)-Y^2_{\eps}(t)|+\eps C_N.
\ec
Sending $\eps\downarrow 0$ first and then $\de\downarrow 0$, and using the almost sure convergence of $Y^1_{\eps}(t)-Y^2_{\eps}(t)$ to $Y^1(t)-Y^2(t)$, we obtain that
\be\label{err1}
\limsup_{\eps\down0}\eps\De\big(\eps^{-1} [Y^1_{\eps}(t)-Y^2_{\eps}(t)]\big)=0\qquad\forall t\geq 0\quad{\rm a.s.}
\ee
Therefore, almost surely, for all $t\geq0$,
\be\label{kerlim}
\lim\limits_{\eps\down0}\eps\bar{A}\big(\eps^{-1}[Y^1_{\eps}(t)-Y^2_{\eps}(t)]\big)=\lim\limits_{\eps\down0}\frac{1}{2\varsigma^2}|Y^1_{\eps}(t)-Y^2_{\eps}(t)|=\frac{1}{2\varsigma^2}|Y^1(t)-Y^2(t)|.
\ee
Recalling (\ref{qvar}) and (\ref{Zeq}), the almost sure  limit (also $L^1$ limit) of (\ref{diffmart1}) is
\be
\frac{1}{2\varsigma^2}|Y^1(t)-Y^2(t)|- \frac{1}{\varsigma^2} \int_0^t1_{\{Y^1(s)=Y^2(s)\}}\di s.
\ee
Together with the Markov property of $(Y^1_\eps, Y^2_\eps)$, it follows that
\be\label{diffmart2}
|Y^1(t)-Y^2(t)|- 2 \int_0^t1_{\{Y^1(s)=Y^2(s)\}}\di s %\frac{2}{m}\int_0^t1_{\{Y^1(s)=Y^2(s)\}}\di s
\ee
is a martingale.

Therefore, the process $(Y_\eps^1,Y_\eps^2,Z_\eps)$ converges in distribution on
$D([0,\infty), \R^3)$ to $(Y^1,Y^2,Z)$, where $(Y^1,Y^2)$ is a pair of sticky Brownian motions with parameter $1$ and
starting point $(y^1,y^2)$ as characterized in Proposition \ref{P:stickMP}, and $Z$ is defined from $(Y^1, Y^2)$ as in \eqref{Zt}.
\epro

\subsection{Proof of the lower bound}\label{S:lowbd}

We are now ready to complete the lower bound on any subsequential weak limit of $S_{\eps}N_\eps$.
\bt[Lower bound]\label{T:lower}
For $\eps>0$, let $N_\eps$ be the discrete net with branching probability $\eps$ and jump distribution $a(\cdot)$ satisfying \eqref{passum}.
Then for any subsequential weak limit $\Ni^\ast$ of $S_{\eps}N_\eps$ as $\eps\down0$, there exists a standard Brownian net $\Ni$ such that almost surely,
\be
\Ni\subseteq\Ni^\ast.
\ee
\et
\bpro
Let $(W_\eps^1, W_\eps^2)$ be the pair of discrete sticky webs defined before \eqref{a2}, which is contained in the the discrete net $N_\eps$.
Thanks to Corollary \ref{C:disweb},  if $\Ni^\ast$ is a subsequential limit of $S_{\eps}N_\eps$ as $\eps\downarrow 0$, then by going to a subsequence if necessary and by Skorohod's representation theorem, we can find a coupling such that
\be\label{WWSkor}
(S_{\eps}W_\eps^1,S_{\eps}W_\eps^2,S_{\eps}N_\eps, (Z_{\eps; z, z'})_{z, z'\in \Di}) \underset{\eps\down0}{\longrightarrow} (\Wi^1,\Wi^2,\Ni^\ast, (Z_{z,z'})_{z,z'\in \Di})\quad{\rm a.s.},
\ee
where $\Wi^1\cup\Wi^2\subseteq\Ni^\ast$, $\Di$ is a deterministic countable subset of $\R^2$, and $Z_{z, z'}$ is the intersection time process
between $Y^1_z\in \Wi^1(z)$ and $Y^2_{z'}\in \Wi^2(z')$ while $Z_{\eps; z, z'}$ is the discrete counterpart. By Theorem~\ref{T:Nhop}, we can construct a standard Brownian net via hopping between paths in $\Wi^1$ and $\Wi^2$, i.e.,
\be\label{lowbd2}
\Ni=\ov{\Hi_{\rm int}(\Wi^1\cup\Wi^2)}.
\ee
To show $\Ni\subseteq\Ni^\ast$, it suffices to show that
\be\label{lowbd3}
\Hi_{\rm int}(\Wi^1\cup\Wi^2)\subseteq\Ni^\ast,
\ee
since $\Ni^\ast$ is a compact subset of the path space $\Pi$. Let $\Di$ be a countable dense subset of $\R^2$. For any path $\pi\in\Wi^i$, we can find
a sequence of $\pi_n \in \Wi^i(\Di)$ such that $\pi_n\to \pi$ and the time of coalescence between $\pi_n$ and $\pi$ converges to the starting time of $\pi$, see \cite[Lemma 3.4]{SS08}. Therefore to prove \eqref{lowbd3}, it suffices to to show $\Hi_{\rm int}(\Wi^1(\Di)\cup\Wi^2(\Di))\subseteq\Ni^\ast$. In particular, by induction, it suffices
to show that for any $Y^1_z\in\Wi^1(z)$ and $Y^2_{z'}\in\Wi^2(z')$ with $z=(x,t), z'=(x', t')\in \Di$, if $u>\max\{t, t'\}$ is an intersection time of $Y^1$ and $Y^2$, then the concatenated path with graph
\be
Y:=\{(Y^1_z(s),s):s\in[t, u]\}\cup\{(Y^2_{z'}(s),s):s\in[u,\infty]\} \in \Ni^\ast.
\ee
By \eqref{WWSkor}, we can find $Y^1_{\eps, z} \in S_{\eps}W_\eps^1$ and $Y^2_{\eps, z'}\in S_{\eps}W_\eps^2$ such that
$$
(Y^1_{\eps, z}, Y^2_{\eps, z'}, Z_{\eps; z, z'}) \to (Y^1_z, Y^2_{z'}, Z_{z, z'}),
$$
where $Z_{\eps; z, z'}$ and $Z_{z, z'}$ are the intersection time processes for $(Y^1_{\eps, z}, Y^2_{\eps, z'})$ and $(Y^1_z, Y^2_{z'})$, respectively.
Note that $Y^1_z(s)-Y^2_{z'}(s)$ is a sticky Brownian motion with stickiness at $0$, which can be obtained from a Brownian motion via a random time change that changes its local time at the origin into real time, see e.g.\ \cite[Section 2.1]{SS08} and \cite[Section 2.1]{Y16}. This ensures that any intersection time $u$ of $Y^1_z$ and $Y^2_{z'}$ lies in the support of the measure $\di Z_{z, z'}(s)$. The convergence $Z_{\eps; z, z'} \to Z_{z, z'}$ then implies that we can find intersection times $u_\eps$ of $Y^1_{\eps, z}$ and $Y^2_{\eps, z'}$ such that $u_\eps \to u$. Therefore we can concatenate $Y^1_{\eps, z}$ and $Y^2_{\eps, z'}$ at time $u_\eps$, which gives a path in $S_\eps N_\eps$ that clearly converges to $Y$.
This shows that $Y$ is an element of $\Ni^\ast$ and completes the proof.
\epro

\section{Auxiliary Bernoulli net}\label{S:bern}

In this section, we introduce an auxiliary system of branching-coalescing random walks which we call the {\em Bernoulli net} $\wt N$.
It will play a crucial role in the upper bound proof in Section \ref{S:upproof} and in the tightness proof in Section~\ref{S:tight}. We refer to
Section \ref{S:crossingstrat} for an outline of how the Bernoulli net fits into our proof strategy. This section is organised as follows. In Subsection \ref{S:coup}, we show that $\wt N$ can be coupled to the discrete net $N_\eps$ such that $N_\eps \subset \wt N$ almost surely. In Subsection \ref{S:braco}, we introduce the branching-coalescing point set $\tilde \xi$ generated by $\wt N$, and its dual process $\tilde \phi$ which is a discrete time biased voter model. In Subsection \ref{S:backbone}, we show that $\tilde \xi$ admits an invariant Bernoulli product law, which
%which is also generated by the {\em backbone} of $\wt N$, i.e., set of paths in $\wt N$ starting at time $-\infty$.
converges to the invariant law of the branching-coalescing point set $\xi$ generated by the Brownian net $\Ni$. This fact is the key to the upper bound proof in Section \ref{S:upproof}. Lastly in Subsection \ref{S:denest}, we give an upper bound on the density of $\tilde \xi$ with maximal initial configuration, which is comparable to the density of coalescing random walks up to the diffusive time scale.

\subsection{Bernoulli net and coupling with discrete net}\label{S:coup}

Recall that the discrete net $N_\eps$ was defined before \eqref{a2} by specifying the set of outgoing arrows from each $(x, t)\in \Z^2$. To define a
more general branching-coalescing system, we can specify a probability law $\mu$ on $\Pc_{\rm fin,+}(\Z)$, the space of non-empty finite subsets of $\Z$.
We then assign i.i.d.\ $\Pc_{\rm fin,+}(\Z)$-valued random variables $\Ai(x,t)$ with common law $\mu$ to each $(x,t)\in \Z^2$. For every $y\in \Ai(x,t)$,
we then draw an arrow from $(x, t)$ to $(x+y, t+1)$. The set of arrows then defines a branching-coalescing system similar to the discrete net $N_\eps$,
with $\mu$ specifying the offspring distribution.

We now introduce the Bernoulli net $\wt N$, which corresponds to a specific class of $\mu$. Given $\psi:\Z\to[0,1)$ with $\sum_{x\in\Z}\psi(x)<\infty$, let $(V_x)_{x\in\Z}$ be independent Bernoulli random variables with parameters $(\psi(x))_{x\in \Z}$ respectively, and let ${\rm Supp}(V) :=\{x\in \Z: V_x=1\}$. We then choose the offspring distribution $\mu$ such that
\be\label{brho}
\mu(A) =\P\big({\rm Supp(V)} =A\,\big| {\rm Supp}(V)\neq \emptyset\big)
\qquad\big(A\in\Pc_{\rm fin,+}(\Z)\big).
\ee
We call the associated branching-coalescing system $\wN_\psi$ the \emph{Bernoulli net} with \emph{kernel} $\psi$. We say that $\psi$ is
\emph{irreducible} if paths in $\wN_\psi(0,0)$ (the set of paths in $\wN_\psi$ starting at $(0,0)$) can reach any point in $\Z$ with positive probability;
and we say $\psi$ is \emph{aperiodic} if the set of times at which $0$ can be reached by paths in $\wN_\psi(0, 0)$ with positive probability has greatest
common divisor $1$.

We now show how the discrete net $N_\eps$ with jump distribution $a(\cdot)$ and branching probability $\eps$, defined as in \eqref{a2}, can be coupled to
a  Bernoulli net $\wN_\psi$ with kernel $\psi$, such that $N_\eps\sub\wN_\psi$ almost surely and the difference between them vanishes as $\eps\downarrow 0$.

\bl[Coupling of Bernoulli and discrete nets]\label{L:Bernbin}
Let $N_\eps$ be the discrete net with jump distribution $a(\cdot)$ and branching probability $\eps>0$. Let $r>0$ be a solution of the equation
\be\label{epsr}
\eps =1 - \frac{r}{e^r-1}.
\ee
Define
\be\label{ala}
\psi_\eps(x):= 1- e^{-ra(x)} \quad(x\in\Z).
\ee
Then the Bernoulli net $\wN_{\psi_\eps}$ with kernel $\psi_\eps$ can be coupled to $N_\eps$ such that $N_\eps\sub\wN_{\psi_\eps}$ almost surely.
\el
\bpro
It suffices to couple the configuration of arrows originating from $(0,0)$ in $N_\eps$ and $\wN_{\psi_\eps}$. The basic idea is to construct
the arrows in $\wN_{\psi_\eps}$ through a Poisson point process. Let $(M_x)_{x\in \Z}$ be independent Poisson random variables with mean
$ra(x)$. We then draw $M_x$ arrows from $(0,0)$ to $(x,1)$, where multiple arrows coalesce with each other at $(x,1)$. Note that this is
equivalent to assigning a single arrow from $(0,0)$ to $(x,1)$ with probability $\psi_\eps(x)=1-e^{-ra(x)}$. To construct the configuration
of arrows in $\wN_{\psi_\eps}$, we are simply conditioning on the event $M:=\sum_x M_x \geq 1$. With our choice of $r$, we have
$$
\P(M=1 |M\geq 1) = \frac{r e^{-r}}{1-e^{-r}} = \frac{r}{e^r-1} = 1-\eps.
$$
Restricted to the event $\{M=1\}$, it is clear that the single arrow from $(0,0)$ is drawn according to distribution $a(\cdot)$, which
coincides with the definition of $a^{(2)}$ in \eqref{a2}. On the event $\{M\geq 2\}$, we draw $M$ arrows independently, each following
the distribution $a(\cdot)$. Retaining only the first 2 arrows already reproduces the arrow configuration in the definition of $a^{(2)}$
in \eqref{a2}. This establishes the desired coupling between $N_\eps$ and $\wN_{\psi_\eps}$.
\epro

\begin{remark}\label{R:reps}
	We are interested in applying Lemma \ref{L:Bernbin} as $\eps\downarrow 0$. By Taylor expansion around $r=0$ in the r.h.s.\ of \eqref{epsr},
	it is easy to see that for $\eps>0$ sufficiently small, there is a unique solution $r(\eps)$ to \eqref{epsr}, and $r(\eps)= (2+o(1))\eps$
	as $\eps\downarrow 0$,  i.e., $\lim_{\eps\downarrow 0} r(\eps)/\eps = 2$.
\end{remark}

\subsection{Branching-coalescing point sets and their duals}\label{S:braco}
Let $\wN :=\wN_\psi$ be a Bernoulli net with kernel $\psi$. We now introduce two associated Markov chains. One is the so-called
branching-coalescing point set $\tilde \xi_t$, the other is its dual $\tilde \phi_t$ which can be interpreted as a biased voter model (see
e.g.\ \cite{SSY19, SSY21}).

Let $\Pc(\Z)$ denote the collection of subsets of $\Z$. Given a Bernoulli net $\wN$, for each $s,t\in\Z$ with $s\leq t$, we define two random maps $\Xi_{s,t}:\Pc(\Z)\to\Pc(\Z)$ and $\Phi_{t, s}:\Pc(\Z)\to\Pc(\Z)$ by
\be\left.\ba{r@{\,}c@{\,}l}\label{XiPhi}
\dis\Xi_{s,t}(A)&:=&\dis
\big\{y\in\Z:\exists x\in A,\ \pi\in \wN(x,s)\mbox{ s.t.\ }\pi(t)=y\big\},\\[5pt]
\dis\Phi_{t, s}(B)&:=&\dis
\big\{x\in\Z:\exists y\in B,\ \pi\in \wN(x,s)\mbox{ s.t.\ }\pi(t)=y\big\},
\ea\right.\quad\big(A,B\in\Pc(\Z)\big),
\ee
where $\wN(x, s)$ is the set of paths in $\wN$ starting from $(x,s)$. In words, $\Xi_{s, t}(A)\subset \Z$ is the set of positions generated at time $t$ by paths in $\wN$ that start from $A$ at time $s$, while
$\Phi_{t, s}(B)$ is the set of positions at time $s$ that can lead to some point in $B$ by following paths in $\wN$.

Fix two times $s, u\in \Z$. We can then define two $\Pc(\Z)$-valued Markov chains
\be\label{xiphi}
\tilde \xi_t:=\Xi_{s,s+t}(A)\quand \tilde\phi_t:=\Phi_{u,u-t}(B)\qquad(t\geq 0),
\ee
where $\tilde\xi$ evolves forward in time and $\tilde \phi$ evolves backward in time. We call $\tilde \xi$ the branching-coalescing
point set generated by $\wN$, and $\tilde \phi$ its dual. More precisely, they are dual in the sense that
\be\label{disdu}
\P\big(\tilde \xi_t\cap B\neq\emptyset \big| \tilde \xi_0=A\big)=\P\big(A\cap \tilde \phi_t\neq\emptyset \big| \tilde \phi_0=B\big),
\qquad\big(A,B\in\Pc(\Z),\ t\geq 0),
\ee
which holds since
$$
\{\Xi_{0,t}(A)\cap B\neq\emptyset\}
=\big\{\exists x\in A,\ y\in B,\ \pi\in \wN(x,0)\mbox{ s.t.\ }\pi(t)=y\big\}
=\{A\cap\Phi_{t, 0}(B)\neq\emptyset\}.
$$

It will be convenient to define the Markov chain $\tilde \xi$ also at half-integer
times in order to keep track of the arrows that are used in each time step.
To this aim, for each $s,t\in\Z$ with $s\leq t$, we define
$\Xi_{s,t+1/2}:\Pc(\Z)\to\Pc(\Z^2)$ by
\be\label{Xihalf}
\Xi_{s,t+1/2}(A):=
\big\{\big(\pi(t),\pi(t+1)\big) \in \Z^2 :\pi\in \wN(x,s)\mbox{ for some }x\in A\big\},
\ee
and we define $\xi_{t+1/2}$ as in (\ref{xiphi}). Then
$(\xi_0,\xi_{1/2},\xi_1,\ldots)$ is a Markov chain taking values alternatively
in $\Pc(\Z)$ and $\Pc(\Z^2)$, where in fact $\xi_{t+1}$ is a deterministic
function of $\xi_{t+1/2}$.
\medskip

In the continuum setting, if $\Ni$ is the standard Brownian net, then for each real $s\leq t$, we
can similarly define a random map $\Xi^\Ni_{s,t}$ acting on closed subsets of the real line by
\be
\Xi_{s,t}(A):=\big\{\pi(t):\pi\in\Ni(x,s)\mbox{ for some }x\in A\big\}.
\ee
Setting $\xi^\Ni_t:=\Xi^\Ni_{s,s+t}(A)$ $(t\geq 0)$ then defines a Markov process
taking values in the space of closed subsets of the real line. This process is known as
the \emph{branching-coalescing point set} of the Brownian net $\Ni$ \cite[Thm~1.11]{SS08}.

\subsection{Product invariant law and the backbone}\label{S:backbone}

For the standard Brownian net $\Ni$, the branching-coalescing point set $\xi^\Ni_t$ is known to have a unique
invariant law given by the law of the Poisson point process on $\R$ with intensity $2$, see \cite[Prop.~1.15]{SS08}.
Furthermore, the stationary process can be constructed by considering $\Ni(*, -\infty)$, the set of paths
in $\Ni$ starting at $t=-\infty$, which is called the  \emph{backbone} of $\Ni$. Similar results hold for the
discrete net $N_\eps$ with jump distribution $a(\cdot)$ supported on $\{\pm 1\}$, see \cite[Prop.~1.14]{SS08}.
In this subsection, we show that analogous results also hold for the Bernoulli nets $\wN$, and we identify explicitly
the product invariant law of the branching-coalescing point set $\tilde \xi$ generated by $\wN$.

In what follows, a {\em Bernoulli subset of $\Z$ with intensity $\alpha(\cdot)$} refers to a random subset of $\Lambda \subset\Z$,
where independently each $x\in \Z$ belongs to $\Lambda$ with probability $\alpha(x)$. A Bernoulli subset of $\Z^2$
with intensity $\alpha(\cdot)$ is defined similarly. We will write $\N_0/2:=\{x/2:x\in\N_0\}$ where $\N_0:=\{0\}\cup\N$.

\bl[Product invariant law]\label{L:prodinv}
Let $\wt N$ be a Bernoulli net with kernel $\psi: \Z \to [0,1)$. Let $\Lambda$ be a Bernoulli subset
of $\Z$ with constant intensity
\be\label{rhopsi}
\rho:=1-\prod_{x\in \Z} (1-\psi(x)).
\ee
Let $(\tilde \xi_t)_{t\in \N_0/2}$ be the branching-coalescing point set defined as in \eqref{xiphi} and \eqref{Xihalf} from $\wN$ with initial condition $\tilde \xi_0=\Lambda$. Then for each $t\in\N$,
\begin{enumerate}
	\item $\tilde \xi_t$ is a Bernoulli subset of $\Z$ with intensity $\rho$;
	
	\item $\tilde\xi_{t+1/2}$ is a Bernoulli subset of $\Z^2$ with intensity $\alpha(x,y)=\psi(y-x)$.
\end{enumerate}
\el
\bpro
Let $\xi_{1/2}$ be a Bernoulli subset of $\Z^2$ with intensity $\alpha(x,y)=\psi(y-x)$.
Define random sets $\xi_0, \xi_1 \subset \Z$ by
\bc
\dis\xi_0&:=&\dis\{x\in\Z:(x,y)\in\xi_{1/2}\mbox{ for some }y\in\Z\},\\[5pt]
\dis\xi_1&:=&\dis\{y\in\Z:(x,y)\in\xi_{1/2}\mbox{ for some }x\in\Z\}.
\ec
Then it is easy to see that both $\xi_0$ and $\xi_1$ are Bernoulli subset of
$\Z$ with constant intensity $\rho$. Moreover, for $x\in\xi_0$, if we define
\be
\Ai(x,0):=\{y\in\Z:(x,x+y)\in\xi_{1/2}\},
\ee
then we note that conditioned on $\xi_0$, the random variables $(\Ai(x,0))_{x\in\xi_0}$ are
i.i.d.\ with common law $\mu$ defined from $\psi$ as in (\ref{brho}). In view of this, the
conditional law of $\xi_{1/2}$ given $\xi_0$ is the same as that for the branching-coalescing
point set $\tilde \xi$ defined from the Bernoulli net $\wN$ with kernel $\psi$, and also
$\xi_1$ is a function of $\xi_{1/2}$ in exactly the same way as $\tilde \xi$. The statements
(i) and (ii) then follow immediately.
\epro

Next we show that the backbone of the Bernoulli net $\wN$ with kernel $\psi$, defined by
\be\label{tildebb}
\tilde \xi^*_t := \{ \pi(t): \pi \in \wN(*, -\infty)\},
\ee
is a stationary process with the law of $\tilde \xi^*_t$ for each $t\in\Z$ given by the product invariant law identified in
Lemma \ref{L:prodinv}. In \cite{SS08}, such a result has been proved for the Brownian net $\Ni$ as
well as for the nearest neighbour discrete net. The proof for the latter, see \cite[Prop.~1.14~(i)]{SS08},
makes essential use of the nearest neighbour property. Here we argue differently by identifying a suitable martingale
for the process $\tilde \phi$ dual to the branching-coalescing point set $\tilde \xi$. We start with a preliminary lemma
for $\tilde \phi$.

\bl[Positive fluctuations]\label{L:fluct}
Let $\wt N$ be a Bernoulli net with kernel $\psi: \Z \to [0,1)$, which is irreducible and aperiodic. Let $\tilde\phi=(\tilde\phi_t)_{t\geq
	0}$ be the dual Markov chain defined in \eqref{xiphi} with initial state $\tilde\phi_0\subset \Z$. Then for all initial states $\tilde \phi_0\subset \Z$ with $0<|\tilde \phi_0|<\infty$, we have
\be\label{fluct}
\P\big(|\tilde\phi_1|\neq|\tilde\phi_0|\big)\geq \delta(\psi) := \frac{1}{\rho} \min\Big\{ \psi(y_1), \psi(y_2) \prod_{w\leq y_1} (1-\psi(w)) \Big\} >0,
\ee
where $\rho$ is defined as in \eqref{rhopsi}, and $y_1<y_2\in\Z$ are any two points with $\psi(y_1), \psi(y_2)\in (0,1)$.

In particular, if $\wN_{\psi_\eps}$ is the Bernoulli net coupled to the discrete net $N_\eps$ with jump distribution $a(\cdot)$ and branching probability $\eps$ as in Lemma \ref{L:Bernbin}, then $\liminf_{\eps\downarrow 0} \delta(\psi_\eps)>0$.
\el
\bpro
By assumption, there exist $y_1 < y_2\in \Z$ with $\psi(y_1), \psi(y_2) \in (0,1)$. Denote $z:=\max \tilde \phi_0$, and let $x=z-y_1$.
Then by the definition of $\tilde\phi$ in \eqref{XiPhi} and \eqref{xiphi} and the definition of $\wN$
with kernel $\psi$ in \eqref{brho}, we have the following bounds:
\be\label{xtphi1}
\P(x\in \tilde \phi_1) \geq \frac{\psi(y_1)}{\rho} \qquad \mbox{and} \qquad \P(x\notin \tilde \phi_1) \geq \frac{1}{\rho} \cdot \psi(y_2) \prod_{w\leq y_1} (1-\psi(w)),
\ee
where the second bound is based on the event that in the random map $\Phi_{1, 0}(\tilde \phi_0)$ defined as in \eqref{XiPhi}, there is no arrow in $\wN$ from $x$ to any $z'\leq z$, and there is an arrow from $x$ to $x+y_2$. Note that conditioned on $\tilde\phi_1\backslash\{x\}$, either $\{ x\in \tilde \phi_1\} \subset \{|\tilde\phi_1|\neq|\tilde\phi_0|\}$ or $\{ x\notin \tilde \phi_1\} \subset \{|\tilde\phi_1|\neq|\tilde\phi_0|\}$. Furthermore, the event $\{x\in \tilde \phi_1\}$ is independent of $\tilde\phi_1\backslash\{x\}$. The claim \eqref{fluct} then follows readily from the bounds in \eqref{xtphi1}.

When $\psi=\psi_\eps = 1 - e^{-ra(x)}$ as in \eqref{ala}, where $r=r(\eps) = (2+o(1))\eps$ as $\eps\downarrow 0$ by Remark \ref{R:reps}, we have
\begin{equation} \label{rhoreps}
	\rho=1 - e^{-r}=(2+o(1))\eps.
\end{equation}
If we choose $y_1<y_2$ with $a(y_1), a(y_2)>0$, then we have
$$
\lim_{\eps\downarrow 0} \psi_\eps(y_1)/\rho =a(y_1) \qquad \mbox{and} \qquad
\lim_{\eps\downarrow 0} \psi_\eps(y_2) \prod_{w\leq y_1} (1-\psi_\eps(y_1))/\rho =a(y_2).
$$
Therefore $\liminf_{\eps\downarrow 0} \delta(\psi_\eps)\geq \min\{a(y_1), a(y_2)\}>0$.
\epro

\bp[The backbone and the product invariant law]\label{P:disdens}
Let $\wt N$ be a Bernoulli net with kernel $\psi$, which is irreducible and aperiodic. Let $(\tilde \xi^*_t)_{t\in\Z}$ be the
backbone of $\wN$ defined in \eqref{tildebb}. Then for each $t\in\Z$, $\tilde \xi^*_t$ is a Bernoulli subset of $\Z$ with constant intensity
$\rho$, defined as in Lemma \ref{L:prodinv}.
\ep
\bpro
Let $\tilde\phi=(\tilde\phi_t)_{t\geq 0}$ be the Markov chain defined from $\wN$ as in
(\ref{xiphi}), starting from a deterministic initial state $\tilde\phi_0$ with
$|\tilde\phi_0|<\infty$. Let $\rho$ be defined from $\psi$ as in \eqref{rhopsi}.
We claim that the process
\be\label{1rhophi}
M_t:=(1-\rho)^{|\tilde\phi_t|}\qquad(t\in\N)
\ee
is a martingale. Indeed, the duality between the branching-coalescing point set
$\tilde \xi$ and $\tilde \phi$ in \eqref{disdu} gives
\be
\E[M_1]=\E\big[(1-\rho)^{|\tilde\phi_1|}\big]=\P\big(\tilde\xi_0\cap\tilde\phi_1=\emptyset\big)
=\P\big(\tilde\xi_1\cap\tilde\phi_0=\emptyset\big)=(1-\rho)^{|\tilde\phi_0|}=M_0,
\ee
where we have chosen $\tilde \xi_0$ according to the product invariant law identified in
Lemma~\ref{L:prodinv}. Since $\tilde \phi$ is a Markov chain, it follows that $M=(M_t)_{t\geq 0}$
is a martingale. Since $|M_t|\leq 1$, it must converge a.s.\ to a limit in $[0,1]$, and hence
$|\tilde\phi_t|$ converges a.s.\ to a limit in $\N\cup\{0, \infty\}$. By Lemma~\ref{L:fluct},
a finite nonzero limit can be ruled out. Therefore almost surely,
\be
(1-\rho)^{|\tilde\phi_t|}\asto{t}1_{\{\tilde\phi_u=\emptyset\mbox{ for some }u\geq 0\}}.
\ee
Using the martingale property, it follows that for any finite $B\subset\Z$,
\be
M_0= (1-\rho)^{|B|}= \E[M_\infty] = \P\big(\tilde\phi_u=\emptyset\mbox{ for some }u\geq 0 \big| \tilde\phi_0=B\big).
\ee
Recalling the almost sure construction of the branching-coalescing point set $\tilde \xi$ and the dual process $\tilde \phi$ from
the Bernoulli net $\wN$ in Section \ref{S:braco}, we note that almost surely,
\be
\big\{\pi(t)\in B\mbox{ for some }\pi\in \wN(\ast,-\infty)\big\}
=\big\{\Phi_{t, s}(B)\neq\emptyset\ \forall s\leq t\big\}.
\ee
It follows that the backbone $\tilde \xi^*$ satisfies
\be
\P(\tilde\xi^*_t\cap B=\emptyset)= \P(\Phi_{t, s}(B)=\emptyset \mbox{ for some } s\leq t) = \P(|\tilde \phi_t|\to 0 \,|\, \tilde\phi_0=B)
= (1-\rho)^{|B|}.
\ee
Since this holds for each finite $B\sub\Z$, the claim follows.
\epro

\subsection{Density estimates}\label{S:denest}
For the discrete net $N_\eps$ with jump distribution $a(\cdot)$ and branching probability $\eps$, let $\wN_\eps:=\wN_{\psi_\eps}$ be the Bernoulli
net with kernel $\psi_\eps$ as in Lemma \ref{L:Bernbin}, coupled to $N_\eps$. The goal of this subsection is to bound the density of the branching-coalescing point set $(\tilde \xi^\eps_t)_{t\in \N_0}$ defined from $\wN_\eps$ as in \eqref{xiphi}, starting from the maximal configuration
$\tilde \xi^\eps_0=\Z$. We will show that the density of $\tilde \xi^\eps_t$ decays at the rate of $1/\sqrt{t}$ up to time of order $\eps^{-2}$,
see Theorem \ref{T:density}, which is comparable to that of coalescing random walks. By the coupling between $N_\eps$ and $\wN_\eps$, this gives an upper bound on the density of the branching-coalescing point set $\xi^\eps_t$ generated by the discrete net $N_\eps$ with $\xi_0=\Z$. In particular, this implies that diffusive scaling limits of $\xi^\eps$ are locally finite subsets of $\R$, which will play an important role in the proof of the upper bound in Section \ref{S:upproof}.

To simplify notation, from now on, we will omit $\eps$ from $\tilde \xi^\eps$ and $\tilde \phi^\eps$.
{ Recall the duality between $\tilde \xi$ and $\tilde \phi$ from (\ref{disdu}).}

\begin{theorem}[Density decay]\label{T:density}
	Let $\tilde \xi$ be the branching-coalescing point set defined from the Bernoulli net $\wN_\eps$ with $\tilde \xi_0=\Z$, and let $\tilde \phi$ be the dual of $\tilde \xi$ with $\tilde \phi_0=\{0\}$. Then there exists $C>0$ such that for all $\eps>0$ sufficiently small and all $T\leq \eps^{-2}$, we have
	\be\label{survival}
	\P\big(0\in\tilde\xi_T\big) = \P\big(|\tilde\phi_T|\neq0\big) \leq \frac{C}{\sqrt T}.
	\ee
	In particular, the same bound holds for the branching-coalescing point set $\xi$ defined from the discrete net $N_\eps$ with $\xi_0=\Z$.
\end{theorem}

We now sketch the key ideas behind the proof of Theorem \ref{T:density}. The first observation is that we can relate the process $|\tilde\phi_t|$ to a martingale. Instead of working with the martingale $M_t$ defined in (\ref{1rhophi}), it is more convenient to work with
\begin{equation}\label{mart}
	\wt M_t:= \frac{1-M_t}{\rho} = \frac{1-(1-\rho)^{|\tilde\phi_t|}}{\rho}.
\end{equation}
Note that $|\tilde \phi_t|=0$ if and only if $\wt M_t=0$, and the time $\tilde \phi$ hits the empty set is exactly the time $\wt M$ hits $0$.  For $\eps$ small, by \eqref{rhoreps}, we have $\rho=1 - e^{-r}=(2+o(1))\eps$. From \eqref{mart}, we see that $\wt M_t$ is close to $|\tilde\phi_t|$ as long as $|\tilde\phi_t| \ll 1/\eps$, and the two remain comparable as long as $|\tilde\phi_t|\leq c/\eps$ for some fixed constant $c>0$. We will show that for $T=O(\eps^{-2})$, the probability of the event $\{\max_{1\leq t\leq T}|\tilde\phi_t| \geq c/\eps\}$ is of order $O(1/\sqrt{T})$ and hence constitutes
part of the bound in \eqref{survival}. On the complementary event, $\wt M_t$ and $|\tilde\phi_t|$ are comparable up to time $T$. Lemma \ref{L:fluct} gives a lower bound on the fluctuations of $|\tilde\phi_t|$, which translates into a lower bound on the fluctuations of $\wt M_t$. This in turn implies
that the probability that $\wt M$ does not hit $0$ before time $T$ is comparable to the probability of the same event for a Brownian motion, which is also of order $O(1/\sqrt{T})$.
To make this last comparison rigorous, we need to Skorohod embed $\wt M_t$ in a Brownian motion (see e.g.\ \cite[Thm~12.16]{K97}), which requires an
additional estimate due to the randomness in the embedding.

\bt[Skorohod embedding]\label{T:Skoremb}
Let $(\wt M_t)_{t\in\N_0}$ be a martingale with $\wt M_0=1$ and induced filtration $(\Gi_t)_{t\in\N_0}$. Then there exist a Brownian motion $B$ with $B_0=1$ and stopping times $0=\tau_0\leq\tau_1\leq\cdots$, such that $(\wt M_t)_{t\in \N_0} \overset{d}{=} (B_{\tau_t})_{t\in \N_0}$ and
\bc
\dis\E[(\wt M_t - \wt M_{t-1})^2|\Gi_{t-1}] &=& \dis\E[\tau_t-\tau_{t-1}|\Fi_{t-1}],
\ec
where $(\Fi_t)_{t\in \N_0}$ is the filtration induced by the process $(\tau_t,B_{\tau_t})_{t\in \N_0}$.
Moreover, conditioned on $\Fi_{t-1}$, there exists a random pair $\alpha_t\leq 0 \leq \beta_t$ independent of $(B_s-B_{\tau_{t-1}})_{s\geq \tau_{t-1}}$
such that
\be\label{optcon}
\dis\tau_t =\inf\big\{s\geq\tau_{t-1}:B_s-B_{\tau_{t-1}}\in\{\alpha_t,\beta_t\} \big\},
\ee
\et

\bpro[of Theorem \ref{T:density}.]
The equality in \eqref{survival} follows from the duality relation \eqref{disdu}. Since $\tilde \phi_t=0$ if and only if $\wt M_t=0$, to prove the bound in \eqref{survival}, it suffices to bound $\P(\wt M_T\neq 0)$, for which we need Theorem \ref{T:Skoremb}.

The idea behind the construction of the stopping times $\tau_t$ in
(\ref{optcon}) is that conditioned on $\wt M_{t-1}=B_{\tau_{t-1}}$, the distribution
of the martingale increment $\wt M_t - \wt M_{t-1}$ can be decomposed into a mixture of two-point
distributions, where the mixture is given by the law of $(\alpha_t, \beta_t)$.
An $\{\alpha_t, \beta_t\}$-valued random variable with mean zero can then be embedded into the Brownian motion
$(B_{\tau_{t-1}+s}-B_{\tau_{t-1}})_{s\geq 0}$ by letting $\tau_t$ be the first time that $B_t$ hits $\{\alpha_t, \beta_t\}$.
Furthermore, since the martingale $(\wt M_t)_{t\in \N_0}$ is non-negative and will stay at $0$ once
it hits $0$, we can conclude that $\wt M$ first hits $0$ at time $t$ if and only if the Brownian motion $B$ (with $B_0=1$) hits $0$ at time
$\tau_t$. Therefore
\be\label{survival2}
\P(\wt M_T\neq0)=\P(B_s>0 \ \forall~0<s<\tau_T).
\ee

To bound the r.h.s.\ of \eqref{survival2}, we choose a constant $C_1>0$ and then define a stopping time $\sigma$ (w.r.t.\ the filtration of $\tilde \phi$ and $\wt M$) by
\be\label{stsig}
\sigma:=\inf\{t\in\Z:|\tilde\phi_t|\geq C_1/\rho\},
\ee
where
\be\label{rhoeps}
\rho=1-\prod_{x\in \Z} (1-\psi_\eps(x)) = 1- e^{-r} = (2+o(1))\eps
\ee
by \eqref{ala} and Remark \ref{R:reps}. Consider three events $E_1, E_2, E_3$, which are respectively the intersection of
the event in the r.h.s.\ of (\ref{survival2}) with each of the following three events
\be\label{intersection}
\{\tau_T\geq c_1 T \},~ \{\sig\leq T \}, \mbox{ and } \{\tau_T<c_1T, \sig>T \},
\ee
where $c_1>0$ is a small constant to be determined in Lemma~\ref{L:denupp1} below. To prove \eqref{survival},
it then suffices to show that $\P(E_i)\leq C/\sqrt{T}$ for each $i$.

For the first event, we have
\be
\P(E_1) \leq \P( B_s>0\ \forall\ 0<s<c_1 T) \leq \frac{1}{\sqrt{c_1 T}},
\ee
while the desired bounds on $\P(E_2)$ and $\P(E_3)$ will follow from Lemmas~\ref{L:denupp2} and \ref{L:denupp1} below.
\epro

\bl\label{L:denupp2}
Let $\wt N_\eps$ and $(\tilde \phi_t)_{t\in \N_0}$ be as in Theorem \ref{T:density}, and let $\sig$ be as in \eqref{stsig}.
Then there exists $C>0$ such that for all $\eps>0$ sufficiently small and all $T\leq \eps^{-2}$, we have
\be\label{optupp}
\P(\sig\leq T)\leq C\rho \leq \frac{3C}{\sqrt T}.
\ee
\el
\bpro By the definition of $\sigma$ in \eqref{stsig} and $\wt M$ in \eqref{mart}, we have
$$
\wt M_{\sig}\geq \rho^{-1}\big[1-(1-\rho)^{C_1/\rho}\big].
$$
Since $\rho=2\eps+o(\eps)$ by \eqref{rhoeps}, we can find $C>0$ such that $\wt M_{\sig}\geq 1/C\rho$ for all $\eps>0$
sufficiently small. Since $\wt M$ is a non-negative martingale, we then have
\be
1=\wt M_0=\E[\wt M_{\sig\wedge T}]\geq \frac{1}{C\rho} \P(\sig\leq T),
\ee
which implies \eqref{optupp}.
\epro

\bl \label{L:denupp1}
Let $\wt N_\eps$ and $(\tilde \phi_t)_{t\in \N_0}$ be as in Theorem \ref{T:density}, and let $\tau_t$ and $\sig$ be defined as in \eqref{optcon} and \eqref{stsig}, respectively. Then there exist $c_1, C>0$ such that for all $\eps>0$ sufficiently small and all $T\leq \eps^{-2}$, we have
\be
\P(\tau_T<c_1T, \sig>T)< \frac{C}{\sqrt T}.
\ee
\el
\bpro
Note that the event $\{\sig>T\}$ implies that $|\tilde\phi_t|\leq C_1/\rho$ for all $t\leq T$. Therefore by the definition
of $\wt M$ in \eqref{mart}, the increments of the martingale $\wt M$ satisfy
$$
\big|\wt M_t - \wt M_{t-1}\big|\geq C_2\big||\tilde\phi_t| - |\tilde \phi_{t-1}|\big| \qquad (1\leq t\leq T)
$$
for some constant $C_2>0$.

By Lemma~\ref{L:fluct}, with probability at least $\delta>0$ for some $\delta$ that is uniform in $\eps>0$ sufficiently
small and uniform in $\tilde \phi_{t-1}$, we have $||\tilde\phi_t|-|\tilde\phi_t||\geq 1$, and hence $\big|\wt M_t - \wt M_{t-1}\big|\geq C_2$.
Recall the construction of the stopping times $\tau_t$ in the Skorohod embedding in \eqref{optcon}, where conditioned on the $\sigma$-algebra $\Gi_{t-1}$
generated by $(\wt M_s)_{0\leq s\leq t-1}$, the law of $\wt M_t- \wt M_{t-1}$ is given by a mixture of the exit distribution of a Brownian motion $B'$ with $B'_0=0$ from $(\alpha_t, \beta_t)$, for some random vector $(\alpha_t, \beta_t)$. We note that the law of $(\alpha_t, \beta_t)$ must satisfy
\be\label{abt}
\P\big(\min\{|\alpha_t|, \beta_t\} \geq \alpha C_2/2\big) \geq \alpha/2,
\ee
since on the event $\{ \min\{|\alpha_t|, \beta_t\} \leq \alpha C_2/2 \}$, we have $\P(|B'_\tau|\geq C_2 |(\alpha_t, \beta_t)) \leq \alpha/2$,
where $\tau$ is the exit time of $B'$, with $B'_0=0$, from $[\alpha_t, \beta_t]$. Therefore \eqref{abt} must hold in order to have $\P(|\wt M_t-\wt M_{t-1}|\geq C_2|\Gi_{t-1})\geq \delta$.

It is then easily seen that $(\tau_t-\tau_{t-1})_{t\in \N}$ is stochastically dominated from below by a sequence of i.i.d.\ non-negative random
variables $(\gamma_{t})_{t\in\N}$ with positive mean, where the law of $\gamma_t$ depends only on $\alpha$ and $C_2$. In particular, we can
couple $(\tau_t)_{t\in \N_0}$ and $(\gamma_t)_{t\in \N}$ such that almost surely, $\tau_t\geq \sum_{s=1}^t \gamma_s$ for all $t\in \N$.
Choose any $0< c_1 <\E[\gamma_1]$. We then have
$$
\P(\tau_T \leq c_1 T) \leq \P\Big( \sum_{t=1}^T \gamma_t \leq c_1 T\Big) \leq C' e^{-T/C'} \leq \frac{C}{\sqrt T}
$$
uniformly in $T\in \N$ by elementary large deviation bounds. Furthermore, this bound is also uniform in $\eps>0$ small, thanks to Lemma~\ref{L:fluct}.
This concludes the proof.
\epro

\section{The upper bound}\label{S:upproof}

Let $N_\eps$ be the discrete net as in Theorem \ref{T:net}, let $\wt N_\eps$ be the Bernoulli net coupled to $N_\eps$ as in Lemma \ref{L:Bernbin},
and let $\Ni$ be the standard Brownian net. In Section \ref{S:lowproof}, it was shown that every subsequential weak limit of the rescaled discrete net
$S_{\eps} N_\eps$ contains a copy of the Brownian net $\Ni$. In this section, we show the matching upper bound that every subsequential weak limit of $S_{\eps}N_\eps$ contains no additional paths besides $\Ni$. { Note that we still need to prove tightness separately in Section \ref{S:tight}, since without tightness, subsequential weak limits may not exist along a particular subsequence, in which case our upper and lower bounds become irrelevant.}

We will start with a sufficient criterion for this last statement in Section \ref{S:crit}, which is based on a family of counting random variables. In Section \ref{S:convbb}, we will show that the subset of $N_\eps$ started at time $0$ from a random subset of $\Z$, distributed according to the product invariant law of the branching-coalescing point set generated by the Bernoulli net $\wt N_\eps$ (see Lemma \ref{L:prodinv}), converges under the scaling map $S_{\eps}$ to the subset of $\Ni$ started at time $0$ according to the stationary law of the branching-coalescing point set
generated by $\Ni$.  This is then used in Section \ref{S:fdd} to deduce the convergence of paths in $S_\eps N_\eps$ starting from finitely many space-time points. Finally, the desired upper bound on subsequential weak limits of
$S_{\eps}N_\eps$ is established in Section \ref{S:upppf} by verifying the criterion in Section \ref{S:crit}.

\subsection{An upper bound criterion}\label{S:crit}

Recall from Section \ref{S:topo} the space of paths $\Pi$, and the space $\Hi$ of compact subsets of $\Pi$. Let $\Pi_t \subset \Pi$ denote the set of
paths starting at time $t$, and let $\Ni_t:=\Ni \cap \Pi_t$. For a path $\pi\in\Pi$ and any $t\geq\sig_\pi$, let $\pi^t:=\{(\pi(s),s):s\geq t\}$ denote the truncation of $\pi$ at time $t$ so that $\sig_{\pi^t}=t$.
Recall that the Brownian net $\Ni$ satisfies the image set property stated in Proposition~\ref{P:image}, which implies that for any $\Hi$-valued random variable $\Xc$, if
$$
\Xc^t:=\{\pi^t:\pi\in\Xc, \sigma_\pi \leq t\}
$$
denotes the set of paths in $\Xc$ starting before or at time $t$ and truncated at time $t$, and almost surely $\cup_{p\in \Xc^t} p \subset \cup_{p\in \Ni_t} p$ for all $t\in \R$, then $\Xc\setminus\Ni= \emptyset$.
Motivated by this observation, we will formulate an upper bound criterion for the Brownian net using a family of {\em counting random variables} $\eta$ to be defined below. Analogous counting random variables were used in the upper bound criteria for the Brownian web $\Wi$ in \cite{NRS05}, which were in turn adapted from \cite{FINR04}.

Similar to the definition of the branching-coalescing point set from $\Ni$, for $\Xc\in\Hi$, we define
\be
\xi_\Xc(t,h;a,b):=\{\pi(t+h)\cap(a,b): \pi\in \Xc, \sigma_\pi\leq t  \}, \quad(t\in\R,h>0,a<b),
\ee
which is the point set generated at time $t+h$ by paths in $\Xc$ that start before or at time $t$.
We then define the family of counting random variables
\be
\eta_\Xc(t,h;a,b):=|\xi_\Xc(t,h;a,b)|, \quad(t\in\R,h>0,a<b).
\ee

Here is our upper bound criterion, which is in fact also a convergence criterion since it already assumes tightness and the desired lower bound.
\bt[Upper bound criterion] \label{T:upper}
Let $(\Xc_n)_{n\in\N}$ be a sequence of $\Hi$-valued random variables whose laws are tight.
If any of its subsequential weak limits $\Xc$ contains a copy of the Brownian net $\Ni$, and moreover for all $t\in\R,h>0,a<b$,
\be\label{countrv}
\E[\eta_\Xc(t,h;a,b)] \leq \E[\eta_\Ni(t,h;a,b)],
\ee
then $\Xc_n$ converges weakly to the Brownian net $\Ni$ as $n\to\infty$.
\et
\br
Due to the assumption that any subsequential weak limit $\Xc$ contains a copy of the Brownian net $\Ni$, it is easily seen that the inequality in \eqref{countrv} actually implies equality.
\er
\bpro
By assumption, $\Xc$ is coupled to a Brownian net $\Ni$ such that $\Ni \subset \Xc$ a.s., which implies that
$\cup_{p\in \Xc^t} p \supset \cup_{p\in \Ni_t} p$ for all $t\in\R$. It remains to show that $\Xc\setminus\Ni=\emptyset$ a.s.
By the image set property Proposition~\ref{P:image}, it suffices to show that a.s.,
\be\label{img}
\bigcup_{p\in \Xc^t} p \subset \bigcup_{p\in \Ni_t} p \quad \quad  \mbox{for all } t\in \R.
\ee
Note that \eqref{countrv} implies that a.s., for all
$t\in \Q$, $h\in \Q\cap (0,\infty)$ and $a<b$, we have $\eta_\Xc(t,h;a,b) \leq \eta_\Ni(t,h;a,b)$, which can be extended to all $t\in\R$ and
$h>0$ since $\Xc$ and $\Ni$ are compact sets of continuous paths. This immediately implies \eqref{img} and concludes the proof.
\epro

\subsection{Convergence of a stationary set of paths}\label{S:convbb}

For the standard Brownian net $\Ni$, recall from \cite[Prop.~1.15]{SS08} that its associated branching-coalescing point set has a unique invariant
law given by that of a Poisson point process on $\R$ with intensity $2$. Let $\Lambda$ be such a Poisson point set independent of $\Ni$, and
let $\Bi:= \Ni(\Lambda \times \{0\})$ denote the subset of $\Ni$
starting from $\Lambda$ at time $0$. Similarly, let $\Lambda_\eps$ be a random subset of $\Z$, whose law is invariant for
the branching-coalescing point set generated by the Bernoulli net $\wt N_\eps$, see Lemma \ref{L:prodinv} and \eqref{rhoeps}. Let
$B_\eps:=N_\eps(\Lambda_\eps \times \{0\})$, the set of paths in the discrete net $N_\eps$ (not $\wt N_\eps$) started from $\Lambda_\eps$ at time $0$. Here is the main result of this subsection.

\bt[Convergence of a stationary set of paths]\label{T:bbconv}
Assume that the law of the family of rescaled discrete nets $S_\eps N_\eps$ is tight. As defined above, let $\eps\Lambda_\eps$ and $\Lambda$ be random subsets of $\eps\Z$ and $\R$, respectively, which are now regarded as random point measures on $\R$, equipped with vague convergence. Then we have $(\eps \Lambda_\eps, S_{\eps}B_\eps) \Rightarrow (\Lambda, \Bi)$ as $\eps\downarrow 0$.
%Furthermore, we can find a coupling such that  $S_\eps \Lambda_\eps \Rightarrow \Lambda$ a.s.\ w.r.t.\ the vague topology for locally finite measures on $\R$, and conditioned on $\Lambda_\eps$ and $\Lambda$, a.s.\ $S_{\eps}B_\eps \Rightarrow \Bi$.
\et
\br
We need to enhance the convergence $S_{\eps}B_\eps \Rightarrow \Bi$ with the convergence of the random measures $\eps \Lambda_\eps \Rightarrow \Lambda$, so that we can find a coupling between $(\eps \Lambda_\eps, S_{\eps}B_\eps)$ and $(\Lambda, \Bi)$ such that not only do the locations of the atoms in $\eps\Lambda_\eps$ converge to that in $\Lambda$, but also the number of atoms converge. This will be used to deduce the convergence of
the rescaled discrete net paths starting from a single point, see Proposition~\ref{P:1dconv} below.
\er
\bpro
Thanks to Lemma \ref{L:prodinv} and \eqref{rhoeps}, the convergence $\eps \Lambda_\eps \Rightarrow \Lambda$ as random measures on $\R$ is clear. By Skorohod's representation theorem, we can choose a coupling such that a.s., $\eps \Lambda\to \Lambda$ w.r.t.\ the vague topology. It remains to show that conditioned on a typical realisation of $\Lambda_\eps$ and $\Lambda$, $S_\eps B_\eps \Rightarrow \Bi$.

The tightness of $(S_{\eps}N_\eps)_{\eps>0}$ implies the tightness of $(S_{\eps}B_\eps)_{\eps>0}$.
It is easily seen that the analogue of Theorem \ref{T:upper} holds if $\Ni$ is replaced by $\Bi=\Ni(\Lambda \times\{0\})$
for any closed $\Lambda \subset \R$ independent of $\Ni$. To show that any subsequential weak limit $\Ai$ of $S_\eps B_\eps$
contains a copy of $\Bi$, we can apply the same discrete approximation argument as in the proof of Theorem \ref{T:lower}, provided we use the fact from
Theorem \ref{T:Nhop} that $\Ni(\Lambda \times \{0\})$ can be constructed via the hopping construction where the first hopping path is chosen from $\Wi^1(\Lambda\times \{0\}) \cup \Wi^2(\Lambda\times \{0\})$. In particular, this implies that a.s.,
\be
\E[\eta_\Ai(t,h;a,b)|\Lambda]\geq \E[\eta_\Bi(t,h;a,b) |\Lambda] \qquad \mbox{for all } t\geq 0, h>0, a<b\in \R,
\ee
To show the reverse inequality, it suffices to show the inequality after taking expectation, i.e., show that
\be\label{countrv2}
\E[\eta_\Ai(t,h;a,b)]\leq \E[\eta_\Bi(t,h;a,b)] \qquad \mbox{for all } t\geq 0, h>0, a<b\in \R,
\ee
where we note that $\E[\eta_\Bi(t,h;a,b)]= 2(b-a)$ by \cite[Prop.~1.15]{SS08}. Recall from Lemma \ref{L:Bernbin} that $N_\eps$ can be coupled to a Bernoulli net $\wt N_\eps$
such that $N_\eps \subset \wt N_\eps$ a.s., and hence also $B_\eps \subset \wt B_\eps:= \wt N_\eps(\Lambda_\eps \times \{0\})$.
By Lemma \ref{L:prodinv}, $\Lambda_\eps$ is an invariant law for the branching-coalescing point set generated by the Bernoulli net $\wt N_\eps$, and $\Lambda_\eps$ has density $\rho=(2+o(1))\eps$ { as $\eps\downarrow 0$ by \eqref{rhoeps}.} It follows by Fatou's lemma that under the scaling map
$S_\eps : (x, t) \to (\eps x, \varsigma^2\eps^2 t)$ with $\eps\downarrow 0$, at each time $t\geq 0$, the point set generated at time $t$ by paths in the subsequential weak limit $\Ai$ has density at most $2$. Therefore $\E[\eta_\Ai(t,h;a,b)]\leq 2(b-a)$, which concludes the proof of \eqref{countrv2}.
\epro

\subsection{Convergence of finite dimensional distributions}\label{S:fdd}

We now use Theorem \ref{T:bbconv} to deduce the convergence of paths in the discrete net $N_\eps$ starting from finitely many space-time points. We first treat the one point case. Recall that for $\Xc\in\Hi$, $A\subseteq\R^2$ and $z\in\R^2$, $\Xc(A)$ denotes the set of paths in $\Xc$ starting from the space-time set $A$, and $\Xc(z) =\Xc(\{z\})$.

\bp[Convergence of one dimensional distribution]\label{P:1dconv}
Assume that the law of the family of rescaled discrete nets $S_\eps N_\eps$ is tight. Let $N_\eps(0,0)$ be the subset of paths in $N_\eps$ starting from the origin at time $0$.
Then we have $S_{\eps}N_\eps(0,0) \Rightarrow \Ni(0,0)$ as $\eps \downarrow 0$.
\ep
\bpro
Recall from Section \ref{S:topo} the space $\Hi$ of compact subsets of the paths space $\Pi$. Let $\Pi(0,0)$ denote the subset of $\Pi$ starting
from $(0,0)$, and let $\Hi(0,0)$ denote the space of compact subsets of $\Pi(0,0)$. We will prove $S_{\eps}N_\eps(0,0) \Rightarrow \Ni(0,0)$ by regarding $S_{\eps}N_\eps(0,0)$ and $\Ni(0,0)$ as $\Hi(0,0)$-valued random variables, which also implies weak convergence if they are regarded as $\Hi$-valued random variables. We need to show that for any bounded continuous $f:\Hi(0,0) \to \R$, we have
\be\label{fNeps}
\E\big[f\big( S_\eps N_\eps(0,0)\big)\big]
\underset{\eps\down0}{\longrightarrow}\E\big[f\big(\Ni(0,0)\big)\big].
\ee

Let $\Pi(\R\times \{0\})$ denote the subset of $\Pi$ starting from $\R$ at time $0$, and let $\Hi(\R\times \{0\})$ denote the space of compact subsets of $\Pi(\R\times \{0\})$. Let $\Hi_{\rm fin}(\R\times\{0\})$ denote the subset of $\Hi(\R\times\{0\})$ such that for $\Xc\in\Hi_{\rm fin}(\R\times\{0\})$, $\Xc_0=\{\pi(0): \pi \in \Xc\}$ is locally finite subset of $\R$. For any bounded continuous $f: \Hi(0,0) \to \R$, we can define $g: \Hi_{\rm fin}(\R\times\{0\}) \to \R$ by
\be\label{test1}
g(\Xc)=\sum_{x\in \Xc_0\cap [0,1]}f(\Xc(x, 0) - x) \qquad  (\, \Xc \in \Hi_{\rm fin}(\R\times\{0\})\, ),
\ee
where $\Xc(x,0)-x$ denotes a spatial shift of $-x$ of paths in $\Xc(x,0)$.
By Theorem \ref{T:bbconv}, we can find a coupling such that a.s.\ $S_\eps B_\eps \to \Bi$ in $\Hi(\R\times \{0\})$ and $\eps \Lambda_\eps \to \Lambda$ w.r.t.\ the vague topology, i.e., the number of atoms in the random set $\eps \Lambda_\eps \cap [0,1]$ and their locations converge to that
of $\Lambda \cap [0,1]$, where we note that a.s.\ $\Lambda$ has no atom on the boundary of $[0,1]$. This is easily seen to imply
\be\label{test2}
\E\big[g\big(S_{\eps}B_\eps \big)\big] \underset{\eps\down0}{\longrightarrow} \E\big[g\big(\Bi\big) \big].
\ee
By translation invariance and the fact that $\Lambda_\eps$ has density $\rho=(2+o(1))\eps$, see \eqref{rhoeps}, and $\Lambda$ has density $2$, \eqref{test2} is equivalent to
$$
\eps^{-1}\rho \E\big[f\big( S_\eps N_\eps(0,0)\big)\big] \to 2 \E[f(\Ni(0,0)].
$$
The desired convergence in \eqref{fNeps} follows immediately.
\epro
\medskip

We now treat the case with finitely many starting points in space and time.
\bt\label{T:fdd}
Assume that the law of the family of rescaled discrete nets $S_\eps N_\eps$ is tight. Let $k\in\N$. Suppose that for each $1\leq i\leq k$,
$z^\eps_i\in \Z^2$ satisfies $S_\eps z^\eps_i \to z_i\in \R^2$ as $\eps\downarrow 0$. Then as random
variables taking values in $\Hi^k$, we have
\be
(S_{\eps} N_\eps (z_1^\eps) ,\cdots, S_{\eps} N_\eps(z_k^\eps)) \underset{\eps\down0}{\Longrightarrow}  (\Ni(z_1),\ldots, \Ni(z_k)).
\ee
\et
\bpro
Tightness follows by assumption. By the same discrete approximation argument as in the proof of Theorem \ref{T:bbconv}, we note that every subsequential weak limit $(\Ai_1,\cdots,\Ai_k)$ of $(S_{\eps} N_\eps (z_1^\eps) ,\cdots, S_{\eps} N_\eps(z_k^\eps))$ can be coupled to a Brownian net $\Ni$
such that $\Ai_i \supset \Ni(z_i)$ for each $1\leq i\leq k$. On the other hand, the convergence in Theorem \ref{P:1dconv} implies that a.s., we must have $\Ai_i=\Ni(z_i)$ for each $i$. This concludes the proof.
\epro

\subsection{Proof of the upper bound}\label{S:upppf}

In this subsection, we prove the desired upper bound on subsequential weak limits of
$S_{\eps}N_\eps$ by verifying condition \eqref{countrv} in Theorem \ref{T:upper}.

\bt[Upper bound]\label{T:Nupper}
Let $N_\eps$ be the discrete net as in Theorem \ref{T:net}. For any subsequential weak limit $\Ni^\ast$ of $S_{\eps}N_\eps$ as $\eps\down0$, there exists a standard Brownian net $\Ni$ such that almost surely, $\Ni=\Ni^\ast$.
\et

\bpro
By going to a subsequence if necessary, we can assume that the rescaled discrete net $S_{\eps}N_\eps$ converges in distribution to a limit $\Ni^\ast$, which by Theorem~\ref{T:lower}, is coupled to a standard Brownian net such that $\Ni\subseteq\Ni^\ast$ a.s.
In order to show $\Ni^\ast\subseteq\Ni$, by Theorem~\ref{T:upper}, it only remains to check (\ref{countrv}) for the counting random variables $\eta_{\Ni^\ast}(t,h;a,b)$ and $\eta_{\Ni}(t,h;a,b)$.
By translation invariance in time, we only need to show that for all $h>0$ and $a<b$,
\be\label{count}
\E[\eta_{\Ni^\ast}(0,h;a,b)]\leq \E[\eta_\Ni(0,h;a,b)].
\ee
We prove (\ref{count}) by an argument similar to the one in \cite[Section~6]{NRS05} in the proof of convergence to the Brownian web.
The strategy is to first show that
\begin{equation}\label{Nstardelta}
	\xi^{\Ni^\ast}_\delta:=\{\pi(\delta): \pi \in \Ni^\ast, \sigma_\pi\leq 0\},
\end{equation}
the branching-coalescing point set generated at time $\delta$ by paths in $\Ni^\ast$ started before or at time $0$, is a locally finite subset of $\R$ (sometimes called {coming down from infinite}). We then apply Theorem \ref{T:fdd} to show that the set of truncated paths $\{\pi^\delta: \pi \in \Ni^\ast, \sigma_\pi\leq 0\}$, truncated at time $\delta$, is distributed as the set of Brownian net paths $\Ni(\xi^{\Ni^\ast}_\delta \times \{\delta\})$ starting from $\xi^{\Ni^\ast}_\delta$ at time $\delta$, which allows us to bound
$$
\E[\eta_{\Ni^\ast}(0,h;a,b)]\leq \E\big[\eta_{\Ni}(\delta, h-\delta; a,b)\big].
$$
Sending $\delta\downarrow 0$ then gives \eqref{count}.

We now give the details of the proof. By the density bound in Theorem \ref{T:density}, there exists a constant $C$ such that for all $0<\delta\leq h \wedge 1$ and $a<b$,
\be\label{count0.5}
\limsup_{\eps\downarrow 0}\E[\eta_{S_{\eps}N_\eps}(0, \delta;a,b)]\leq \frac{C\varsigma (b-a)}{\sqrt{\delta}} <\infty.
\ee
Since we have assumed $S_{\eps}N_\eps\Rightarrow \Ni^\ast$, by Fatou's Lemma, we must have $\E[\eta_{\Ni^\ast}(0,\delta;a,b)]<\infty$,
which implies that $\xi^{\Ni^\ast}_\delta$, defined in \eqref{Nstardelta} above, is a.s.~locally finite. Let $\xi^\eps_t$ denote the
branching-coalescing point set generated by $N_\eps$, starting from $\xi^\eps_0=\Z$. Then thanks to the uniform density bound in \eqref{count0.5}
and the convergence $S_{\eps}N_\eps\Rightarrow \Ni^\ast$ as $\Hi$-valued random variables, by going to a further subsequence if necessary, we have $\eps \xi^\eps_{\delta/\varsigma^2\eps^2} \Rightarrow \nu^{\Ni^\ast}_\delta$ as random variables taking values in the space of locally finite measures
on $\R$, equipped with vague convergence. We can assume a coupling such that $\eps \xi^\eps_{\delta/\varsigma^2\eps^2} \to \nu^{\Ni^\ast}_\delta$ a.s., and the support of $\nu^{\Ni^\ast}_\delta$ equals $\xi^{\Ni^\ast}_\delta$, where $\nu^{\Ni^\ast}_\delta$ may assign integer mass greater than $1$ to some points in $\xi^{\Ni^\ast}_\delta$. For any finite interval $(c,d)$,  the a.s.\ vague convergence $\eps \xi^\eps_{\delta/\varsigma^2\eps^2} \to \nu^{\Ni^\ast}_\delta$ and the fact that $\nu^{\Ni^\ast}_\delta$ assigns no mass to the boundary of $(c, d)$ imply
that the number of points in $\eps \xi^\eps_{\delta/\varsigma^2\eps^2} \cap (c, d)$ will converge, and their locations will also converge. We can now apply Theorem \ref{T:fdd} to $N_\eps(\xi^\eps_{\delta/\varsigma^2\eps^2} \cap (c/\eps, d/\eps) \times \{\delta/\varsigma^2\eps^2\})$,\footnote{The argument here differs from the one in \cite[Section~6]{NRS05} for the Brownian web. Translated into our setting, it would mean using the weaker convergence of $\eps \xi^\eps_{\delta/\varsigma^2\eps^2} \to \xi^{\Ni^\ast}_\delta$ as closed subsets of $\R\cup\{\pm \infty\}$ w.r.t.\ the Hausdorff distance. The rest of the argument would then require a stronger version of Proposition \ref{P:1dconv} (cf.\ \cite[Lemma 6.5]{NRS05}), which would establish the convergence of $N_\eps(D_\eps)\Rightarrow \Ni(0,0)$ for any closed $D_\eps \subset \R^2$ that converges in Hausdorff distance to $\{(0,0)\}$.} which by the Markov property of paths in $N_\eps$, is distributed as branching-coalescing random walks starting from the set $\xi^\eps_{\delta/\varsigma^2\eps^2}\cap (c/\eps, d/\eps)$ at time $\delta/\varsigma^2\eps^2$. In particular, conditioned on a realisation of the coupling
$\eps \xi^\eps_{\delta/\varsigma^2\eps^2} \to \nu^{\Ni^\ast}_\delta$, we have
$$
S_\eps \Big( N_\eps\big(\xi^\eps_{\delta/\varsigma^2\eps^2} \cap (c/\eps, d/\eps) \times \{\delta/\varsigma^2\eps^2\}\big) \Big) \Rightarrow \Ni(\xi^{\Ni^\ast}_\delta \cap (c, d) \times\{\delta\}) \qquad (c<d\in \R),
$$
where $\Ni$ is a standard Brownian net. This implies that we can couple $\Ni^\ast$ with a standard Brownian net $\Ni$ such that a.s., the set of
truncated paths
$$
\{\pi^\delta:\pi \in \Ni^\ast, \sigma_\pi\leq0\} = \{ (\pi(t), t)_{t\geq \delta}: \pi \in \Ni^\ast, \sigma_\pi\leq 0\}
$$
is contained in $\Ni(\xi^{\Ni^\ast}_\delta \times\{\delta\})$. Furthermore, by comparing with paths in $\Ni(\R\times\{\delta\})$, we have
\be\label{count1}
\eta_{\Ni^\ast}(0,h;a,b)\leq \eta_{\Ni(\xi^{\Ni^\ast}_\delta \times\{\delta\})}(\delta,h-\delta;a,b) \leq \eta_{\Ni}(\delta,h-\delta;a,b).
\ee
Take expectation on both sides and apply the density calculation for $\Ni$ in \cite[Prop.~1.12]{SS08}, we obtain
\be\label{count2}
\E[\eta_{\Ni^\ast}(0,h;a,b)]\leq \E\big[\eta_{\Ni}(\delta, h-\delta; a,b)\big]=(b-a)\cdot\Big(\frac{e^{-(h-\delta)}}{\sqrt{\pi (h-\delta)}}+2\Phi(\sqrt{2(h-\delta)})\Big),
\ee
where $\Phi$ is the distribution function of the standard normal. Sending $\delta\downarrow 0$ then gives \eqref{count} and concludes the proof.
\epro

\section{Tightness}\label{S:tight}
In this section, we prove that the laws of the rescaled discrete nets $(S_\eps N_\eps)_{\eps\in (0,1)}$ are tight. In Section \ref{S:tcrit}, we recall the standard tightness criterion for a family of $\Hi$-valued random variables. To address the fundamental difficulty created by branching, we introduce in Section \ref{S:RBP} the so-called relevant branching points of the branching-coalescing random walks. In Sections \ref{S:1RBP} and \ref{S:kRBP}, we show that such relevant branching points are rare on small time scales, which then allows us to implement in Section \ref{S:multi}
the multiscale argument used in \cite{BMSV06} to verify the tightness criterion for the rescaled discrete webs.

\subsection{Tightness criterion} \label{S:tcrit}

We recall here the standard tightness criterion for a sequence of $\Hi$-valued random variables $(\Xc_n)_{n\in\N}$, which was first formulated in \cite[Proposition~B.1]{FINR04}.

\bp[Tightness criterion]\label{P:tig}
A sequence of $\Hi$-valued random variables $(\Xc_n)_{n\in\N}$ is tight if for any space-time box $\Lambda_{L,T}:=[-L,L]\times[-T,T]$ with $L,T>0$
and every $M>0$,
\begin{equation}\label{tight}
	\lim_{\delta\downarrow0}\frac{1}{\delta}\limsup_{n\to\infty}\sup_{(x,t)\in
		\Lambda_{L,T}} \P\big(\Xc_n\in A_{M,\delta}(x,t)\big)=0,
\end{equation}
where $A_{M,\delta}(x,t)\subset \Hi$ consists of compact sets of paths $K\in \Hi$, such that $K$ contains some path which intersects the box
$[x-M,x+M]\times[t,t+\delta]$, and at a later time, intersects the left or right boundary of the larger box $[x-2M,x+2M]\times[t,t+2\delta]$ $($see Figure \ref{fig:tight}$)$.
\ep

\begin{figure}[ht!]
	\begin{center}
	\begin{tikzpicture}[ 
		x=0.7cm,
		y=0.7cm,
		>=stealth,
		line cap=round,
		line join=round
		]
		
		% Outer rectangle: width 4M, height 2\delta
		\draw[very thick] (0,0) rectangle (12,6);
		
		% Inner rectangle: width 2M, height \delta
		\draw[very thick] (3,0) rectangle (9,3);
		
		% Labels for the rectangles
		\node[above] at (6,6) {$4M$};
		\node[left]  at (0,3) {$2\delta$};
		
		\node[above] at (6,3) {$2M$};
		\node[left]  at (3,1.5) {$\delta$};
		
		% Distinguished point (x,t)
		\fill (6,0) circle (3pt);
		\node[below=3pt] at (6,0) {$(x,t)$};
		
		% Curve \pi
		\draw[very thick]
		(7.5,0.9)
		.. controls (8.8,1.35) and (10.6,1.75) .. (10.9,2.45)
		.. controls (11.1,2.9) and (10.1,3.05) .. (10.15,3.65)
		.. controls (10.2,4.45) and (11.5,4.85) .. (13.8,5.85);
		
		\node[right] at (10.1,1.55) {$\pi$};
		
		% Space-time direction indicator
		\draw[->,thick] (-1.,-0.7) -- (-0.1,-0.7)
		node[right] {\small Space};
		\draw[->,thick] (-1.,-0.7) -- (-1.,0.2)
		node[above] {\small Time};
		
	\end{tikzpicture}	
		\caption{A path $\pi$ causing the unlikely event $\{\Xc_n\in A_{M,\delta}(x,t)\}$ to occur.}
		\label{fig:tight}
	\end{center}
\end{figure}
\begin{remark}
	Proving tightness of $(\Xc_n)_{n\in\N}$ amounts to controlling the uniform modulus of continuity of paths in $\Xc_n$ in such a way that the probability
	bounds are uniform in $n\in\N$. Thanks to the compactification of space-time in Section \ref{S:topo}, we can restrict to a large but finite space-time window. Dividing the time window into intervals of length $\delta$ and controlling the uniform modulus of continuity on each time interval then leads to the condition \eqref{tight}. For more details, see e.g.\ \cite[Section 6.6]{SSS17}.
\end{remark}

\subsection{Relevant branching points}\label{S:RBP}

To verify the tightness condition (\ref{tight}) for the rescaled discrete nets $S_\eps N_\eps$, we will adapt the multiscale argument
developed in \cite{BMSV06}  to prove tightness for the rescaled discrete webs. A key ingredient is to show that the probability that some of the random walks and their descendants (due to branching) will exit the left or right boundary of the larger box $[x-2M,x+2M]\times[t,t+2\delta]$ in Figure \ref{fig:tight} is very small. Branching creates fundamental difficulties for this estimate, which motivates us to consider the so-called {\em relevant branching points (RBP)} of the discrete net $N_\eps$. Controlling the number of RBP along paths in $N_\eps(0,0)$ then allows us to show that the probability that some path in $N_\eps(0,0)$ makes a large excursion as depicted in Figure \ref{fig:tight} is comparable to the probability that a single random walk starting from $(0,0)$ makes a large excursion. We now give the precise definition of RBP.

Let $A\subset \Z^2$ be a set of space-time points. Given $S, U\in \Z$ with $S<U$, a space-time point $z=(x,t)\in\Z^2$ with $t\in(S,U)$ is called a {\em $(S,U)$-relevant branching point $($RBP$)$} of $N_\eps(A)$ if there exist two paths $\pi_1,\pi_2\in N_\eps$ starting from $A$, such that
\be\label{RBPdef}
\pi_1(s)=\pi_2(s)~\forall s\in[S,t],\quad \pi_1(t)=\pi_2(t)=x,\quand \pi_1(s)\neq\pi_2(s)~\forall s\in(t,U]\cap \Z.
\ee
We call $t$ an $(S,U)$-relevant branching time of $N_\eps(A)$. When things are clear from the context, we will simply call $z$ (resp.\ $t$) a relevant branching point (resp.\ time).  We say that a path $\pi$ passes through $z=(x,t)$ if $\pi$ starts before or at time $t$ and $\pi(t)=x$.

We define the branching-coalescing point set generated by $N_\eps(A)$ by
\be
\xi_t^A := \{\pi(t): \pi \in N_\eps(A)\}.
\ee
To determine how paths in $N_\eps(A)$ lead from the set $\xi^A_S\subset \Z$ at time $S$ to the set $\xi^A_U$ at time $U$, we consider the set of $(S, U)$-RBP of $N_\eps(A)$ and how each path $\pi \in N_\eps(A)$ passes through these RBP's. More precisely, we build a directed graph $G_{S, U}(N_\eps(A))= (V, E)$ with the vertex set $V$ consisting of $\xi^A_S \times \{S\}$, the set of $(S, U)$-RBP of $N_\eps(A)$, and $\xi^A_U\times \{U\}$.
For any two vertices $z_1=(x_1, t_1), z_2=(x_2, t_2)\in V$ with $S\leq t_1<t_2\leq U$, we draw a directed edge $e$ from $z_1$ to $z_2$ if there exists
$\pi \in N_\eps(A)$ with starting time $\sigma_\pi\leq S$, and $\pi$ pass through $z_1$ and $z_2$, but no $(S, U)$-RBP in between. Given such a
directed edge $e=(z_1, z_2)$, we call  $e$ an outgoing edge at $z_1$ and an incoming edge at $z_2$.

The definition of RBP above for the discrete net $N_\eps$ is motivated by a similar definition of {\em relevant separation points} for the Brownian net $\Ni$ in \cite[Section 2.3]{SSS09} and \cite[Section 6.1]{SSS17}. The Brownian net analogue of the directed graph $G_{S, U}(N_\eps(A))$ was studied in \cite[Section 6.2]{SSS17} and was called the {\em finite graph representation} of the Brownian net $\Ni$ in the time interval $[S, U]$, which played a crucial role in \cite{SSS17} in showing that the coupling between the discrete sticky webs and the discrete net converge to a coupling between the
sticky Brownian webs and the Brownian net.

For the discrete net $N_\eps$, we show that the RBP and the directed graph $G_{S, U}$ share similar properties to their continuum counterparts for the Brownian net $\Ni$.

Our first result shows that if $\pi\in N_\eps(A)$ starting at or before time $S$ passes through $z_1=(x_1, t_1)$ and $z_2=(x_2, t_2)$ for some directed edge $(z_1, z_2)$ in $G_{S, U}(N_\eps(A))$, then descendants of $\pi$ resulting from branching in the time interval $(t_1, t_2)$ must all pass through $z_2$.

\bl[No effective branching] \label{L:nobran}
For $S<U$ and $A\subset \Z^2$, let $G_{S, U}(N_\eps(A))=:G=(V, E)$ be defined as above. Let $e=(z_1,z_2)$, with $z_i=(x_i, t_i)\in\Z^2$, be a directed edge of $G$. If $\pi\in N_\eps(A)$ with $\sigma_\pi\leq S$ passes through $z_1$ and $z_2$ and no $(S,U)$-RBP in between, then there is no effective branching along $\pi$ in the time interval $(t_1,t_2)$, that is, for any $\varpi\in N_\eps(A)$ that passes through $z_1$ and $\varpi(s)=\pi(s)$ for some $s\in(t_1,t_2)\cap \Z$, $\varpi$ must pass through $z_2$ and no $(S, U)$-RBP in the time interval $(s, U)$.
\el
\bpro
We will implicitly assume that all times are integers. By assumption, either $t_2=U$ or $z_2$ is an $(S, U)$-RBP. If $t_2=U$ and $\varpi$ does not pass through $z_2$, then a.s.\ $t_0:=\sup\{n\leq U: \varpi(n)=\pi(n)\}<U$. By hopping from $\pi$ to $\varpi$ at time $t_0$, which still defines a path in
$N_\eps(A)$, it is easily seen that $z_0=(\pi(t_0),t_0)$ is also an $(S, U)$-RBP for $N_\eps(A)$. This contradicts the assumption that $\pi$ does not
pass through any RBP between $z_1$ and $z_2$.

\begin{figure}[ht!]
	\begin{center}
		\includegraphics[width=0.8\textwidth]{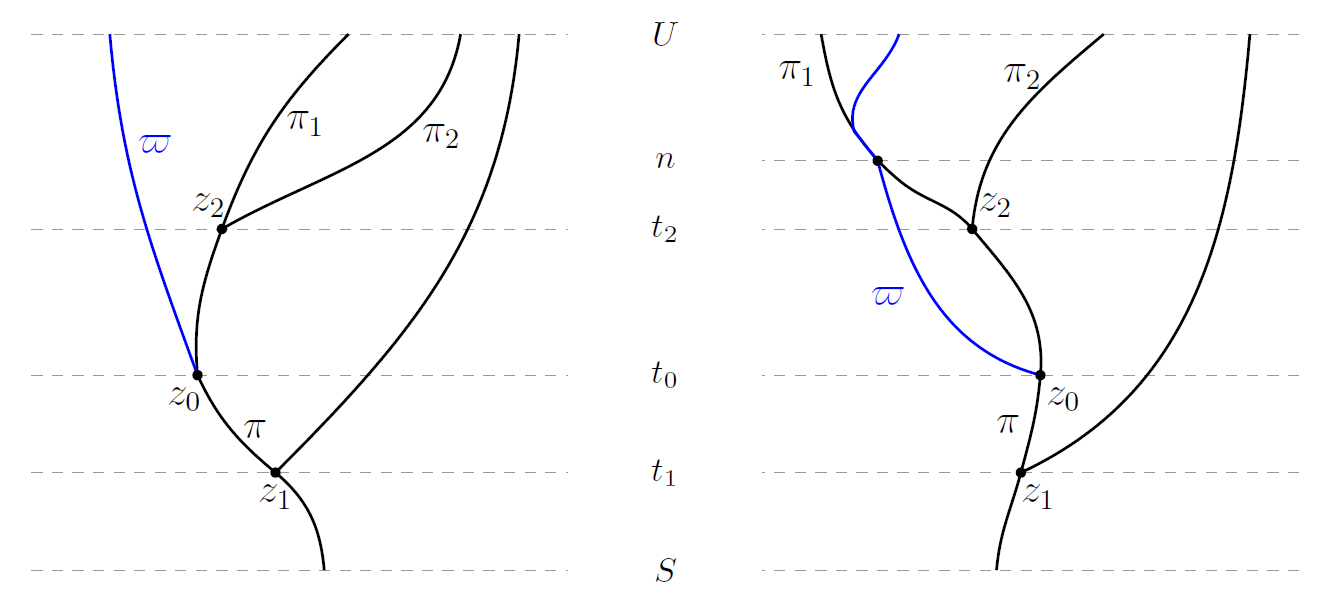}
		\caption{Two cases of $\varpi$ branching off from $\pi$ at $z_0$, but does not pass through the RBP $z_2$.}\label{fig:nobranching}
	\end{center}
\end{figure}
If $z_2$ is an $(S, U)$-RBP, then there exist two paths $\pi_1, \pi_2 \in N_\eps(z_2)$ which do not meet on the time interval $(t_2, U]$,
see Figure \ref{fig:nobranching}. If $\varpi$ does not pass through $z_2$, then we have $t_0:=\sup\{n\leq t_2: \varpi(n)=\pi(n)\}\in(s,t_2)$.
Denote $z_0=(\pi(t_0),t_0)$. We claim that $z_0$ must be a RBP, which again leads to a contradiction. Indeed, if $\varpi$ intersects
neither $\pi_1$ nor $\pi_2$ on the time interval $(t_2, U]\cap \Z$, then hopping from $\pi$ to $\pi_1$ at time $t_2$, and hopping from
$\pi$ to $\varpi$ at time $t_0$, identifies two paths in $N_\eps(A)$ which shows that $z_0$ is a RBP. If $\varpi$ intersects $\pi_1\cup\pi_2$
for the first time at some $n\in (t_2, U]$ (say the intersection is with $\pi_1$), then hopping from $\pi$ to $\varpi$ at time $t_0$ and then
from $\varpi$ to $\pi_1$ at time $n$, and hopping from $\pi$ to $\pi_2$ at time $t_2$, identifies two paths in $N_\eps(A)$ which shows
that $z_0$ is a RBP.

Lastly we note that since all $\varpi\in N_\eps$ that coincides with $\pi$ at some time $s\in (t_1, t_2)$ must pass through $z_2$, it means that
$z_2$ is a bottleneck and hence $\varpi$ cannot pass through any $(S, U)$-RBP before reaching $z_2$.
\epro

As a corollary of Lemma \ref{L:nobran}, we show that the reason a branching point along some $\pi \in N_\eps(A)$ with $\sigma_\pi\leq S$
fails to be $(S, U)$-RBP is because all its descendants must pass through a common bottleneck before or at time $U$.

\begin{corollary}[Non-relevant branching points]\label{C:nonrel}
	Let $S<U$. If some $\pi \in N_\eps$ with starting time $\sigma_\pi\leq S$ passes through $z=(x,t)\in \Z^2$ with $t\in (S, U)$, and $z$ is a
	branching point but not an $(S, U)$-RBP for $N_\eps$, then there exists $z'=(x', t')$ with $t'\in (t, U]$ such that all
	paths in $N_\eps(z)$ pass through $z'$.
\end{corollary}
\bpro
Let $z_1$ be the last vertex in $G_{S, U}(N_\eps)$ that $\pi$ passes through before time $t$, and let $z_2$ be the first vertex in
$G_{S, U}(N_\eps)$ that $\pi$ passes through after time $t$. Then $(z_1, z_2)$ is a directed edge in $G_{S, U}$, and hence Lemma \ref{L:nobran}
implies that all paths in $N_\eps(z)$ must pass through $z_2$.
\epro

Our next result shows that for each $(S, U)$-RBP $z$, the two branches of paths starting from $z$ will not intersect until one branch visits
the next RBP or reaches time $U$.

\bl[Disjoint paths from relevant branching points] \label{L:nocoa}
For $S<U$ and $A\subset \Z^2$, let $G_{S, U}(N_\eps(A))$ be defined as before. Let $z=(x,t)$ be a $(S,U)$-RBP.
If $\pi_1, \pi_2\in N_\eps(z)$ follow the two different branches from $z$, i.e., $\pi_1(t+1)\neq \pi_2(t+1)$, then $\pi_1$
and $\pi_2$ follow two different edges $(z, z_1)$ and $(z, z_2)$ in the directed graph $G_{S, U}$, where $z_1=(x_1, t_1)\neq z_2=(x_2, t_2)$, and
\be
\pi_1(s)\neq\pi_2(s)\quad  \mbox{for all } t< s\leq \min\{t_1, t_2\}.
\ee
\el
\bpro
Without loss of generality, we assume that $t_1\leq t_2$. Suppose that $\pi_1(t_0)=\pi_2(t_0)$ for some $t<t_0<t_1$. Then
Lemma~\ref{L:nobran} implies that $\pi_2$ must also pass through $z_1=(x_1, t_1)$ and follows the same directed $(z, z_1)$
in $G_{S, U}$. Since there is only binary branching in the discrete net $N_\eps$, any path in $N_\eps(z)$ must either coincide
with $\pi_1$ or $\pi_2$ at time $t+1$. Again, Lemma~\ref{L:nobran} implies that all such paths must pass through $z_1$, which
contradicts the assumption that $z$ is a RBP. Therefore $\pi_1$ and $\pi_2$ must be disjoint until $z_1$, and hence $\pi_2$
follows a different directed edge $(z, z_2)$ in $G_{S, U}$.
\epro

Combing Lemmas \ref{L:nobran} and \ref{L:nocoa}, we obtain the following result on $G_{S, U}$.
\bp\label{P:fgr}
For $S<U$ and $A\subset \Z^2$, let $G_{S, U}(N_\eps(A))$ be defined as before. If $z=(x, t)$ is an $(S, U)$-RBP for $N_\eps(A)$, then there are
exactly two outgoing edges at $z$ in $G_{S, U}$. If $z\in \xi^A_S \times\{S\}$, then there is either one or two outgoing edges at $z$ in $G_{S, U}$.
\ep
\bpro
If $z=(x,t)$ is an $(S, U)$-RBP, then we can find $\pi_1, \pi_2\in N_\eps$ passing through $z$ such that $\pi(s)\neq\pi_2(s)$ for $s\in (t, U]$.
By Lemma~\ref{L:nocoa}, $\pi_1$ and $\pi_2$ follow two different directed edges $(z, z_1)$ and $(z, z_2)$ in $G_{S, U}$. Any other path $\varpi$ passing though $z$ coincides with either $\pi_1$ or $\pi_2$ at time $t+1$, and hence by Lemma~\ref{L:nobran} must follow the same directed edge in $G_{S, U}$.
This implies that there are exactly two outgoing edges at $z$ in $G_{S, U}$. If $z\in \xi^A_S\times \{S\}$, the argument is similar, and there is either one or two outgoing edges at $z$ because there is either one or two branches in $N_\eps$ coming out of $z$.
\epro

\subsection{Probability of having at least $1$ relevant branching point} \label{S:1RBP}

We now consider the branching-coalescing point set $\xi$ generated by the discrete net $N_\eps$, starting from $\xi_0= \{0\}$.
Let $\Ri_T$ denote the number of $(0,T)$-relevant branching points for $N_\eps(0,0)$. We give here an upper bound on $\P(\Ri_T\geq 1)$, or equivalently, $\P(|\xi_T|\geq2)$. We will bootstrap this bound to obtain bounds on $\P(\Ri_T\geq K)$ for $K\geq 2$ in the next subsection.

\bp[At least $1$ RBP]\label{T:1RBP}
Let $(\xi_t)_{t\geq0}$ be the branching-coalescing point set generated by the discrete net $N_\eps$ with $\xi_0=\{0\}$. Let $\Ri_T$ be the number of $(0,T)$-RBP for $N_\eps(0,0)$. Then there exist $C, \delta_0, \eps_0>0$ such that for all $\eps\leq \eps_0$ and $1 \leq T\leq \delta_0\eps^{-2}$,
\be\label{1RBP}
\P(\Ri_T\geq 1) = \P(|\xi_T|\geq 2)  \leq C\eps\sqrt{T}.
\ee
\ep
\bpro[strategy]
The proof starts with the fact that the i.i.d.\ Bernoulli product law with parameter $\rho_\eps=2\eps+O(\eps^2)$ as in \eqref{rhoeps} is invariant for the branching-coalescing point set $(\tilde\xi_t)_{t\geq0}$ generated by the Bernoulli net $\wt N_\eps$ (see Lemma \ref{L:prodinv}), which thanks to the coupling between $N_\eps$ and $\wt N_\eps$ in Lemma \ref{L:Bernbin}, gives an upper bound on the branching-coalescing point set $(\xi_t)_{t\geq 0}$ generated from $N_\eps$ with $\xi_0=\tilde \xi_0$. The stationarity of $\tilde \xi$ implies that for $\xi$, the reduction of particle density up to time $T$ due to coalescence must be greater than the increase of density due to branching. The reduction in density can be computed for the collection of coalescing random walks starting from the same Bernoulli product law. By the reasoning above, this gives an upper bound on the number of offsprings of each particle (with branching turned on for that particle and its offsprings) that differ from all particles in the coalescing system. Restricted to the event that there is initially no other particles nearby in the coalescing system, this then leads to an upper bound on the probability that a single particle has more than 1 offspring at time $T$.
\epro

We now make the above proof sketch rigorous. First we consider the coalescing random walks (CRW) $(\eta_t)_{t\geq0}$ with $\eta_0$ following the i.i.d.\ Bernoulli product law with $\P(0\in \eta_0)=\rho_\eps$ as in \eqref{rhoeps}, and the random walk increment distribution $a(\cdot)$ is the same as that for the discrete net $N_\eps$. Note that there is a natural coupling between $\eta$ and the branching-coalescing point set $\xi$ generated by $N_\eps$, such that almost surely, $\eta_0=\xi_0$, and $\eta_t\subset \xi_t$ for all $t\geq 0$. By translation invariance, the density of $\eta_t$ equals
\begin{equation}
	p_{\eta}(t):=\P(0\in\eta_t) = \P(x\in\eta_t) \quad \mbox{for all } x\in\Z.
\end{equation}
We now bound the decrease of density $p_\eta(0)-p_\eta(t)$ due to coalescence.
\bl[Density reduction for CRW]\label{L:pdec}
Let $(\eta_t)_{t\geq0}$ be the coalescing random walks defined above, starting from an i.i.d.\ Bernoulli product law with density $\rho_\eps$ as in \eqref{rhoeps}. Then there exists a constant $C>0$ such that for any $\eps>0$ and $T>0$,
\begin{equation}\label{eq:coaden}
	%\Delta p_{\eta}(T)=
	p_\eta(0)-p_\eta(T)\leq C\eps^2\sqrt{T}.
\end{equation}
\el
\bpro
For each $x\in\Z$, let $(\pi^x_t)_{t\geq 0}$ denote the coalescing random walk starting from $x$ if $x\in \eta_0$; otherwise let $\pi^x_t\equiv *$ for some cemetery state $*$. Note that for any $x, y\in \Z$, conditioned on $x, y\in \eta_0$, $(\pi^x, \pi^y)$ is a pair of coalescing random walks, which
up to the time they coalesce, has the same distribution as a pair of independent random walks $(S^x, S^y)$ up to their first meeting time $\tau_{x, y}:=\inf\{t\geq 0: S^x_t=S^y_t\}$. Let us fix a total order $\prec$ on $\Z$ with $0$ being the minimal element, and let $\sum_{x\in\Z}$ denote summation w.r.t.\ this total order. We can then bound

\begin{equation}\label{eq:coaden1}
	\begin{aligned}
		\dis p_\eta(T) & =\P(\pi^x_T=0~\text{for some}~x\in\Z)= \sum_{x\in \Z} \E\Big[ 1_{\{\pi^x_T=0\}} 1_{\{\pi^x_T\neq \pi^y_T\, \forall\, y\prec x\}}\Big] \\
		\dis & \geq \sum_{x\in\Z} \Big(\E[ 1_{\{\pi^x_T=0\}}] - \sum_{y\prec x} \E\Big[1_{\{\pi^x_T = \pi^y_T=0\}}\Big]\Big) \\
		\dis & \geq \rho_\eps \sum_{x\in\Z} \P(S^x_T=0) - \rho_\eps^2\sum_{x\neq y} \P(S^y_T=0, \tau_{x,y}\leq T) \\
		\dis & = \rho_\eps- \rho_\eps^2\sum_{z\neq 0}\P(\tau_{0,z}\leq T),
	\end{aligned}
\end{equation}
where in the last step, we made the change of variable $z=y-x$ and used translation invariance. Let $\tilde S:=S-S'$ where $S$ and $S'$ are two independent random walks with the same jump distribution $a(\cdot)$. Then $\tau_{0,z}$ has the same distribution as
$\tilde \tau_{z}$, the first hitting time of $z$ by $\tilde S$ with $\tilde S_0=0$. Therefore,
\begin{equation}\ba{l}
	\dis  \sum_{z\neq 0}\P(\tau_{0,z}\leq T)
	\leq \sum_{z\in\Z}\P(\tilde \tau_z\leq T) = \E\big[\big|\cup_{t\in\N}\{\tilde S_t: t\leq T\} \big| \big]
	\leq 2\E[\sup_{0\leq t\leq T}\widetilde S_T] \leq C\sqrt{T},
	\ea\end{equation}
where the last bound follows from Donsker's invariance principle. Substituting the bound
into \eqref{eq:coaden1} and using $\rho_\eps=2\eps+O(\eps^2)$ then gives \eqref{eq:coaden}.
\epro
\medskip

\bpro[of Proposition~\ref{T:1RBP}] We now turn on branching along the coalescing random walks $\eta$ in Lemma \ref{L:pdec}, namely, independently for each space-time point, we turn it into a branching point with probability $\eps$, where a new random walk is born and makes an independent increment with distribution $a(\cdot)$ satisfying \eqref{passum}. This leads to a branching-coalescing point set $(\xi_t)_{t\geq 0}$ with $\xi_0=\eta_0$, and $\eta_t\subset \xi_t$ a.s.\ for all $t\geq 0$. Using the coupling between the discrete net $N_\eps$ and the Bernoulli net $\wt N_\eps$ in Lemma \ref{L:Bernbin}, we can further couple $(\xi_t)_{t\geq 0}$ with $(\tilde \xi_t)_{t\geq 0}$, the branching-coalescing point set generated from $\wt N_\eps$ with $\tilde \xi_0=\xi_0=\eta_0$. Then $(\tilde \xi_t)_{t\geq 0}$ is a stationary process by Lemma \ref{L:prodinv}, and a.s.\ $\xi_t\subset \tilde \xi_t$ for all $t\geq 0$.

Let $J\in\Z$ denote the site in $\eta_0$ closest to $0$ (pick the positive $J$ if there is a tie). If we replace the coalescing random walk path $\pi^J$
starting from $J$ by branching-coalescing random walks starting from $J$, and denote the resulting particle configuration by $(\bar\eta_t)_{t\geq 0}$, then there is a natural coupling such that almost surely,
$$
\eta_t\subset \bar\eta_t \subset \xi_t \qquad \mbox{for all } t\geq 0.
$$
Since  $\xi_t\subset\tilde\xi_t$, and $\tilde\xi_t$ is stationary with density $\rho_\eps=p_\eta(0)$, we can compare the mean number of particles in $\xi_T\cap [-5\eps^{-1}, 5\eps^{-1}]$ and $\eta_T\cap [-5\eps^{-1}, 5\eps^{-1}]$ and obtain the bound
\begin{equation}\label{eq:Kpeta}
	\begin{aligned}
		(10\eps^{-1}+O(1)) (p_\eta(0)-p_\eta(T)) \geq\ & \E[\xi_T \cap [-5\eps^{-1}, 5 \eps^{-1}]] - \E[\eta_T \cap [-5\eps^{-1}, 5 \eps^{-1}]] \\
		\geq \ & \E[\bar\eta_T \cap [-5 \eps^{-1}, 5 \eps^{-1}]] - \E[\eta_T \cap [-5 \eps^{-1}, 5 \eps^{-1}]] \\
		\geq \ & \P((\bar\eta_T\backslash \eta_T) \cap [-5 \eps^{-1}, 5 \eps^{-1}] \neq \emptyset) =: \P(A).
	\end{aligned}
\end{equation}
If $\xi^J_t$ denotes the configuration of branching-coalescing random walks in the system $\bar\eta$ that start from the single particle at position $J$,
then $\bar\eta_T\backslash \eta_T = \xi^J_T\backslash \eta_T$, and $A$ is just the event that $\xi^J_t$ adds some new particles (compared to $\eta_T$) in $[-5\eps^{-1}, 5\eps^{-1}]$ at time $T$. We will lower bound $\P(A)$ as follows.

First we restrict to the event
$$
A_1:=\big\{ |\eta_0 \cap [-\eps^{-1}, \eps^{-1}]| =  |\eta_0 \cap [-6\eps^{-1}, 6\eps^{-1}]| =1\big\},
$$
namely that initially there is only one particle in $[-\eps^{-1}, \eps^{-1}]$ (whose position is denoted by $J$) and no particle in $[-6\eps^{-1}, 6\eps^{-1}]\backslash [-\eps^{-1}, \eps^{-1}]$. Note that uniformly in $\eps$ small,
\begin{equation}\label{eq:A1}
	\P(A_1) = (1+o(1)) 2\eps^{-1} \cdot \rho_\eps \cdot (1-\rho_\eps)^{12\eps^{-1}} \geq c_0 >0.
\end{equation}

We then further restrict to the event $A_2$ that the coalescing random walk paths $(\pi^x_t)_{t\geq 0}$ starting from all $x\in \eta_0 \backslash [-6\eps^{-1}, 6\eps^{-1}]$ do not enter $[-5\eps^{-1}, 5\eps^{-1}]$ up to time $T$, namely
$$
A_2 := \Big\{  \inf_{x\in \eta_0\backslash [-6\eps^{-1}, 6\eps^{-1}]} \inf_{t\in [0,T]} |\pi^x_T|  > 5\eps^{-1} \Big\} \subset \big\{\eta_T \backslash\{ \pi^J_T\} \cap [-5\eps^{-1}, 5\eps^{-1}]=\emptyset \big\}.
$$
Recall that $T\leq \delta_0 \eps^{-2}$. Then uniformly in $\eta_0\in A_1$, we have
\begin{equation}\label{PA2c}
	\lim_{\delta_0\downarrow 0} \lim_{\eps\downarrow 0} \P(A_2^c|\eta_0) \leq 2 \lim_{\delta_0\downarrow 0} \lim_{\eps\downarrow 0} \P(\sup_{x\leq 0} \sup_{t\leq T} \pi^x_t \geq \eps^{-1}) \leq 4 \lim_{\delta_0\downarrow 0} \lim_{\eps\downarrow 0} \P(\sup_{x\leq 0} \pi^x_{\delta_0\eps^{-2}} \geq \tfrac12 \eps^{-1}) = 0,
\end{equation}
where the second inequality follows by a stopping time argument (the time when one of the $\pi^x$ exceeds $\eps^{-1}$) and Donsker's invariance principle. The last equality above follows from known results on the interface tightness of the dual voter model  \cite{BMV07} and a CLT for the location of the interface \cite{AS11}. More precisely, the last event in \eqref{PA2c} is exactly the event that a dual voter model $\phi_t$ starting with $\phi_0(x) = 1_{\{x\leq 0\}}$ satisfies $\phi_{\delta_0\eps^{-2}}(x)=1$ for some $x\geq 1/2\eps$. {The interface tightness result of \cite{BMV07} shows that the width of the interface between the leftmost $0$ and rightmost $1$ in the configurations $(\phi_t)_{t\geq 0}$ is a tight family of random variables, and hence the interface in the voter configuration $\phi_{\delta_0\eps^{-2}}$ at time $\delta_0\eps^{-2}$ converges to a single point under diffusive scaling of space and time $(y, t) \to (\eps y, \eps^2 t)$ as $\eps\downarrow 0$, while \cite{AS11} shows that the diffusive scaling limit of the location of the interface in the voter configuration $\phi_{\delta_0\eps^{-2}}$ is distributed as a Brownian motion $B_{\delta_0}$ at time $\delta_0$. The last equality above then follows from the fact that $\lim_{\delta_0\downarrow 0} \P(B_{\delta_0} \geq 1/2)=0$. }
Therefore uniformly in $\eta_0\in A_1$, we can make $\P(A_2|\eta_0)$ arbitrarily close to $1$ by choosing $\delta_0>0$ and $\eps\leq \eps_0$ small.

Lastly, we restrict to the event
$$
A_3 :=\big\{ \sup_{t\in [0,T]} |\pi^J_t| < 5\eps^{-1}\} \cap  \big\{\exists\ \tilde \pi^J \in \xi^J \ \text{s.t.}\, \sup_{t\in [0,T]} |\tilde \pi^J_t| < 5\eps^{-1} \ \&\ \tilde \pi^J_T\neq \pi^J_T\big\},
$$
where $\tilde \pi^J \in \xi^J$ means $\tilde \pi^J$ is a path among the collection of branching-coalescing random walks $(\xi^J_t)_{t\geq 0}$. Note that $A_1\cap A_2 \cap A_3 \subset A$, the event defined in \eqref{eq:Kpeta}, and by translating $J$ to $0$, uniformly in $\eta_0\in A_1$, we have
\begin{equation}
	\P(A_3|\eta_0) \geq \P(A_4) :=\P\big(\sup_{t\in [0,T]} |\pi^0_t| < 4\eps^{-1} \ \& \ \exists\ \tilde \pi^0 \in \xi^0 \ \text{s.t.}\sup_{t\in [0,T]} |\tilde \pi^0_t| < 4\eps^{-1}\ \&\ \tilde \pi^0_T\neq \pi^0_T\big).
\end{equation}

Note that conditioned on $\eta_0\in A_1$, the events $A_2$ and $A_3$ are independent since they depend on randomness in disjoint space-time regions. Therefore applying Lemma \ref{L:pdec} to \eqref{eq:Kpeta}, we have
$$
C\eps \sqrt{T} \geq \P(A) \geq \P(A_1) \inf_{\eta_0\in A_1}  \P(A_2|\eta_0) \inf_{\eta_0\in A_1} \P(A_3|\eta_0) \geq \frac{c_0}{2} \P(A_4),
$$
and hence
\begin{equation}
	\P(A_4) \leq C \eps \sqrt{T}.
\end{equation}

Note that $A_4 \subset \{|\xi^0_T|\geq 2\}$. To conclude the proof of Proposition \ref{T:1RBP}, it only remains to prove the same probability bound for the event
$$
A_5:= \{|\xi^0_T|\geq 2\} \backslash A_4 \ \subset \ A_6 \cup A_7,
$$
where
$$
\begin{aligned}
	A_6 & :=\{ \sup_{t\in [0,T]} |\pi^0_t| \geq \eps^{-1}\}, \\
	A_7 & := A_6^c \cap \{ \exists\, \tilde \pi^0 \ {\rm s.t.}\ \tilde \pi^0_T\neq \pi^0_T\} \cap
	\big\{ \forall\, \tilde \pi^0\in \xi^0\ {\rm s.t.}\ \tilde \pi^0_T\neq \pi^0_T, \ \sup_{t\in [0,T]} |\tilde \pi^0_t|\geq 4 \eps^{-1}\big\}.
\end{aligned}
$$
By a standard large deviation estimate for random walks whose increments $a(\cdot)$ have finite $\gamma$-th moment, see e.g.~\cite[Corollary~1.8]{N79} (the continuous-time random walk analogue is provided in \cite[Lemma~5.1]{BMSV06}),
\begin{equation}\label{eq:A6bd}
	%\P(A_6) \leq C\big(e^{-c\,\eps^{-1/2}}+e^{-c/(\eps^{2}T)}+ T\eps^\gamma\big) \leq C \eps \sqrt{T},
	\P(A_6) \leq C\big(e^{-c/(\eps^{2}T)}+ T\eps^\gamma\big) \leq C \eps \sqrt{T},
\end{equation}
where the second inequality is easily seen to hold for $\delta$ and $\eps$ small if $\gamma\geq 2$.

Lastly we show that $\P(A_7)\leq C \eps \sqrt{T}$, which will imply the same bound for $\P(A_5)$ and conclude the proof of Proposition \ref{T:1RBP}.
Note that on the event $A_7$, at least one branching-coalescing random walk in $(\xi^0_t)_{t\geq 0}$ must exit $[-2\eps^{-1}, 2\eps^{-1}]$ before time $T$.
Let $\tau$ denote the stopping time when this happens, and let $A_8$ denote the event that all branching-coalescing random walks in
$(\xi^0_t)_{t\geq 0}$ that exit $[-2\eps^{-1}, 2\eps^{-1}]$ at time $\tau$ actually exit $[-3\eps^{-1}, 3\eps^{-1}]$. Using the graphical construction of $(\xi^0_t)_{t\geq 0}$, where we sample the realisation of the branchings and random walk jumps originating at each space-time point, we note that $A_8$ occurs only if one of the random walk jumps originating from the space-time window $[-2\eps^{-1}, 2\eps^{-1}]\times [0, T]$ jumps outside $[-3\eps^{-1}, 3\eps^{-1}]$. Let $\Delta$ denote a random walk increment. Then by a union bound,
\begin{equation}\label{eq:A8bd}
	\P(A_8) \leq 2 T \sum_{x\geq 0} \P(|\Delta|\geq \eps^{-1}+x) \leq 2 T \sum_{x\geq 0} \frac{\E[|\Delta|^\gamma]}{(\eps^{-1}+x)^\gamma} \leq C \eps^{\gamma-1}T \leq  C\eps\sqrt{T},
\end{equation}
where the last bound holds if $\gamma \geq 3$. On the other hand, we have
\begin{equation}\label{eq:A78bd}
	\P(A_7\cap A_8^c) = \E[ 1_{A_8^c} \, \P(A_7 | (\xi^0_t)_{t\leq \tau}) ] \leq C\eps\sqrt{T},
\end{equation}
because on the event $A_8^c$, there exists some $x\in \xi^0_\tau \cap [-3\eps^{-1}, 3\eps^{-1}]\backslash [-2\eps^{-1}, 2\eps^{-1}]$. We can then sample a random walk $S=(S_t)_{t\geq \tau}$ in $(\xi^0_t)_{t\geq \tau}$ with $S_\tau=x$, whose law is still that of a random walk with transition kernel $a(\cdot)$
since $\tau$ is a stopping time. The event $A_7$ implies that either $S$ returns to the interval $[-\eps^{-1}, \eps^{-1}]$ before time $T$ to ensure $S_T=\pi^0_T$, or $\sup_{t\in [\tau, T]} |S_t|\geq 4\eps^{-1}$. Either way, we must have $\sup_{t\in [\tau, T]} |S_t-S_\tau|\geq \eps^{-1}$, the probability of which is bounded by $C\eps\sqrt{T}$ by \eqref{eq:A6bd}. The bounds \eqref{eq:A8bd} and \eqref{eq:A78bd} imply $\P(A_7)\leq C\eps\sqrt{T}$, which concludes the proof of Proposition
\ref{T:1RBP}.
\epro

\subsection{Probability of having more than $K$ relevant branching points}\label{S:kRBP}

We now bootstrap the bound on $\P(\Ri_T\geq 1)$ in Proposition \ref{T:1RBP} to get a bound on $\P(\Ri_T\geq K)$ for $K\geq 2$.

\begin{figure}[ht!]
	\begin{center}
		\includegraphics[width=0.7\textwidth]{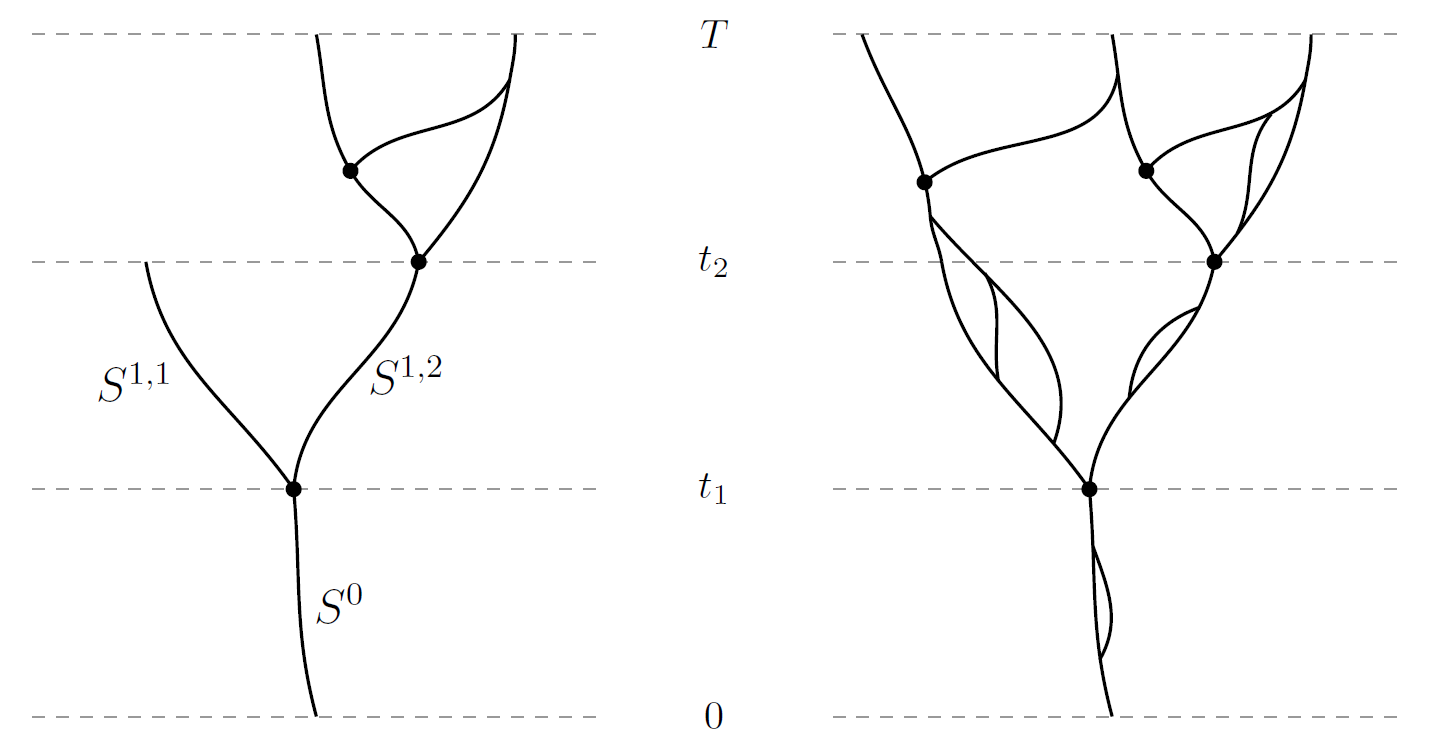}
		\caption{The image on the left illustrates the exploration procedure for $\{\Ri_T\geq2\}$. The image on the right
			illustrates the branching-coalescing point set, with $\bullet$ identifying the relevant branching points (RBP's).}\label{fig:kRBP}
	\end{center}
\end{figure}

\bt[More than $K$ RBP]\label{T:kRBP}
Let $(\xi_t)_{t\geq0}$ be the branching-coalescing point set generated by the discrete net $N_\eps$ with $\xi_0=\{0\}$. Let $\Ri_T$ be the number of $(0,T)$-RBP for $N_\eps(0,0)$. Then for any $K\in\N$, there exist $C_K, \delta_0 =\delta_0(K), \eps_0 =\eps_0(K)>0$ such that for all $\eps\leq \eps_0$
and $1 \leq T\leq \delta_0\eps^{-2}$,
\be\label{nRBP}
\P(\Ri_T\geq K)  \leq C_K\,(\eps\sqrt{T})^K.
\ee
\et
\bpro
We first outline the proof strategy, which starts with a decomposition according to the times $t_1, \ldots, t_K$ of the $K$ earliest $(0,T)$-RBP's
$z_i=(x_i, t_i)$, $1\leq i\leq K$. Assume for the moment the $t_i$'s are all distinct. Let us include half integer times in $(\xi_t)_{t\geq 0}$ to keep track of the increments of the branching-coalescing random walks. We can construct $(\xi_t)_{t\geq 0}$ through an exploration procedure as follows (see Figure~\ref{fig:kRBP}). First sample a random walk $S^0$ starting at $(0,0)$ and run till time $t_1$, where at each branching point of $N_\eps$ the walk encounters, it chooses one of the two outgoing jumps with equal probability. This generates a path in $(\xi_t)_{t\geq 0}$. Then $z_1=(t_1, S^0_{t_1})$ must be the first $(0,T)$-RBP. In particular, $z_1$ is a branching point, an event which we denote by $B_1$ and has probability $\eps$. We next sample two random walks $S^{1,1}$ and $S^{1, 2}$ starting at $z_1$, whose first jumps are independent due to the binary branching at $z_1$, but then interact as coalescing random walks. The paths $(S^{1,1}_t)_{t_1\leq t\leq t_2}$ and $(S^{1,2}_t)_{t_1\leq t\leq t_2}$ are also paths in $(\xi_t)_{t\geq 0}$, and by Lemma \ref{L:nocoa}, $S^{1,1}$ and $S^{1,2}$ cannot coalesce before or at time $t_2$ because $z_1$ is a $(0,T)$-relevant branching point. We denote this event by $E_1$, the probability of which can be bounded by $C/\sqrt{t_2-t_1}$ by standard random walk estimates. At time $t_2$, the next RBP $z_2$ is either $(t_2, S^{1,1}_{t_2})$ or  $(t_2, S^{1,2}_{t_2})$. Once a choice has been made, that point must be a branching point, an event which we denote by $B_2$. We then sample two random walks $S^{2,1}$ and $S^{2,2}$ starting at $z_2$, whose first jumps are independent and then interact as coalescing random walks (also coalesce with previously sampled random walks if they have not terminated as a result of reaching a RBP). Denote by $E_2$ the event that $S^{2,1}$ and $S^{2,2}$ do not coalesce by time $t_3$. We then repeat this construction so that at each time $t_k$, there are at most $k$ choices (the number of locations of the previously sampled coalescing random walks) for the RBP $z_k$, there is an event $B_k$ for $z_k$ to be a branching point, and there is an event $E_k$ for two sampled random walks $S^{k,1}$ and $S^{k,2}$ starting from $z_k$ not to coalesce by time $t_{k+1}$. The event for the time interval $[t_K, T]$ will be different because $t_{K+1}$ is unspecified, and we avoided
a full decomposition of all $(0,T)$-RBP's because we are unable to control the growing combinatorial factor $k$ associated with the choice of $z_k$ at time $t_k$. Instead of restricting to the event $E_K$ on the two sampled random walks $S^{K,1}$ and $S^{K,2}$ starting from $z_K$, we simply restrict to the event $G(z_K; t_K, T)$ that there are two paths in $N_\eps(z_K)$ which are disjoint on the time interval $(t_K, T]$. Thanks to Proposition \ref{T:1RBP}, we can control $\P(G(z_K; t_K, T))$ averaged over $t_K$, which leads to a version of \eqref{nRBP} with $T$ therein replaced by $T'$ averaged over $[T, 2T]$. A monotonicity argument then gives \eqref{nRBP}.

We now make the above proof sketch more precise. First we restrict to the event $D$ that the first $K$ relevant branching points all occur at distinct times. Then by successively conditioning on the sampled coalescing random walks ($S^0$, $S^{1,1}$, $S^{1,2}$, \ldots) up to time $t_K$, $t_{K-1}$, \ldots and applying the Markov property, we have
\begin{equation}\label{DRTK}
	\begin{aligned}
		\P(D \cap \{\Ri_T\geq K\})  & \leq \sum_{0< t_1<\cdots <t_K<T}  \prod_{k=1}^{K-1} (\eps k \P(E_k)) \, \cdot K \, \P(G(z_K; t_K, T)) \\
		& = \eps^{K-1} K! \sum_{0\leq t_1<\cdots <t_K<T} \prod_{i=1}^{K-1} \P(E_k)  \cdot \phi(T-t_K),
	\end{aligned}
\end{equation}
where $\phi(T-t_K):=\P(G(z_K; t_K, T))$ depends only on $T-t_K$ by translation invariance.

To bound $\P(E_k)$, note that if $S$ and $S'$ are two independent random walks with kernel $a(\cdot)$ starting at $(0, 0)$, and $\tau^{S, S'}:=\min\{t\geq 1: S_t=S'_t\}$. Since $S'-S$ is a symmetric random walk with finite second moments, by a standard random walk hitting time estimate (see e.g.~\cite[Lemma 2.2]{NRS05}),
\begin{equation}
	\P(E_k) = \P(\tau^{S, S'} > t_{k+1}-t_k) \leq  \E\Big[C\, \frac{|S_1-S'_1|}{\sqrt{t_{k+1}-t_k}} \Big] \leq \frac{C}{\sqrt{t_{k+1}-t_k}}.
\end{equation}
Substitute this bound into \eqref{DRTK}, replace $T$ by $T'$ and average over $T\leq T'\leq 2T$, we obtain
\begin{equation}\label{SDRTK}
	\begin{aligned}
		\frac{1}{T} \sum_T^{2T} \P(D \cap \{\Ri_{T'}\geq K\})  & \leq \eps^{K-1} K! \Big(\sum_{t=1}^{2T} \frac{C}{\sqrt t}\Big)^{K-1} \sum_{t=1}^{2T} \phi(t) \\
		& \leq K! \, (C\eps \sqrt{T})^{K-1} \sum_{t=1}^{2T} \phi(t),
	\end{aligned}
\end{equation}
where in the first inequality, the sum over $t_1\in (0, T)$ compensates the factor $\frac{1}{T}$. To bound the last sum in \eqref{SDRTK},
note that for $K=1$, we have a similar decomposition according to the time of the first relevant branching point:
$$
\P(\Ri_{2T}\geq 1) = \sum_{t=0}^{2T-1} \P(|\xi_t|=1) \phi(2T-t) \geq \frac{1}{2} \sum_{t=1}^{2T} \phi(t),
$$
where we applied Proposition \ref{T:1RBP} to obtain $\P(|\xi_t|=1)\geq 1/2$. Again by Proposition \ref{T:1RBP},
\begin{equation}
	\sum_{t=1}^{2T} \phi(t) \leq 2 \P(\Ri_{2T}\geq 1) \leq C \eps\sqrt{T}.
\end{equation}
Substituting this bound into \eqref{SDRTK} then gives
\begin{equation}\label{SDRTK2}
	\frac{1}{T} \sum_T^{2T} \P(D \cap \{\Ri_{T'}\geq K\})   \leq K! \, (C\eps \sqrt{T})^K \leq C_K (\eps \sqrt{T})^K,
\end{equation}
which is an averaged version of \eqref{nRBP}.

Before removing the averaging in \eqref{SDRTK2}, let us first remove the restriction $D$ that the first $K$ relevant branching points all occur at distinct times. Assume $D^c$, so that $|\{t_1, t_2, \ldots, t_K\}|=M<K$. Let us relabel these distinct times by $0\leq s_1<\cdots <s_M$, with multiplicities $r_1, \ldots, r_M$ so that $\sum r_i= K$. We can apply the same argument as in the bound for $\P(D\cap \{\Ri_T\geq M\})$, except that at each time $s_k$ with $r_k>1$, we only restrict to the event $E_k$ for a pair of sampled random walks starting from one of the $r_k$ RBP's. However, the $r_k$ RBP's at time $s_k$ must still all be branching points, which contribute a probability factor of $\eps^{r_k}$ instead of $\eps$. This argument gives
\begin{equation}\label{SDRTK3}
	\begin{aligned}
		\frac{1}{T} \sum_T^{2T} \P(D^c \cap \{\Ri_{T'}\geq K\})   & \leq \sum_{M=1}^{K-1} C_{K,M} \eps^{K-M} (C\eps \sqrt{T})^M \\
		& = \sum_{M=1}^{K-1} C_{K,M}  (C\sqrt{T})^{M-K} (C\eps \sqrt{T})^K \leq C_K (\eps \sqrt{T})^K,
	\end{aligned}
\end{equation}
where $C_{K,M}$ are combinatorial factors depending on $K$ and $M$. Therefore \eqref{SDRTK2} can be replaced by
\begin{equation}\label{SDRTK4}
	\frac{1}{T} \sum_T^{2T} \P(\Ri_{T'}\geq K)  \leq C_K (\eps \sqrt{T})^K.
\end{equation}
To remove the averaging in $T'$, note that for any $T'\geq T$, the event $\{\Ri_{T'}\geq K\}$ contains the intersection of the events that $\xi_{T'-T}=:\{x\}$ is a singleton and the branching-coalescing random walks starting from $x$ at time $T'-T$ has at least $K$ relevant branching points for the time interval $(T'-T, T')$. Using the Markov property of $(\xi^\eps_t)_{t\geq 0}$, we then have
\begin{equation}
	\begin{aligned}
		\P(\Ri_{T'}\geq K) \geq \P(|\xi_{T'-T}|=1) \P(\Ri_T\geq K) \geq \frac{1}{2}  \P(\Ri_T\geq K),
	\end{aligned}
\end{equation}
where the last bound follows from Proposition \ref{T:1RBP}. Substituting this bound into \eqref{SDRTK4} then gives \eqref{nRBP}.
\epro

\subsection{A multiscale argument for tightness}\label{S:multi}

We are now ready to prove tightness. 
\bt[Tightness]\label{T:tight}
The family of rescaled discrete nets $(S_{\eps} N_\eps)_{\eps\in(0, 1)}$ in Theorem \ref{T:net} is tight.
\et
We will verify the tightness condition \eqref{tight} for the rescaled discrete nets $S_\eps N_\eps$ by adapting a multiscale argument from \cite{BMSV06} for coalescing random walks. More precisely, we will prove the following:
\begin{equation}
	\lim_{\delta\downarrow0}\frac{1}{\delta}\limsup_{\eps\downarrow0}\sup_{(\eps x, \eps^2 t)\in
		\Lambda_{L,T}} \P\big(N_\eps \in A_{M\eps^{-1},\delta\eps^{-2}}(x,t)\big)=0, \quad\text{for all~}L,T, M>0,
\end{equation}
where $A_{M\eps^{-1},\delta\eps^{-2}}(x,t)$ is defined as in \eqref{tight}. By translation invariance, it suffices to show
\begin{equation}\label{tight2}
	\lim_{\delta\downarrow0}\frac{1}{\delta}\limsup_{\eps\downarrow0} \P\big(N_\eps \in A_{M\eps^{-1},\delta\eps^{-2}}(0,0)\big)=0 \quad \mbox{for all } M>0.
\end{equation}
The key ingredients of the multiscale argument used in \cite{BMSV06} are: density reduction through
coalescence, and large deviation estimates for a single random walk. In our case, we need large deviation estimates for a random
walk and all its branching descendants, which can be achieved thanks to our bounds on the number of relevant branching points (RBP's)
in Theorem \ref{T:kRBP}. This is the content of the following lemma, where the three terms in \eqref{0Tfar} originate respectively from
large deviation from the sum of many small jumps, large deviation from rare large jumps, and excessive branching.

\bl[Large deviation estimate]\label{L:0Tfar}
Let $(\xi_t)_{t\geq0}$ be the branching-coalescing point set generated by the discrete net $N_\eps$ with $\xi_0=\{0\}$.
Then for any $K\in\N$, there exist $c_K, \delta_0, \eps_0>0$ such that for all $\eps\leq \eps_0$, $\ell >0$, and
$1 \leq T\leq\delta_0\eps^{-2}$,
\begin{equation}\label{0Tfar}
	\P\Big(\sup_{0\leq t\leq T} \max \xi_t \geq \ell\Big)\leq e^{-\frac{\ell^2}{c_K T}}+  c_K T \ell^{-\ga}+ c_K \eps^K T^{K/2},
\end{equation}
where the jump distribution $a(\cdot)$ satisfies \eqref{passum} with $\sum_x |x|^\gamma a(x)<\infty$ for some $\gamma>3$.
\el
\bpro
We first prove the deviation bound for $\xi$ at the terminal time $T$, i.e.,
\begin{equation}\label{Tfar}
	\P\big(\max \xi_T \geq \ell \big)\leq e^{-\frac{\ell^2}{c_K T}}+ c_K T\ell^{-\ga} +c_K \eps^K T^{K/2}.
\end{equation}
Let $\Ri_T$ be the number of $(0,T)$-RBP for $N_\eps(0,0)$, defined as in \eqref{RBPdef}. Then by Theorem~\ref{T:kRBP},
\begin{equation}
	\P(\Ri_T\geq K)\leq c_K \eps^K T^{K/2},
\end{equation}
which gives the last term in \eqref{Tfar}.

On the complementary event $\{\Ri_T\leq K-1\}$, we have $|\xi_T|\leq K$ by the properties of RBP established in Section \ref{S:RBP}.
Let us include half integer times in $(\xi_t)_{t\geq 0}$ to keep track of the increments of the branching-coalescing random walks. We can then
sample a random walk $S$ in $\xi$ such that $S_0=0$, and whenever $S$ encounters a branching point of $N_\eps$, it chooses one of the two
outgoing jumps with equal probability. Clearly $S$ has the law of a random walk with jump distribution $a(\cdot)$. On the other hand,
the position $S_T$ depends only on the choices made by $S$ at each $(0, T)$-RBP it encounters along the way. Therefore on the event $\{\Ri_T\leq K-1\}$,
each point in $\xi_T$ (including $\max \xi_T$) will coincide with $S_T$ with probability at least $1/2^{K-1}$, and hence
$$
\P(\max \xi_T \geq \ell, \Ri_T\leq K-1, S_T \geq \ell) \geq \frac{1}{2^{K-1}} \P(\max \xi_T \geq \ell, \Ri_T\leq K-1).
$$
It follows that
\begin{equation}
	\P(\max \xi_T \geq \ell, \Ri_T\leq K-1) \leq 2^{K-1} \P(S_T\geq \ell) \leq e^{-\frac{\ell^2}{c_K T}}+c_K T\ell^{-\ga},
\end{equation}
where the last inequality is a standard random walk large deviation estimate (see e.g.\ \cite[Corollary~1.8]{N79}).

Lastly, to show (\ref{0Tfar}), one only needs to note that
\begin{equation}
	\tau:=\inf\{t\in \N: \max \xi_t \geq \ell\}
\end{equation}
is a stopping time. On the event $\tau\leq T$, we can sample a random walk $S$ in $\xi$ starting from position $\max \xi_\tau$ at time $\tau$, which then leads to $\xi_T \ni S_T\geq \ell$ with a probability that is bounded away from $0$ uniformly in $T-\tau$ (by the central limit theorem). The bound \eqref{0Tfar} then follows from \eqref{Tfar} with an adjusted value of $c_K$.
\epro

{\yu
We now state the density reduction result needed in the proof of Theorem \ref{T:tight}. It extends \cite[Prop.~5.4]{BMSV06} from purely coalescing random walks to branching-coalescing random walks. Roughly speaking, it says that if the initial density of
the random walks is $\upsilon\in (0,1)$ so that the average distance between neighbouring walks is of order $\upsilon^{-1}$, then after evolving for a time $R\upsilon^{-2}$, the density of the branching-coalescing random walks would have decreased by a fixed factor depending on $R$.
}

\bl[Density decay of branching-coalescing random walks]\label{L:denbc}
Let $p\in (1/2, 1)$ be fixed. Let $(\xi_t)_{t\geq 0}$ be a branching-coalescing point set generated by the discrete net $N_\eps$ with initial condition
$\xi_0\subset [0, L]$. Then there exist $R_0, \eps_0, \delta_0, C, c>0$ such that for all $\eps \in (0, \eps_0)$, if
$|\xi_0| \leq \upsilon L$ for some $\upsilon\in (0, 1)$ and $T:= R_0\upsilon^{-2} \leq 2\delta_0 \eps^{-2}$, then
\be\label{denbcbd}
\P( |\xi_T| \geq p\upsilon L) \leq Ce^{-c\upsilon L}+ C\upsilon L\eps^4 T^2.
\ee
\el 

\bpro
By adding new walks if necessary, we can assume that we start with exactly $\upsilon L$ walks at positions $0\leq x_1 <x_2<\cdots <x_{\upsilon L}\leq L$.
Let $\Ri^{x_k}_T$ be the number of $(0, T)$-relevant branching points of $N_\eps(x_k, 0)$, and let $\Ri_T:=\max_{1\leq k\leq \upsilon L}\Ri^{x_k}_T$.
First note that by Theorem~\ref{T:kRBP} with $K=4$, we have
\begin{equation}\label{denbc0}
	\P(\Ri_T\geq 4)\leq C\upsilon L \eps^4 T^2,
\end{equation}
which is the source of the second term in the bound in \eqref{denbcbd}.

Note that for any $L_0\in \N$ (to be chosen later), the number of consecutive pairs $(x_k, x_{k+1})$ with $x_{k+1}-x_k > L_0$ is at most $\frac{L}{L_0}$. Dividing $(x_i, \cdots, x_j, x_{j+1}, \cdots, x_k)$ into two blocks $(x_i, \cdots, x_j)$ and $(x_{j+1}, \cdots, x_k)$ whenever $x_{j+1}-x_j> L_0$, we can divide $(x_1, \ldots, x_{\upsilon L})$ into at most $1+L/L_0$ blocks. In each block, we can match adjacent $(x_i, x_{i+1})$ into pairs, which leaves at most
one unmatched particle. This gives at least $\kappa:=\frac{1}{2}(\upsilon L - 1 - L/L_0)$ disjoint pairs of particles where each pair is within distance $L_0$ of each other. For such a pair $(x, y)$ with $0< y-x \leq L_0$, if we consider the branching-coalescing point sets $\xi^x$ and $\xi^y$  generated by the same discrete net $N_\eps$ with initial conditions $\xi^x_0=\{x\}$ and $\xi^y_0=\{y\}$, then we have
\begin{equation}\label{denbc1}
	\begin{aligned}
		\P(|\xi^x_T\cup\xi^y_T|=1) & \geq \P(\Ri^x_T=\Ri^y_T=0,\xi^x_T = \xi^y_T) \\
		& = \P(\Ri^x_T=\Ri^y_T=0)-\P(\Ri^x_T=\Ri^y_T=0,\xi^x_T\neq \xi^y_T).
	\end{aligned}
\end{equation}
On the event $\{\Ri^x_T=\Ri^y_T=0\}$, $|\xi^x_T| = |\xi^y_T|=1$. Therefore by sampling the random walk jumps at branching points of $N_\eps$, we can couple $\xi^{x}$ and $\xi^y$ with two coalescing random walks $S^x$ and $S^y$ such that $\xi^{x}_T=\{S^x_T\}$ and $\xi^{y}_T=\{S^y_T\}$. This implies that
$\{\xi^x_T\neq \xi^y_T\}=\{S^x_T\neq S^y_T\}$ and
\begin{equation}\label{denbc2}
	\P(\Ri^x_T=\Ri^y_T=0,\xi^x_T\neq \xi^y_T) \leq
	\P(\tau_{x,y}>T)\leq C \frac{L_0}{\sqrt{T}},
\end{equation}
where $\tau_{x, y}$ is the first collision time between $S^x$ and $S^y$, and the last inequality is a standard estimate, see e.g.\ \cite[Lemma 2.2]{NRS05}.

On the other hand, by Theorem~\ref{T:kRBP}, we have
\begin{equation}\label{denbc3}
	\P(\Ri^x_T=\Ri^y_T=0)\geq 1-\sum_{w=x,y}\P(\Ri^w_T\geq 1)\geq 1-2C\eps\sqrt{T}.
\end{equation}
Substituting the above bounds in to \eqref{denbc1} then gives
\be\label{denbc4}
\P(|\xi^x_T\cup\xi^y_T|=1) \geq 1-2C\eps\sqrt{T} - C \frac{L_0}{\sqrt{T}} =: 1-p_0.
\ee

Let us denote the $\kappa =\frac{1}{2}(\upsilon L - 1 - L/L_0)$ pairs of particles by $(u_i, v_i)_{1\leq i\leq \kappa}$, with $v_i \in (u_i, u_i+L_0]$.
We observe that the branching coalescing point set $\xi^{\vec u, \vec v}$ generated by the discrete net $N_\eps$ with $\xi^{\vec u, \vec v}_0=\{u_1, v_1, \ldots, u_\kappa, v_\kappa\}$ can be coupled to $\kappa$ independent branching-coalescing point sets $\xi^{u_i, v_i}$, $1\leq i\leq \kappa$, with $\xi^{u_i, v_i}_0= \{u_i, v_i\}$, such that almost surely,
\be\label{coupling}
|\xi^{\vec u, \vec v}_T| \leq \sum_{i=1}^\kappa |\xi^{u_i, v_i}_T|.
\ee
This is an analogue of \cite[Lemma 5.3]{BMSV06} for coalescing random walks and can be shown as follows. First sample a realisation of the jointly independent branching-coalescing random walks associated with $(\xi^{u_i, v_i})_{1\leq i\leq \kappa}$, keeping track of the random walk jumps
by defining $\xi^{u_i, v_i}_t$ at the half-integer times $t\in \N_0+1/2$. We then construct $\xi^{\vec u, \vec v}$ by recursively defining the branching-coalescing random walks starting from $(u_i, v_i)_{1\leq i\leq \kappa}$. For $i=1$, we construct the branching-coalescing random walks starting from $\{u_1, v_1\}$ by following exactly $\xi^{u_1, v_1}$. We then add random walk paths starting from $\{u_2, v_2\}$ by following paths in $\xi^{u_2, v_2}$, except when a random walk path meets an existing walk, it coalesces with that walk and discards its own descendants. The same procedure is then iterated by adding random walk paths in $\xi^{u_3, v_3}, \ldots, \xi^{u_\kappa, v_\kappa}$. Clearly \eqref{coupling} holds, while the $\xi^{\vec u, \vec v}$ constructed this way has the correct distribution because our procedure recursively explores the realisation of the random walk branching and increments at different lattice sites.

For $\xi$ with $\xi_0=\{x_1, \ldots, x_{\upsilon L}\}$, we can now use \eqref{coupling} to bound
\begin{align}
	\P(|\xi_T|\geq p \upsilon L) & \leq \P(\Ri_T\geq 4) + \P(|\xi_T|\geq p \upsilon L, \Ri_T\leq 3) \notag \\
	& \leq C\upsilon L \eps^4 T^2 + \P\big(|\xi^{\vec u, \vec v}_T| + 4(1+L/L_0) \geq p \upsilon L\big) \notag \\
	& \leq C\upsilon L \eps^4 T^2 + \P\Big(\sum_{i=1}^\kappa |\xi^{u_i, v_i}_T| \geq p \upsilon L - 4(1+L/L_0)\Big), \label{denbc5}
\end{align}
where we used the bound \eqref{denbc0} and the fact that for each $x_i \notin \{u_1, v_1, \ldots, u_\kappa, v_\kappa\}$,
the event $\{\Ri^{x_i}_T\leq 3\}$ implies that $|\xi^{x_i}_T|\leq 4$, and there are at most $1+L/L_0$ many such $x_i$'s.

We now choose $L_0=(\alpha \upsilon)^{-1}$ with $\alpha>0$ small and depending only on $p\in (1/2, 1)$ such that
\begin{align}\label{p'}
	p \upsilon L - 4(1+L/L_0) \geq p'\upsilon L -4, \quad \mbox{with } p'=\frac{1}{2}(p+1/2) > 1/2.
\end{align}
Note that $|\xi^{u_i, v_i}_T|$ are independent random variables, and $|\xi^{u_i, v_i}_T| 1_{\{\Ri^{u_i}_T, \Ri^{v_i}_T\leq 3\}}$
can be dominated by i.i.d.\ random variables $\zeta_i$ with
\begin{align*}
	\P(\zeta_i=1) =1- p_0, \qquad \P(\zeta_i=8) = p_0,
\end{align*}
where as in \eqref{denbc4},
$$
p_0 := 2C\eps\sqrt{T} + C \frac{L_0}{\sqrt{T}} \leq 2C \sqrt{2\delta_0} + \frac{C}{\alpha\upsilon \sqrt{T}}.
$$
We now choose $\delta_0$ sufficiently small and $T= R_0\upsilon^{-2} \leq 2\delta_0 \eps^{-2}$ with $R_0$ sufficiently large such
that $\E[\zeta_i] = 1+7p_0 < 2p'$ with $p'$ defined in \eqref{p'}.

We can then continue the bound in \eqref{denbc5} to obtain
\be
\begin{aligned}
	\P(|\xi_T|\geq p \upsilon L) & \leq C\upsilon L \eps^4 T^2 + \P( \cup_{i=1}^\kappa \{\Ri^{u_i}_T \vee \Ri^{v_i}_T\geq 4\})
	+ \P\Big(\sum_{i=1}^{\gamma L/2} \zeta_i \geq p' \gamma L -4\Big) \\
	& \leq 2 C\upsilon L \eps^4 T^2 + C e^{-c \upsilon L},
\end{aligned}
\ee
where we used $\kappa \leq \upsilon L/2$ and applied a standard large deviation bound for $\sum\zeta_i$. This is exactly
\eqref{denbcbd}, which concludes the proof of Lemma \ref{L:denbc}.
\epro

We are now finally ready to prove the tightness of the family of rescaled discrete nets.
\smallskip

\bpro[of Theorem \ref{T:tight}] 
As discussed above, we need to verify the condition \eqref{tight2} for all $M>0$. As we will see, the constant $M$ enters the proof only through the large deviation estimate in Lemma \ref{L:0Tfar} when we take $\ell=cM\eps^{-1}$ for constants $c>0$ independent of $M$. For any fixed $M>0$,
the last two terms in the right hand side of \eqref{0Tfar} will always be of order $o(\eps^3)$, while the first term will always be of order $o(1/\delta)$ when we take $T=\delta \eps^{-2}$ with $\delta\downarrow 0$. These are the only bounds we need. To simplify the notation, we may therefore assume that $M= 10$ and $\varsigma^2$=1. We also note that to bound $\P(N_\eps \in A_{M\eps^{-1}, \delta \eps^{-2}}(x,t))$ in \eqref{tight2}, we can discard the event that there are random walk increments in $N_\eps$ that jump across the interval $[-M\eps^{-1},M\eps^{-1}]$ in the time interval $[0,\delta\eps^{-2}]$, so that we only need to restrict our attention to random walks that originate inside the space-time box $[-M\eps^{-1}, M\eps^{-1}]\times[0,\delta\eps^{-2}]$. The same issue arises in the proof of tightness for the family of rescaled coalescing random walks, and the proof is identical since the branching at most doubles the number of random walk increments, see \cite[Section 4]{NRS05}.

We now apply the multiscale argument as in the proof of \cite[Proposition~2.3]{BMSV06}, where for convenience we actually consider all random walks originating from a larger space-time box $[-M\eps^{-1}, M\eps^{-1}]\times[0,2\delta\eps^{-2}]$. {\yu The basic strategy is to first let random walks starting
from inside the space-time box evolve until a time of the form $iR$ for a suitably chosen $R>0$ and some $i\in\N$. In this step, we use Lemma \ref{L:0Tfar} to  show that with high probability, no random walk travels very far. For each $i\in\N$, we then consider the resulting random walks at times $(2i-2)R$ and $(2i-1)R$ and let them evolve together till time $2iR$. We will show that with high probability, coalescence leads to a density reduction when these random walks reach time $2iR$ (by Lemma \ref{L:denbc}), and no random walk travels very far (by Lemma \ref{L:0Tfar}). This argument is then repeated until the only random walks remaining are at the final time $2\delta\eps^{-2}$. The tightness condition \eqref{tight2} is verified by showing that the aggregate probability
of all the bad events are small enough. 
}

{\yu \bf Step 1.} We introduce a first time scale by choosing a large fixed constant $R>0$ independent of $\eps>0$ (to be determined later using Lemma \ref{L:denbc}) and decompose the space-time box $[-M\eps^{-1},M\eps^{-1}]\times[0,2\delta\eps^{-2}]$ in the time direction into $I_0:=2\delta\eps^{-2}/R$ boxes of the form
\be\label{B0i}
B_{0,i}:=[-M\eps^{-1},M\eps^{-1}]\times[(i-1)R,iR], \qquad i=1,\ldots,I_0.
\ee
To simplify notation later, we may assume $I_0=2^J$ for some $J\in \N$; otherwise we can increase $\delta$ and $I_0$ by a factor of at most $2$ such that this property holds.
For the branching-coalescing random walks in $N_\eps$ starting from each space-time box $B_{0, i}$, we let them evolve up to time  $T^{(0)}_i:=iR$ and discard the event $E_{0, i}$ that some particle exits the interval
$$
[-\wt M_0, \wt M_0] :=[-(M+1)\eps^{-1},(M+1)\eps^{-1}]
$$
before time $T^{(0)}_i$. Indeed, by the large deviation estimate in Lemma~\ref{L:0Tfar} with $K=4$ and a union bound over all starting points in $\cup_{i=1}^{I_0} B_{0, i}$,
we have
\begin{equation}\label{multi1}
	\begin{aligned}
		\dis\P(\cup_{i=1}^{I_0} E_{0, i})
		& \dis\leq 4M\delta \eps^{-3} \big( e^{-\frac{cM^2\eps^{-2}}{R}}+ cR M^{-\gamma}\eps^\gamma +c\eps^4 R^2 \big),
	\end{aligned}
\end{equation}
where $\gamma=3+\eta >3$ by the assumption on the jump distribution $a(\cdot)$ in \eqref{passum}. Denote $E_0:=\cup_{i=1}^{I_0} E_{0, i}$.
Given $M$ and $R$, clearly
\begin{equation}\label{multi1.5}
	\lim_{\delta\downarrow 0} \frac{1}{\delta} \limsup_{\eps\downarrow 0} \dis\P(E_0) = 0.
\end{equation}
Therefore to verify \eqref{tight2}, we can discard the event $E_0$, which ensures that no walk has exited
$[-\wt M_0, \wt M_0]$ in this first step, and at the end of this step, all remaining walks lie in $[-\wt M_0, \wt M_0]$
at times $T^{(0)}_i=iR$, $1\leq i\leq I_0$, whose configurations we denote by $\xi^{(0)}_{i} \subset [-\wt M_0, \wt M_0]$.

{{\bf Step 2.} We now consider the second time scale.} Let $I_1:= I_0/2$ and let $T^{(1)}_i := T^{(0)}_{2i} = 2iR$ for $1\leq i\leq I_1$. For each $1\leq i\leq I_1$, we consider the remaining random walks starting from $\xi^{(0)}_{2i-2} \times \{T^{(0)}_{2i-2}\}$
and $\xi^{(0)}_{2i-1} \times \{T^{(0)}_{2i-1}\}$ and let them evolve together till time $T^{(1)}_i=T^{(0)}_{2i}$, and denote the resulting particle configuration at time $T^{(1)}_i$ by $\xi^{(1)}_{i}$. We will show that we can also discard the events $E_{1,i}:=D_{1,i}\cup F_{1,i}$, $1\leq i\leq I_1$, where
\begin{itemize}
	\item [$D_{1,i}$:] Some random walk starting from $\xi^{(0)}_{2i-2}$ at time $T^{(0)}_{2i-2}$, or from $\xi^{(0)}_{2i-1}$ at time $T^{(0)}_{2i-1}$, exits the interval
	$[-\wt M_1,\wt M_1]$ before time $T^{(1)}_i=T^{(0)}_{2i}$, where $\wt M_1:= \wt M_0 + \beta_1\eps^{-1}$ for some $\beta_1$ to be chosen later.
	
	\item [$F_{1,i}$:]  The particle density $\frac{1}{4M \eps^{-1}} |\xi^{(1)}_i| > p:=\frac{1}{\sqrt{2}}$.
\end{itemize}
Let us denote $E_1 := \cup_{i=1}^{I_1}E_{1,i}$. Similar to \eqref{multi1.5}, we will show that
\begin{equation}\label{multi2}
	\lim_{\delta\downarrow 0} \frac{1}{\delta} \limsup_{\eps\downarrow 0} \dis\P(E_1 \cap E_0^c) = 0,
\end{equation}
which allows us to discard the event $E_1$ in the verification of \eqref{tight2}.

{{\bf Step 3.} We now iterate the argument above and inductively define what is $E_j$ for $j\geq 2$.} After defining and discarding the events
$E_0\cup \cdots \cup E_{j-1}$ for some $j\geq 2$ as above, we are left with random walk configurations $\xi^{(j-1)}_i$ at times $T^{(j-1)}_i= T^{(j-2)}_{2i} = i\cdot 2^{j-1} R$ for $1\leq i \leq I_{j-1}:= 2^{-(j-1)} I_0$, such that
\be\label{xijbdd}
\xi^{(j-1)}_i \subset [- \wt M_{j-1}, \wt M_{j-1}] \qquad \mbox{and} \qquad \frac{1}{4M\eps^{-1}} |\xi^{(j-1)}_{i}| \leq p^{j-1} = 2^{-(j-1)/2},
\ee
where $\wt M_{j-1} =\wt M_0 + (\beta_1 +\cdots +\beta_{j-1})\eps^{-1}$. Note that we choose $p=1/\sqrt{2}$ because by diffusive scaling, this is exactly the rate at which the density of coalescing random walks on $\Z$ decays when time is doubled. This allows us to show inductively in $j$ that the events in \eqref{xijbdd} occur with high probability.

\begin{figure}[ht!]
	\begin{center}
		\includegraphics[width=1\textwidth]{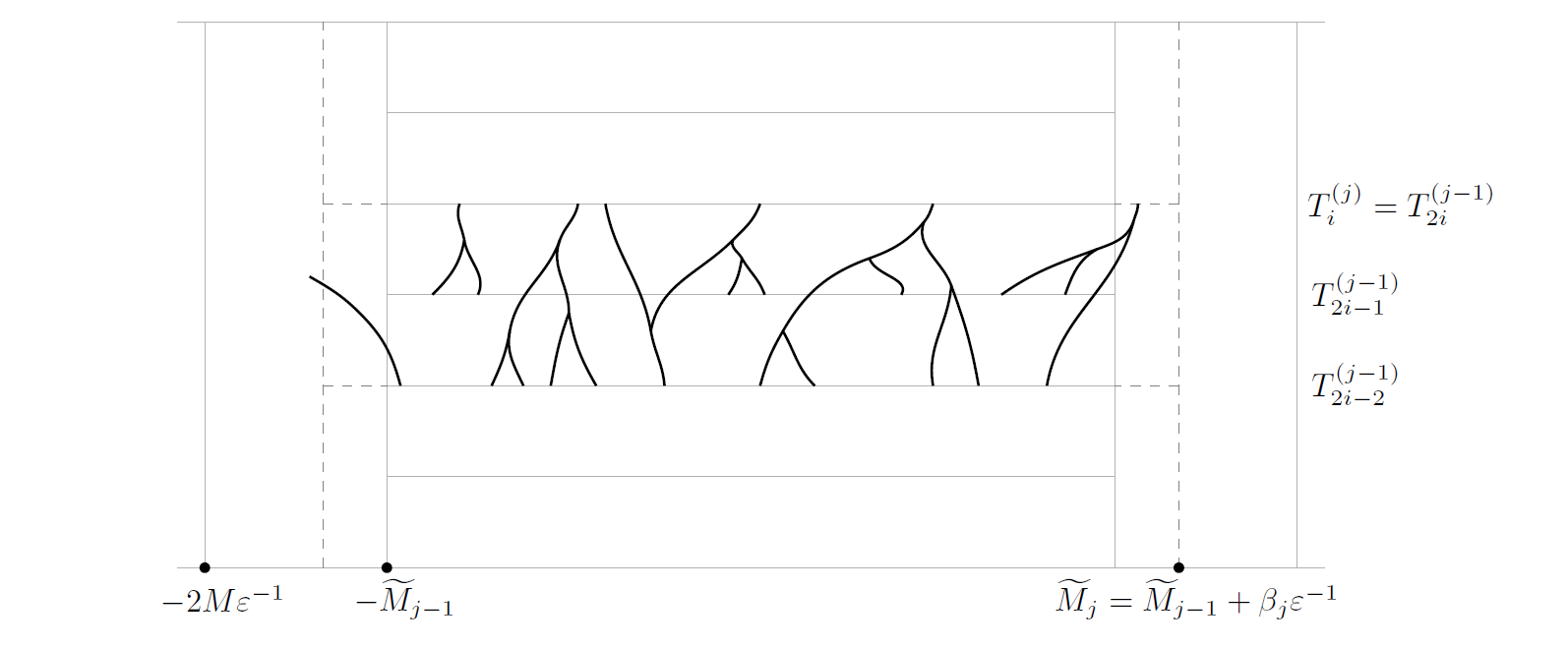}
		\caption{Illustration of the $j$-th iteration of the multiscale argument.}\label{fig:multiscale}
	\end{center}
\end{figure}

To apply the $j$-th iteration of our inductive argument (see Figure \ref{fig:multiscale}), for each $1\leq i \leq I_j:= 2^{-j}I_0$, we consider the remaining random walks starting from $\xi^{(j-1)}_{2i-2} \times \{T^{(j-1)}_{2i-2}\}$
and $\xi^{(j-1)}_{2i-1} \times \{T^{(j-1)}_{2i-1}\}$ and let them evolve together till time $T^{(j)}_{i}:= T^{(j-1)}_{2i} =i\cdot 2^jR$, and denote the resulting particle configuration at time $T^{(j)}_i$ by $\xi^{(j)}_i$. For $1\leq i\leq I_j$, define the events
\begin{itemize}
	\item [$D_{j,i}$:] Some random walk starting from $\xi^{(j-1)}_{2i-2}$ at time $T^{(j-1)}_{2i-2}$, or from $\xi^{(j-1)}_{2i-1}$ at time $T^{(j-1)}_{2i-1}$, exits the interval
	$[-\wt M_j,\wt M_j]$ before time $T^{(j)}_i=T^{(j-1)}_{2i}$, where $\wt M_j:= \wt M_{j-1} + \beta_j\eps^{-1}$ for some $\beta_j>0$ to be chosen.
	
	\item [$F_{j,i}$:]  The particle density $\frac{1}{4M\eps^{-1}} |\xi^{(j)}_i| > p^j=2^{-j/2}$.
\end{itemize}
As before, we denote $E_{j, i} := D_{j, i} \cup F_{j, i}$, and $E_j := \cup_{i=1}^{I_j}E_{j, i}$.

{{\bf Step 4.} We now collect all the bad events and control their probability via suitable choices of the parameters $(\beta_j)_{1\leq j\leq J}$ and $R$.}
We first focus on the choice of $(\beta_j)_{1\leq j\leq J}$. By the assumption after \eqref{B0i}, $I_0= 2^J$ for some $J\in \N$. After $J$ iterations, all walks starting from the box $[-M\eps^{-1}, M\eps^{-1}]\times[0,2\delta\eps^{-2}]$ have evolved till time $2\delta \eps^{-2}$. To verify \eqref{tight2} and bound the probability of the event
that some random walk exits the interval $[-2M \eps^{-1}, 2M\eps^{-1}]$ before time $2\delta \eps^{-2}$, it suffices to choose $(\beta_j)_{1\leq j\leq J}$
such that $\wt M_J < 2 M \eps^{-1}$, and we can discard $\cup_{j=0}^J E_j$ in the sense that
$$
\lim_{\delta\downarrow 0} \frac{1}{\delta} \limsup_{\eps\downarrow 0} \P\big(\cup_{j=0}^J E_{j}\big) = 0,
$$
which would follow if we show that
\begin{equation}\label{multi2.5}
	\lim_{\delta\downarrow 0} \frac{1}{\delta} \limsup_{\eps\downarrow 0} \dis \sum_{j=0}^J \P\big(E_{j} \bigcap \cap_{k=0}^{j-1} E_k^c\big) = 0.
\end{equation}
We will verify \eqref{multi2.5} for $\beta_j:=\frac{1}{(J+1-j)^2}$, a choice that ensures
$$
\wt M_J = \wt M_0 + \eps^{-1} \sum_{j=1}^J \beta_j < 2M \eps^{-1}.
$$

For $1\leq j\leq J$, we can write $E_j= D_j \cup F_j$, where $D_j := \cup_{i=1}^{I_j}D_{j, i}$ and $F_j := \cup_{i=1}^{I_j}F_{j, i}$. Since \eqref{multi1.5} already
gives the desired bound on $\P(E_0)$, to show \eqref{multi2.5}, we will bound  $\P(D_j \cap \cap_{k=0}^{j-1} E_k^c)$ and $\P(F_j \cap \cap_{k=0}^{j-1} E_k^c)$ for $1\leq j\leq J$.

We first consider each event $D_{j, i} \subset D_j= \cup_{i=1}^{I_j}D_{j, i}$, which concerns the random walks starting from $\xi^{(j-1)}_{2i-2}\times \{T^{(j-1)}_{2i-2}\}$ and $\xi^{(j-1)}_{2i-1}\times \{T^{(j-1)}_{2i-1}\}$. Restricted to the event $\cap_{k=0}^{j-1} E_k^c$, the particle configurations $\xi^{(j-1)}_{2i-2}$ and $\xi^{(j-1)}_{2i-1}$ satisfy the conditions in \eqref{xijbdd}. For the event $D_{j, i}$ to occur, one of the random walks starting from either $\xi^{(j-1)}_{2i-2}$ at time $T^{(j-1)}_{2i-2}$, or from $\xi^{(j-1)}_{2i-1}$ at time $T^{(j-1)}_{2i-1}$, must travel a distance of at least $\beta_j \eps^{-1}$ within a time interval of length $2^j R$. Therefore by conditioning on $\xi^{(j-1)}_{2i-2}$, resp., $\xi^{(j-1)}_{2i-1}$, and performing a union bound and applying Lemma~\ref{L:0Tfar} with $\ell =\beta_j\eps^{-1}$ and $T=2^j R$, we obtain
\begin{align}
	\P\big(D_j \bigcap \cap_{k=0}^{j-1} E_k^c\big)  & \leq \sum_{i=1}^{I_j} \P\big(D_{j, i} \bigcap \cap_{k=0}^{j-1} E_k^c\big) \notag \\ 	
	& \leq I_j 2^{2-\frac{j-1}{2}} M \eps^{-1}
	\Big( \exp\Big\{-\frac{c\beta_j^2\eps^{-2}}{2^j R}\Big\}+ c2^jR(\beta_j\eps^{-1})^{-\gamma} +c\,\eps^4 2^{2j}R^2 \Big) \notag \\
	& \leq 2^{3-\frac{3}{2}(j-1)} \frac{M\delta}{R\,\eps^3} \Big( \exp\Big\{-\frac{c\beta_j^2\eps^{-2}}{2^j R}\Big\}+ c2^jR(\beta_j\eps^{-1})^{-\gamma} +c\,\eps^4 2^{2j}R^2 \Big) \notag\\
	& =: 8(B^1_j + B^2_j + B^3_j),  \label{multi3}
\end{align}
where we have expanded the product.

Note that
$$
\sum_{j=1}^J B^3_j = c \sum_{j=1}^J 2^{\frac{j+3}{2}} MR \delta \eps \leq C 2^{\frac{J}{2}} MR \delta  \eps \leq C \sqrt{\frac{\delta}{\eps^2 R}} MR \delta  \eps = C M R^{\frac{1}{2}} \delta^{\frac{3}{2}},
$$
where we used that $I_0=2\delta\eps^{-2}/R =2^J$. Therefore
\be\label{B3j}
\limsup_{\delta\downarrow 0} \frac{1}{\delta} \limsup_{\eps \downarrow 0} \sum_{j=1}^J B^3_j =0.
\ee

Also note that because $\beta_j\geq \beta_1 = J^{-2}$,
$$
\sum_{j=1}^J B^2_j \leq  C \sum_{j=1}^J 2^{-\frac{j}{2}} \beta_j^{-\gamma} M \delta \eps^{\gamma -3} \leq C M \delta \eps^{\gamma -3} \beta_1^{-\gamma}  \sum_{j=1}^J 2^{-\frac{j}{2}} \leq C M \delta \eps^{\gamma -3} J^{2\gamma}.
$$
Since $J= \log_2 I_0=\log_2 \frac{2\delta\eps^{-2}}{R} =o(\eps^{\frac{3-\gamma}{2\gamma}})$ as $\eps\downarrow 0$ because $\gamma>3$, we have
\be\label{B2j}
\limsup_{\delta\downarrow 0} \frac{1}{\delta} \limsup_{\eps \downarrow 0} \sum_{j=1}^J B^2_j =0.
\ee

Lastly using $\beta_j:=\frac{1}{(J+1-j)^2}$ and $2^J=I_0= \frac{2\delta\eps^{-2}}{R}$, we note that
\begin{align*}
	\sum_{j=1}^J B^1_j & \leq  C \frac{M\delta}{R\, \eps^3} \sum_{j=1}^J 2^{-\frac{3 j}{2}} e^{-\frac{c\beta_j^2}{2^jR \eps^2}} \\
	& = C \frac{M\delta}{R\, \eps^3} \sum_{j=1}^J 2^{-\frac{3 j}{2}} e^{-\frac{c}{2^j(J+1-j)^2R \eps^2}} = C \frac{M\delta}{R\, \eps^3} 2^{-\frac{3 (J+1)}{2}}
	\sum_{\tilde j=1}^J 2^{\frac{3 \tilde j}{2}} e^{-\frac{c}{2^{J+1-\tilde j}\, \tilde j^2R \eps^2}} \\
	& \leq C M R^{\frac{1}{2}} \delta^{-\frac{1}{2}} \sum_{\tilde j=1}^J 2^{\frac{3 \tilde j}{2}} e^{-c  \cdot \frac{2^{\tilde j-1}}{\tilde j^2} \cdot \frac{1}{\delta} }
	\leq C M R^{\frac{1}{2}} \delta^{-\frac{1}{2}} e^{-\frac{c}{2\delta}} \sum_{\tilde j=1}^J 2^{\frac{3 \tilde j}{2}} e^{-c  \big( \frac{2^{\tilde j-1}}{\tilde j^2}-\frac{1}{2}\big) \frac{1}{\delta} },
\end{align*}
where we note that the last sum is uniformly bounded in $J$ and in $\delta>0$ sufficiently small, because $c$ comes from Lemma \ref{L:0Tfar} and does not depend on $M, R, \delta$ and $\eps$. It follows that
\be\label{B1j}
\limsup_{\delta\downarrow 0} \frac{1}{\delta} \limsup_{\eps \downarrow 0} \sum_{j=1}^J B^1_j =0.
\ee
Combining \eqref{B3j}-\eqref{B1j} gives
\begin{equation}\label{multi3.5}
	\lim_{\delta\downarrow 0} \frac{1}{\delta} \limsup_{\eps\downarrow 0} \dis \sum_{j=0}^J \P\big(D_{j} \bigcap \cap_{k=0}^{j-1} E_k^c\big) = 0,
\end{equation}
which accounts for the contribution of $D_j \subset E_j$ in \eqref{multi2.5}.

{{\bf Step 5.} We are now ready to choose $R$ in \eqref{B0i}, which will allow us to bound the probability of the events $F_j= \cup_{i=1}^{I_j}F_{j, i}$ on the density of the branching-coalescing random walks.} Let $R_0$ be as in Lemma \ref{L:denbc} for $p=1/\sqrt{2}$, and set $R :=7R_0$.

\begin{figure}[ht!]
	\begin{center}
		\includegraphics[width=1\textwidth]{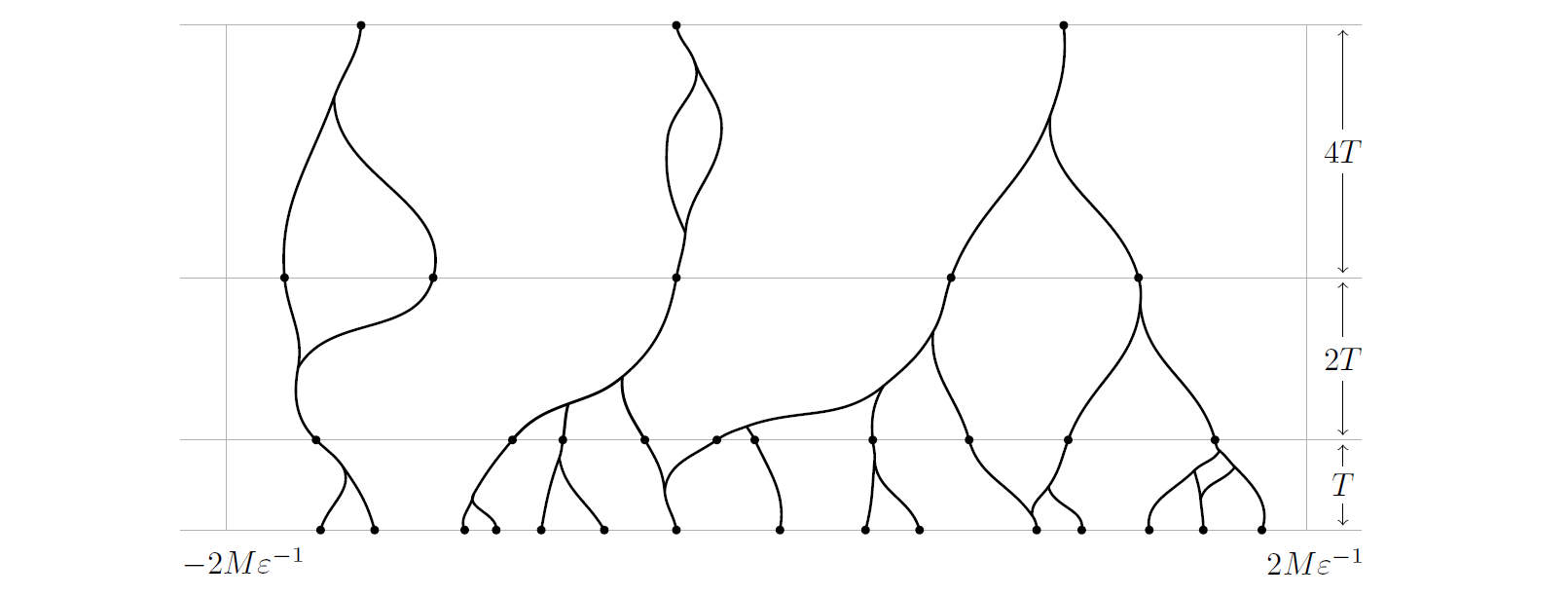}
		\caption{For random walks starting from $\xi^{(j-1)}_{2i-1}$ at time $T^{(j-1)}_{2i-1}$, Lemma~\ref{L:denbc} is applied three times with $p=1/\sqrt{2}$ and successive time durations $T=2^{j-1}R_0$, $2T$, and $4T$.}\label{fig:dendec}
	\end{center}
\end{figure}

Assuming Lemma \ref{L:denbc} for the moment, we now bound $\P(F_{j, i} \cap \cap_{k=0}^{j-1} E_k^c)$ for $1\leq j\leq J$ and $1\leq i\leq I_j$.
Note that restricted to the event $\cap_{k=0}^{j-1} E_k^c$, the particle configurations $\xi^{(j-1)}_{2i-2}$ and $\xi^{(j-1)}_{2i-1}$ satisfy the conditions in \eqref{xijbdd}.  To bound $\P(F_{j, i} \cap \cap_{k=0}^{j-1} E_k^c)$, we first apply Lemma \ref{L:denbc} with
$$
p=1/\sqrt{2}, \quad  L=4 M \eps^{-1}, \quad \upsilon=2^{-\frac{j-1}{2}}=p^{j-1}, \quad T=R_0\upsilon^{-2} = 2^{j-1}R_0
$$
to the random walks starting from $\xi^{(j-1)}_{2i-2}$ at time $T^{(j-1)}_{2i-2}$, which implies that
\begin{align}
	\P\big( |\xi_{2^{j-1}R_0}| \geq p^j L  \,\big|\, \xi_0 = \xi^{(j-1)}_{2i-2}\big)
	& \leq  Ce^{-cp^{j-1} L}+Cp^{j-1}L \eps^4 2^{2(j-1)} R_0^2 \notag \\
	& \leq  Ce^{-cp^j \eps^{-1}}+ C p^{-3j} \eps^3, \label{denpr1}
\end{align}
where we have adjusted the values of $C$ and $c$ to absorb factors of $M, R_0$ and $p$.
Thanks to \eqref{multi3.5}, we may restrict to the event that the random walks starting from $\xi^{(j-1)}_{2i-1}$ at time $T^{(j-1)}_{2i-1}$
is contained in $(-2M \eps^{-1}, 2M \eps^{-1})$ for a time duration of $2^{j-1} R=2^{j-1}\cdot 7 R_0$. Therefore, after running the walks for
a time duration of $T=2^{j-1}R_0$, which reduces the density by a factor of $p$, we can evolve the random walks further for a time duration of $2T$ and
then $4T$ and apply Lemma \ref{L:denbc} each time (see Figure~\ref{fig:dendec}), so that the density of the walks in $(-2M \eps^{-1}, 2M\eps^{-1})$ at time $T^{(j-1)}_{2i}=T^{(j-1)}_{2i-1}+ 7T= i\cdot 2^j R$ is bounded by $p^{j+2}$ (the probability of the complementary event is bounded by the r.h.s.\ of \eqref{denpr1} with adjusted values of $C$ and $c$).
Similarly, for the walks starting at $\xi^{(j-1)}_{2i-2}$ at time $T^{(j-1)}_{2i-2}$, we can first evolve them till time $T^{(j-1)}_{2i-1}$ and
then till time $T^{(j)}_i = T^{(j-1)}_{2i}$, and we can apply Lemma \ref{L:denbc} to conclude that the density of these walks is also bounded by $p^{j+2}$ with high probability. Taking a union bound over the walks starting from $\xi^{(j-1)}_{2i-2}\times T^{(j-1)}_{2i-2}$ and from $\xi^{(j-1)}_{2i-1}\times T^{(j-1)}_{2i-1}$, we conclude that the walks that remain at time $T^{(j)}_i$
have density at most $2p^{j+2} = p^j$ in the interval $(-2M\eps^{-1}, 2M\eps^{-1})$, which is exactly the event $F_{j, i}^c$.
Therefore using the bound in \eqref{denpr1}, we have
\begin{align}\label{multi3.6}
	\sum_{j=1}^J \P\big(F_j \bigcap \cap_{k=0}^{j-1} E_k^c\big) & \leq C \sum_{j=1}^J I_j \big(e^{-cp^j \eps^{-1}} \!\!\! + p^{-3j} \eps^3\big)
	\leq \frac{C \delta}{R \eps^2} \sum_{j=1}^J  2^{-j} \big(e^{-cp^j \eps^{-1}} \!\!\! + p^{-3j} \eps^3\big).
\end{align}
Writing the above as two sums and using $2^J=I_0= \frac{2\delta\eps^{-2}}{R}$, we can bound the first sum by
\begin{align*}
	& \frac{C \delta}{R \eps^2} \sum_{j=1}^J  2^{-j} e^{-cp^j \eps^{-1}}  = \frac{C \delta}{R \eps^2} \sum_{j=1}^J  p^{2j} e^{-cp^j \eps^{-1}} \\
	\leq\ &  \frac{C \delta}{R \eps^2} \cdot p^{2J} \sum_{\ti j=1}^J  p^{-2\ti j} e^{-cp^{J-\ti j} \eps^{-1}}
	\leq C e^{-c p^J \eps^{-1}} \sum_{\ti j=1}^J  2^{\ti j} e^{-c p^J \eps^{-1}(p^{-\ti j} -1)} \\
	\leq & C e^{-\frac{c}{\sqrt \delta}}  \sum_{\ti j=1}^J  2^{\ti j} e^{-\frac{c}{\sqrt \delta}(2^{\ti j/2} -1)},
\end{align*}
where the last sum is uniformly bounded in $J$ and in $\delta>0$ small.

The second sum from \eqref{multi3.6} can be bounded by
\begin{align*}
	\frac{C \delta}{R \eps^2} \sum_{j=1}^J  2^{-j} p^{-3j} \eps^3 = \frac{C \delta \eps}{R} \sum_{j=1}^J p^{-j} \leq \frac{C \delta \eps}{R} \cdot 2^{\frac{J}{2}} \leq \frac{C \delta \eps}{R} \frac{\sqrt{\delta}}{\eps \sqrt{R}} \leq C \delta^{3/2}.
\end{align*}
Substituting the above bounds into \eqref{multi3.6} then gives
\be\label{multi4}
\limsup_{\delta\downarrow 0} \frac{1}{\delta} \limsup_{\eps\downarrow 0} \sum_{j=1}^J \P\big(F_j \bigcap \cap_{k=0}^{j-1} E_k^c\big) =0.
\ee
Together with \eqref{multi3.5}, this concludes the proof of \eqref{multi2.5} and Theorem \ref{T:tight}.
\epro

%%%%%%%%%%%%%%%%%%%%%%%%%%%%%%%%%%%%%%%%%%%%%%
%% Single Appendix:                         %%
%%%%%%%%%%%%%%%%%%%%%%%%%%%%%%%%%%%%%%%%%%%%%%
\begin{appendix}
\section*{Proof of some facts for the Brownian web and net}%\label{S:App}
In this appendix, we sketch the proof of Proposition \ref{P:equiv} and Theorem \ref{T:Nhop}, which follow
readily from results in the literature.
\medskip

\noindent
\bpro[sketch for Proposition \ref{P:equiv}]
As in \eqref{a2}, let us consider a pair jump distribution
$$
a^{(2)}(x_1,x_2):= (1-2\theta\eps)1_{\{x_1=x_2\}}a(x_1)+2\theta\eps a(x_1)a(x_2)~~\mbox{with}~~ a(x)=\frac{1}{2}(1_{\{x=1\}} + 1_{\{x=-1\}}).
$$
Since the simple random walk has period $2$, we consider the pair of discrete webs $(W^1, W^2)$ defined as before \eqref{a2}, but with paths starting only from the even lattice $\Z^2_{\rm even}=\{(x,t):x+t~\text{is even}\}$. Since $\varsigma^2= \sum_x x^2 a(x)=1$, we write $S_\eps$ for the scaling map $S_\eps$. We will show that $(S_{\eps}W^1,S_{\eps}W^2)$ converges to a pair of sticky Brownian webs $(\Wi^1,\Wi^2)$ with parameter $\theta$ both in the sense of marking construction and in the sense of martingale characterisation, which then yields the equivalence of the two characterisations.

The convergence $(S_{\eps}W^1,S_{\eps}W^2)$ to $(\Wi^1,\Wi^2)$ with parameter $\theta$ is a slight extension of \cite[Theorem~6.15]{SSS14} and
follows the same proof. Assume $\theta=1$ to simplify notation. Theorem 6.15 of \cite{SSS14} establishes the convergence of $(S_\eps W^1,S_\eps N_\eps)$, the rescaling of a nearest neighbour discrete web $W^1$ coupled to a nearest neighbour discrete net $N_\eps$, to a coupled Brownian web and net $(\Wi, \Ni)$.  The discrete
web $W^1$ is obtained from $N_\eps$ by sampling at branching points, and the pair $(\Wi, \Ni)$ is constructed in two ways: either start with $\Wi$ and
construct $\Ni$ by a marking construction similar to the marking construction of the sticky webs, or start with $\Ni$ and construct $\Wi$ by sampling at the so-called relevant separation points of $\Ni$. The key step in the proof of \cite[Theorem~6.15]{SSS14} is that the relevant branching points of $N_\eps$ converges to the relevant separation points of $\Ni$, while the sampling distribution at relevant branching points of $N_\eps$ converges
to the sampling distribution at the relevant separation points of $\Ni$. As noted in (7.23) in \cite{SSS14}, the convergence of $(S_\eps W^1,S_\eps N_\eps)\Rightarrow (\Wi, \Ni)$ can be extended to include multiple sticky discrete webs sampled from the same discrete net $N_\eps$. In our case,
this means $(S_\eps W^1, S_\eps W^2, S_\eps N_\eps)\Rightarrow (\Wi^1, \Wi^2, \Ni)$, where the discrete net $N_\eps$ is generated from the pair of discrete webs $(W^1, W^2)$ by switching between paths in $W^1$ and $W^2$. This is equivalent to starting from a discrete net $N_\eps$ and then sampling the random walk increments in $(W^1, W^2)$ with a sampling distribution such that at each branching point of $N_\eps$, $W^1$ samples one of the two outgoing arrows in $N_\eps$ with equal probability while $W^2$ always samples the other arrow. The same sampling distribution can be defined at the relevant separation points of $\Ni$, which in the notation of \cite[Lem.~6.16]{SSS14}, corresponds to parameters $p_{--}=p_{++}=0$, $p_{-+}=p_{+-}=1/2$. The limit $(\Wi^1, \Wi^2)$ obtained this way is defined exactly through the marking construction in \cite[Theorem 3.5]{SSS14}.

On the other hand, the convergence of $(S_{\eps}W^1,S_{\eps}W^2)$ to a limit $(\Wi^1, \Wi^2)$ whose law satisfies the conditions in Theorem \ref{T:martc} is a special case of our Theorem \ref{T:stick2}. For the special nearest neighbour case, it also follows from results in the literature. Conditions (i) and (iii) in Theorem \ref{T:martc} follow by the same argument as in the proof of Theorem \ref{T:stick2}, while the key condition (ii)
on the convergence of a pair of nearest neighbour sticky random walks to sticky Brownian motions was proved in \cite[Proposition A.6]{SSS14}
(the continuous time version was proved earlier in \cite[Theorem 8.1]{HW09a}. This concludes the proof sketch.
\epro
\medskip

%\noindent
\bpro[sketch for Theorem \ref{T:Nhop}] Assume $\theta=1$ to simplify the notation. In \cite{SS08}, the Brownian net $\Ni$ is defined via hopping
between a pair of sticky Brownian webs $(\Wi^{\rm l}, \Wi^{\rm r})$ with drifts $-1$ and $+1$ respectively, called the {\em left-right web}, where
paths in $\Wi^{\rm l}$ cannot cross paths in $\Wi^{\rm r}$ from left to right and vice versa. To show that $\Ni$ constructed
via hopping from a pair of sticky Brownian webs $(\Wi^1, \Wi^2)$ with parameter $\theta=1$ is the same as the $\Ni$ constructed via hopping from a
pair of left-right webs $(\Wi^{\rm l}, \Wi^{\rm r})$, it suffices to construct a coupling $(\Wi^1, \Wi^2, \Wi^{\rm l}, \Wi^{\rm r}, \Ni)$ that makes this equivalence obvious. In the discrete setting, such a coupling is given precisely by $(W^1, W^2, N_\eps)$ defined around \eqref{a2}, and then define $W^{\rm l}_\eps$ (resp.\ $W^{\rm r}_\eps$) by always choosing the leftmost (resp.\ rightmost) arrow coming out of each point in $\Z_{\rm even}^2:=\{(x, t)\in \Z^2: x+t \mbox{ is even}\}$. Equivalently, we can start with $N_\eps$ and then sample the arrows followed by $(W^1, W^2, W^{\rm l}_\eps, W^{\rm r}_\eps)$ at each branching point of $N_\eps$, such that $W^{\rm l}$ (resp.\ $W^{\rm r}$) always samples the leftmost (resp.\ rightmost) arrow,
$W^1$ samples one of the two outgoing arrows with equal probability while $W^2$ samples the other arrow. The scaling limit then gives the desired coupling between $(\Wi^1, \Wi^2, \Wi^{\rm l}, \Wi^{\rm r}, \Ni)$, where the same sampling distribution is applied independently at each relevant separation point of $\Ni$ to recover the quadruple $(\Wi^1, \Wi^2, \Wi^{\rm l}, \Wi^{\rm r})$. Such a convergence follows from the proof of Theorem 6.15 and (7.23) in \cite{SSS14}, and the coupling $(\Wi^1, \Wi^2, \Ni)$ was formulated in \cite[Lemma 6.16]{SSS14}, which can be easily extended to include also $(\Wi^{\rm l}, \Wi^{\rm r})$. Lastly, to see why hopping
among $\Wi^1 \cup \Wi^2$ and hopping among $\Wi^{\rm l} \cup \Wi^{\rm r}$ will both lead to $\Ni$, we note that this is because in both cases, the hopping construction will exhaust both branches of paths coming out of each relevant separation point of $\Ni$, just as in $N_\eps$, paths in
both $W^1\cup W^2$ and $W^{\rm l}\cup W^{\rm r}$ use up all arrows coming out of each branching point of $N_\eps$. Using the finite graph representation of paths in the Brownian net (see \cite[Section 6.2]{SSS14}), which is the continuum analogue of the directed graph connecting relevant branching points of $N_\eps$ as introduced in Section \ref{S:RBP}, every path $\pi \in \Ni$ can be approximated by a hopping path $\pi_\delta$ constructed from either $\Wi^1 \cup \Wi^2$ or $\Wi^{\rm l} \cup \Wi^{\rm r}$ that visits the same set of $(\sigma_\pi + (n-1)\delta, \sigma_\pi+n \delta)$-relevant separation points as $\pi$ for every $n\in\N$. In fact, it suffices to hop between paths at the relevant separation points of $\Ni$. As the time spacing $\delta\downarrow 0$, the approximating path $\pi_\delta$ then converges to $\pi$ thanks to the almost sure compactness of $\Ni$.
\epro
\end{appendix}

%%%%%%%%%%%%%%%%%%%%%%%%%%%%%%%%%%%%%%%%%%%%%%
%% Multiple Appendixes:                     %%
%%%%%%%%%%%%%%%%%%%%%%%%%%%%%%%%%%%%%%%%%%%%%%
%\begin{appendix}
%\section{???}
%
%\section{???}
%
%\end{appendix}

%%%%%%%%%%%%%%%%%%%%%%%%%%%%%%%%%%%%%%%%%%%%%%
%% Support information, if any,             %%
%% should be provided in the                %%
%% Acknowledgements section.                %%
%%%%%%%%%%%%%%%%%%%%%%%%%%%%%%%%%%%%%%%%%%%%%%
\begin{acks}[Acknowledgments]
We thank the referees for their helpful comments to improve the exposition of the paper.
J.~Yu acknowledges the support of NYU-ECNU Institute of Mathematical Sciences at NYU Shanghai.
\end{acks}
%%%%%%%%%%%%%%%%%%%%%%%%%%%%%%%%%%%%%%%%%%%%%%
%% Funding information, if any,             %%
%% should be provided in the                %%
%% funding section.                         %%
%%%%%%%%%%%%%%%%%%%%%%%%%%%%%%%%%%%%%%%%%%%%%%
\begin{funding}
 R.~Sun is supported by NUS Tier 1 grants A-8001448-00-00 and A-8004725-00-00, and NSFC grant 12271475.
 J.M.~Swart is supported by grant 20-08468S of the Czech Science Foundation (GA CR).
 J.~Yu is supported by National Key R\&D Program of China (No.~2021YFA1002700) and NSFC (No.~12101238 and 12271010).
\end{funding}

\bibliographystyle{imsart-number}
\bibliography{bibliography}

\end{document}